\documentclass[12pt]{article}
\usepackage{amssymb}
\usepackage{times}
\usepackage{epsfig}
\usepackage{amssymb,amsmath,amsthm}
\usepackage{amsfonts}%,amscd}

\newtheorem{defn}{Definition}[section]

\newtheorem{lem}[defn]{Lemma}

\newtheorem{prop}[defn]{Proposition}

\newtheorem{thm}[defn]{Theorem}

\newtheorem*{thmN}{Theorem 6.2}   % theorem NO nonumber, named Theorem 6.2

\newcommand{\be}{\begin{equation}}
\newcommand{\ee}{\end{equation}}
\newcommand{\bea}{\begin{eqnarray}}
\newcommand{\eea}{\end{eqnarray}}
\newcommand{\beas}{\begin{eqnarray*}}
\newcommand{\eeas}{\end{eqnarray*}}

\newcommand{\goto}{\rightarrow}

\newcommand{\ds}{\displaystyle}
\newcommand{\ts}{\textstyle}

\newcommand{\noi}{\noindent}

\newcommand{\ve}{\varepsilon}

\newcommand{\skp}{\vspace{\baselineskip}}

\newcommand{\R}{{\mathbb R}}

\newcommand{\N}{\mathbb N}

\newcommand{\barx}{\bar{x}}

\newcommand{\bn}{\beta_n}

\newcommand{\kn}{K_n}

\newcommand{\half}{\ts\frac{1}{2}}

\newcommand{\bc}{\beta_c}
\newcommand{\kc}{K(\bc)}

\newcounter{bean}
\newcommand{\benuma}{\setlength{\labelwidth}{.25in}
\begin{list}%
{(\alph{bean})}{\usecounter{bean}}}
\newcommand{\eenuma}{\end{list}}

\def\theequation{\thesection.\arabic{equation}}
\def\theequation{\arabic{section}.\arabic{equation}}
\def\thedefn{\arabic{section}.\arabic{defn}}
\newcommand{\beginsec}{\setcounter{equation}{0}}

\begin{document}

% \pagewiselinenumbers
% use with "\usepackage{lineno}"

\title{Refined Asymptotics of the Finite-Size \\
Magnetization via a New Conditional \\
Limit Theorem for the Spin}
\author{Richard S.\ Ellis\normalsize \vspace{-.1in} \\
\small{rsellis@math.umass.edu} \vspace{-.125in} \\ \\
Jingran Li \normalsize \vspace{-.1in}\\
\small{jingran@math.umass.edu} \vspace{-.125in} \\ \\
\normalsize{Department of Mathematics and Statistics} \vspace{-.05in} \\
\normalsize{University of Massachusetts} \vspace{-.05in} \\
\normalsize{Amherst, MA 01003} \vspace{-.1in}}
\maketitle

\begin{abstract}
We study the fluctuations of
the spin per site around the thermodynamic magnetization in the mean-field Blume-Capel model.
Our main theorem generalizes the main result in a previous paper \cite{EllMacOtt2010}
in which the first rigorous confirmation of the
statistical mechanical theory of finite-size scaling for a mean-field model is given. In that paper
our goal is to determine whether the thermodynamic magnetization is
a physically relevant estimator of the finite-size magnetization. This is done by comparing the asymptotic
behaviors of these two quantities along parameter
sequences converging to either a second-order point or the tricritical point in the mean-field
Blume-Capel model. The main result is that the thermodynamic magnetization and the finite-size magnetization
are asymptotic when the parameter $\alpha$ governing the speed at which
the sequence approaches criticality is below a certain threshold $\alpha_0$. Our main theorem in the present
paper on the fluctuations of the spin per site around the thermodynamic magnetization is based on a new
conditional limit theorem for the spin, which is closely related to a new conditional central
limit theorem for the spin.
\end{abstract}

\noi
{\it American Mathematical Society 2000 Subject Classifications.}  Primary 60F05, 60F10, Secondary 82B20
\skp

\noi

%===================== new section =====================================
\section{Introduction}
\label{section:Intro}
\beginsec

The purpose of this paper is to analyze the asymptotic behavior of
the fluctuations of the spin per site around the thermodynamic
magnetization along parameter sequences having physical relevance in
the mean-field Blume-Capel model. This research culminates a series
of papers that study the phase-transition structure of the model via
analytic techniques and probabilistic limit theorems
\cite{EllOttTou, CosEllOtt, EllMacOtt1, EllMacOtt2010}. \iffalse As
a mean-field version of an important lattice spin model, the
mean-field Blume-Capel model was due to Blume and Capel, to which we
refer as the B-C model \cite[6-8]{Blu}. We refer to this model as
the mean-field B-C model \cite[6-8]{Blu}. On the complete graph on
$n$ vertices, this mean-field model is equivalent to the B-C model.
\fi The mean-field Blume-Capel model is a mean-field version of an
important lattice model due to Blume and Capel, to which we will
refer as the B--C model \cite{Blu,Cap1,Cap2,Cap3}. The mean-field
B-C model is an important object of study because it is one of the
simplest models that exhibits the following complicated
phase-transition structure: a curve of second-order points; a curve
of first-order points; and a tricritical point, which separates the
two curves.

\iffalse (this is deleted on Oct. 12.) The main finding in this
dissertation is that the speed of the fluctuations of the $|S_n/n|$
around the thermodynamic magnetization converge to 0 is governed by
the parameter $\kappa$, which is related to a parameter $\alpha$ in
\cite{EllMacOtt2010}. Specifically, our main result Theorem
\ref{thm:main} gives that for $0<\alpha<\alpha_0$ this fluctuation
converges to 0 at rate $\bar{z}/n^{\kappa}$, where $\bar{z}$ is a
certain explicit number calculated by the value of the second
derivative of Ginzburg-Landau polynomials at their global minimum
points. This main finding is a consequence of the other main result,
Theorem \ref{thm:6.1}, in this dissertation, proving a non-classical
conditioned central limit theorem. \fi

The main theorem in this paper generalizes the main result in
\cite{EllMacOtt2010}. The goal of \cite{EllMacOtt2010} is to compare
the asymptotic behaviors of the thermodynamic magnetization and the
finite-size magnetization along parameter sequences converging to
either a second-order point or the tricritical point of the
mean-field B-C model. Theorem 4.1 in that paper shows that these two
quantities are asymptotic when the parameter $\alpha$ governing the
speed at which the sequence approaches criticality is below a
certain threshold $\alpha_0$. However, when $\alpha$ exceeds
$\alpha_0$, the thermodynamic magnetization converges to 0 much
faster than the finite-size magnetization. These results in
\cite{EllMacOtt2010} are worthwhile because they are the first
rigorous confirmations of the statistical mechanical theory of
finite-size scaling for a mean-field model \cite{Barber}, \cite[\S
6]{EllMacOtt2010}.

The importance of both the theory of finite-size scaling and the
mean-field B-C model motivate us in this paper to refine Theorem 4.1
in \cite{EllMacOtt2010}. We do this by studying the fluctuations of
the spin per site around the thermodynamic magnetization for $0 <
\alpha < \alpha_0$, obtaining a more refined asymptotic estimate
that yields the conclusion of Theorem 4.1 in \cite{EllMacOtt2010} as
a corollary. This refined asymptotic estimate is stated in
(\ref{eqn:refined}) and is proved in part (a) of Theorem \ref{thm:main} below.
While Theorem 4.1 in \cite{EllMacOtt2010} is obtained from a
moderate deviation principle, the refinement of that theorem in this
paper is obtained from the conditional limit theorem stated in (\ref{eqn:condlimit})
and proved in part (b) of Theorem \ref{thm:6.1}.

The mean-field B-C model is defined by a canonical ensemble that we
denote by $P_{N,\beta,K}$; $N$ is the number of vertices, $\beta>0$
is the inverse temperature, and $K>0$ is the interaction strength.
$P_{N,\beta,K}$ is defined in (\ref{eqn:P}) in terms of the
Hamiltonian
\[ H_{N,K}(\omega) = \sum_{j=1}^N \omega_j^2 - \frac{K}{N} \left( \sum_{j=1}^N
\omega_j \right)^2.
\]
In this formula $\omega_j$ is the spin at site $j \in
\{1,2,\ldots,N\}$ and takes values in $\Lambda = \{-1,0,1\}$. The
configuration space for the model is the set $\Lambda^N$ containing
all sequences $\omega=(\omega_1,\ldots,\omega_N)$ with each
$\omega_j \in \Lambda$. Expectation with respect to $P_{N,\beta,K}$
is denoted by $E_{N,\beta,K}$. The finite-size magnetization is
defined by $E_{N,\beta,K}\{|S_N/N|\}$, where $S_N$ equals the total
spin $\sum_{j=1}^N \omega_j$.

Before discussing the results in this paper, we first summarize the
phase-transition structure of the mean-field B-C model as derived in
\cite{EllOttTou}. For $\beta >0$ and $K>0$, we denote by
$\mathcal{M}_{\beta,K}$ the set of equilibrium values of the
magnetization. The set $\mathcal{M}_{\beta,K}$ coincides with the
set of global minimum points of the free-energy function
$G_{\beta,K}$, which is defined in (\ref{eqn:Cbeta})--(\ref{eqn:GbetaK}). The critical
inverse temperature of the mean-field B-C model is $\beta_c=\log 4$.
For $0< \beta \le \beta_c$ there exists a quantity $K(\beta)$ and
for $\beta > \beta_c$ there exists a quantity $K_1(\beta)$ having
the following properties. The positive quantity $m(\beta,K)$
appearing in this list is the thermodynamic magnetization.
\begin{enumerate}
\item Fix $0 < \beta \leq \bc$. Then for $0 < K \leq K(\beta)$,
$\mathcal{M}_{\beta,K}$ consists of the unique pure phase 0, and for
$K > K(\beta)$, $\mathcal{M}_{\beta,K}$ consists of two nonzero
values $\pm m(\beta,K)$.
\item For $0 < \beta \leq \beta_c$, $\mathcal{M}_{\beta,K}$
undergoes a continuous bifurcation at $K = K(\beta)$, changing
continuously from $\{0\}$ for $K \leq K(\beta)$ to $\{\pm
m(\beta,K)\}$ for $K > K(\beta)$. This continuous bifurcation
corresponds to a second-order phase transition.
\item Fix $\beta > \bc$. Then for $0 < K < K_1(\beta)$, $\mathcal{M}_{\beta,K}$ consists of the unique
pure phase 0; for $K = K_1(\beta)$, $\mathcal{M}_{\beta,K}$ consists
of 0 and two nonzero values $\pm m(\beta,K_1(\beta))$; and for $K >
K_1(\beta)$, $\mathcal{M}_{\beta,K}$ consists of two nonzero values
$\pm m(\beta,K)$.
\item For $\beta > \beta_c$, $\mathcal{M}_{\beta,K}$
undergoes a discontinuous bifurcation at $K = K_1(\beta)$, changing
discontinuously from $\{0\}$ for $K < K(\beta)$ to $\{0, \pm
m(\beta,K)\}$ for $K = K_1(\beta)$ to $\{\pm m(\beta,K)\}$ for $K >
K_1(\beta)$. This discontinuous bifurcation corresponds to a
first-order phase transition.
\end{enumerate}

Because of item 2, we refer to the curve
$\{(\beta,K(\beta)),0<\beta<\beta_c\}$ as the second-order curve
and points on this curve as second-order points.
Because of item 4, we refer to the curve
$\{(\beta,K_1(\beta)),\beta>\beta_c\}$ as the first-order curve
and points on this curve as first-order points. The
point $(\beta_c, K(\beta_c))=(\log 4, 3/2\log4)$, called the
tricritical point, separates the second-order curve from the
first-order curve. The phase-coexistence region is defined as the
set of all points in the positive $\beta$-$K$ quadrant for which
${\mathcal{M}}_{\beta,K}$ consists of more than one value. Therefore
the phase-coexistence region consists of all points above the
second-order curve, above the tricritical point, on the first-order
curve, and above the first-order curve; that is,
\[ \{(\beta,K): 0 < \beta \le \beta_c, K > K(\beta) \ \textrm{and} \ \beta > \beta_c, K \ge K_1(\beta)
\}.
\]
Figure 1 exhibits the sets that describe the phase-transition structure of mean-field B-C
model.

% RSE deleted [phase] in \begin{figure}[phase] and added [h].
\begin{figure}[h]
\begin{center}
\epsfig{file=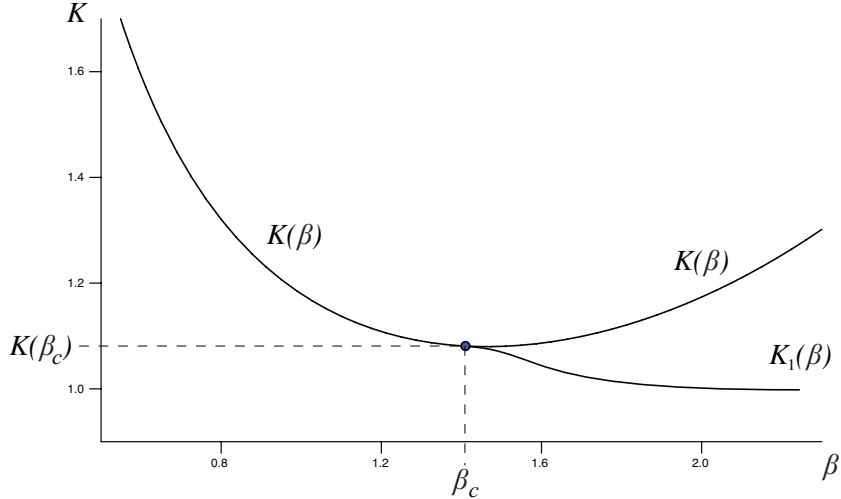,width=12cm} \caption{\footnotesize
The sets that describe the phase-transition structure of the
mean-field B-C model: the second-order curve $\{(\beta,K(\beta)), 0
< \beta < \bc\}$, the first-order curve $\{(\beta,K_1(\beta)), \beta
> \bc\}$, and the tricritical point $(\bc,\kc)$.  The
phase-coexistence region consists of all $(\beta,K)$ above the
second-order curve, above the tricritical point, on the first-order
curve, and above the first-order curve. \iffalse RSE deleted the
last sentence in the caption. The extension of the second-order
curve to $\beta > \bc$ is called the spinodal curve. \fi}
\end{center}
\end{figure}

In order to discuss the contributions of this paper, it is helpful
first to explain the main results in \cite{EllMacOtt1} and
\cite{EllMacOtt2010}. Those papers focus on positive sequences
$(\beta_n,K_n)$ that lie in the phase-coexistence region for all
sufficiently large $n$, converge to either a second-order point or
the tricritical point, and satisfy the hypotheses of Theorem 3.2 in
\cite{EllMacOtt1}. These sequences are parameterized by $\alpha >0$
in the sense that the limits
\[
b = \lim_{n \goto \infty} n^\alpha (\bn - \beta) \ \mbox{ and } \ k
= \lim_{n \goto \infty} n^\alpha (\kn - K(\beta))
\]
are assumed to exist and are not both 0. Six specific such sequences
are introduced in section 4 of that paper. Theorem 3.2 in
\cite{EllMacOtt1} states that for any $\alpha >0$, $m(\beta_n,K_n)$
has the asymptotic behavior \be \label{eqn:mbnkn} m(\beta_n,K_n)
\sim {\bar{x}}/{n^{\theta\alpha}}, \ee where $\theta >0$ and
$\bar{x}$ is the positive global minimum point of a certain
polynomial $g(x)$ called the Ginzburg-Landau polynomial.  This
polynomial is defined in terms of the free-energy function
$G_{\beta,K}$ in hypothesis (iii)(a) of Theorem \ref{thm:3.1} below.

One of the surprises in our study of the mean-field B-C model is the
appearance of the Ginzburg-Landau polynomial in a number of basic
results. These include the asymptotic formula (\ref{eqn:mbnkn}),
the quantity $\bar{z}$ in the asymptotic formula (\ref{eqn:refined})
and the conditional limit theorem (\ref{eqn:condlimit}),
the limiting variance in the conditional central
limit theorem (\ref{eqn:condclt}), and the rate function in
the moderate deviation principle in Theorem \ref{thm:MDP}. As we will explain,
this conditional central limit theorem is closely related to the main result in
this paper, which is the asymptotic formula (\ref{eqn:refined}).

A straightforward large-deviation calculation summarized in \cite[p.\
2120]{EllMacOtt2010} shows that for fixed $(\beta,K)$ lying in the
phase-coexistence region the spin per site $S_N/N$ has the weak-convergence
limit
\be
\label{eqn:spin}
P_{N, \beta, K} \{S_N/N \in dx\}
\Longrightarrow
\left(\frac{1}{2}\delta_{m(\beta,K)}+\frac{1}{2}\delta_{-m(\beta,
K)}\right)\!(dx).
\ee
This implies that
\[ \lim_{N \goto \infty} E_{N, \beta, K} \{|S_N/N|\}=m(\beta,K).
\]
Because the thermodynamic magnetization $m(\beta, K)$ is the limit,
as the number of spins goes to $\infty$, of the finite-size
magnetization $E_{N, \beta, K} \{|S_N/N|\}$, the thermodynamic
magnetization is a physically relevant estimator of the
finite-size magnetization, at least when evaluated at fixed
$(\beta,K)$ in the phase-coexistence region.

The main focus of \cite{EllMacOtt2010} is to determine whether the
thermodynamic magnetization is a physically relevant estimator of
the finite-size magnetization in a more general sense, namely, when
evaluated along positive sequences that lie in the phase-coexistence
region for all sufficiently large $n$, converge to a second-order
point or the tricritical point, and satisfy a set of hypotheses
including those of Theorem 3.2 in \cite{EllMacOtt1}. The criterion
for determining whether $m(\bn,\kn)$ is a physically relevant
estimator is that as $n \goto \infty$, $m(\bn,\kn)$ is asymptotic to
the finite-size magnetization $E_{n,\bn,\kn}\{|S_n/n|\}$, both of
which converge to 0. In this formulation we let $N = n$ in the
finite-size magnetization; i.e., we let the number of spins $N$
coincide with the index $n$ parametrizing the sequence $(\bn,\kn)$.

As summarized in Theorems 4.1 and 4.2 in \cite{EllMacOtt2010}, the
main finding is that $m(\bn,\kn)$ is a physically relevant estimator
when the parameter $\alpha$ governing the speed at which
$(\bn,\kn)$ approaches criticality is below a certain threshold
$\alpha_0$. The value of $\alpha_0$ depends on the type of the phase
transition --- first-order, second-order, or tricritical --- that
influences the sequence, an issue addressed in section 6 of
\cite{EllMacOtt1}. For $0 < \alpha < \alpha_0$ this finding is
summarized by the limit \be \label{eqn:enbnknmbnkn} \lim_{n \goto
\infty} n^{\theta\alpha} \left|E_{n, \beta_n, K_n}\{|S_n/n|\} -
m(\beta_n, K_n) \right| = 0, \ee which in combination with
(\ref{eqn:mbnkn}) implies that
\[
E_{n, \beta_n, K_n}\{|S_n/n|\} \sim \bar{x}/n^{\theta\alpha} \sim
m(\bn,\kn).
\]
By contrast, when $\alpha > \alpha_0$, $m(\bn,\kn)$ converges to 0
much faster than $E_{n, \beta_n, K_n}\{|S_n/n|\}$. The sequences for
which these asymptotic results are valid include the six sequences
introduced in \cite[\S 4]{EllMacOtt1}.

\iffalse Theorem 4.1 in \cite{EllMacOtt2010} shows that when the
parameter $\alpha$ governing the speed at which $(\beta_n, K_n)$
approaches criticality is below a certain threshold $\alpha_0$, the
thermodynamic magnetization $m(\beta_n,K_n)$ is a physically
relevant estimator of the finite-size magnetization
$E_{n,\beta_n,K_n}\{|S_n/n|\}$, namely,
\[ E_{n,\beta_n,K_n}\{|S_n/n|\} \sim m(\beta_n,K_n) \sim
{\bar{x}}/{n^{\theta\alpha}}.
\]
On the other hand, when $\alpha > \alpha_0$, the thermodynamic
magnetization converges to 0 asymptotically faster than the
finite-size magnetization, namely,
\[E_{n,\beta_n,K_n}\{|S_n/n|\} \sim {\bar{y}}/{n^{\theta\alpha_0}}  \ \ \ \textrm{and} \
\]
\[E_{n,\beta_n,K_n}\{|S_n/n|\} \gg m(\beta_n,K_n) \sim
{\bar{x}}/{n^{\theta\alpha}},
\]
where $\bar{y}$ is a positive quantity. In this dissertation, we
focus on the case $0<\alpha<\alpha_0$. \fi

We now turn to the main focus of this paper, which is a refined
analysis of the fluctuations of $S_n/n$ around $m(\beta_n,K_n)$ for
$0<\alpha<\alpha_0$. \iffalse This refined analysis is based on the
conditional central limit theorem stated in part (a) of Theorem
\ref{thm:6.1}. The central limit theorem leads to the following
improvement of the asymptotic formula (\ref{eqn:enbnknmbnkn}). \fi
Define $\kappa = \frac{1}{2} (1-\alpha/\alpha_0)+\theta\alpha$.
As shown in part (a) of Theorem \ref{thm:main}, for $0 < \alpha <
\alpha_0$ and for a class of sequences that includes the first five
sequences introduced in \cite[\S 4]{EllMacOtt1}
\be
\label{eqn:refined}
E_{n,\beta_n,K_n}\{||S_n/n| - m(\beta_n,K_n)|\}
\sim \bar{z}/n^{\kappa}.
\ee
In this formula $\bar{z} = (2 / [\pi g^{(2)}(\bar{x})] )^{1/2}$, where
$g^{(2)}(\bar{x})$ denotes the positive second derivative of the
Ginzburg-Landau polynomial $g$ evaluated at its unique positive
global minimum point $\bar{x}$. For all $0<\alpha<\alpha_0$,
$\kappa$ is larger than $\theta\alpha$. Thus the rate
$\bar{z}/n^\kappa$ at which $E_{n, \beta_n,K_n}\{ | |S_n/n| -
m(\beta_n, K_n)  | \}$ converges to 0 is asymptotically faster than
the rate $\bar{x}/n^{{\theta\alpha}}$ at which $E_{n, \beta_n,K_n}\{
|S_n/n| \}$ and $m(\beta_n, K_n)$ converge separately to 0.

This asymptotic result generalizes (\ref{eqn:enbnknmbnkn}), which is the
conclusion of Theorem 4.1 in \cite{EllMacOtt2010}.  To see this,
define $A_n = E_{n,\beta_n,K_n}\{||S_n/n| - m(\beta_n,K_n)|\}$ and
note that
\[  |E_{n,\beta_n,K_n} \{ |S_n/n|\} - m(\beta_n,K_n)| \le  A_n.
\]
Equation (\ref{eqn:refined}) states that $\lim_{n \goto \infty} n^\kappa A_n = \bar{z}$.
Since $\kappa > \theta\alpha$, this implies that
\[
0 = \lim_{n \goto \infty} n^{\theta\alpha} A_n \geq \lim_{n \goto
\infty} n^{\theta\alpha}  |E_{n,\beta_n,K_n} \{ |S_n/n|\} -
m(\beta_n,K_n)| = 0.
\]
The fact that this second limit equals 0 yields
(\ref{eqn:enbnknmbnkn}), which is the conclusion of Theorem 4.1 in
\cite{EllMacOtt2010}.

The proof of our main result (\ref{eqn:refined}) is based on the
following new
conditional limit stated in part (b) of Theorem \ref{thm:6.1} for $0
< \alpha < \alpha_0$: 
\be 
\label{eqn:condlimit} \lim_{n \goto \infty}
n^\kappa E_{n,\beta_n,K_n} \{ |S_n/n-m(\beta_n,K_n)|  \ \big| \
S_n/n > \delta m(\beta_n,K_n)\} = \bar{z}. 
\ee
The conditioning is on the event
$\{S_n/n>\delta m(\beta_n,K_n)\}$, where $\delta \in (0,1)$ is
sufficiently close to 1. This conditioning allows us to study the
asymptotic behavior of the system in a neighborhood of the pure
states having thermodynamic magnetization $m(\beta_n, K_n)$. According
to Lemma \ref{lem:Sym1/2Small}
\bea
\label{eqn:half} \lefteqn{
\lim_{n \goto \infty} P_{n,\bn,\kn}\{S_n/n > \delta m(\beta_n,K_n)\} } \\
\nonumber && = \lim_{n \goto \infty} P_{n,\bn,\kn}\{S_n/n < -\delta
m(\beta_n,K_n)\} = 1/2.
\eea
\iffalse
----------------------------------
which is under the same hypotheses of the
main result Theorem \ref{thm:main}. This result involves a
conditioned event $S_n/n>\delta m(\beta_n,K_n)$, for $\delta \in
(0,1)$ sufficiently close to 1, which means that $S_n/n$ is
conditioned to lie in a suitable neighborhood of $m(\beta_n,K_n)$.
--------------------\fi
These limits are the analog of the weak
convergence limit (\ref{eqn:spin}), showing that as $n \goto \infty$
the mass of the $P_{n,\beta_n,K_n}$-distribution of $S_n/n$
concentrates at $\pm m(\beta_n,K_n)$. 
\iffalse
where
\be
\label{eqn:barz}
\bar{z} =
E\{|N(0,(1/g^{(2)}(\bar{x}))^{1/2})|\} =
(2 / [\pi g^{(2)}(\bar{x})] )^{1/2}.
\ee
In (\ref{eqn:barz})
$g^{(2)}(\bar{x})$ denotes the positive second derivative of the
Ginzburg-Landau polynomial $g$ evaluated at its unique positive
global minimum point $\bar{x}$.
\fi
As we show in section 6, the limits (\ref{eqn:half}) and (\ref{eqn:condlimit})
and a moderate deviation estimate on the
probability $P_{n,\bn,\kn}\{\delta m(\bn,\kn) \geq S_n/n \geq
-\delta m(\bn,\kn)\}$ yield
\[
\lim_{n \goto \infty} n^\kappa E_{n,\beta_n,K_n}\{||S_n/n| -
m(\beta_n,K_n)|\} = \bar{z}.
\]
This limit is equivalent to (\ref{eqn:refined}).

The main result in part (a) of Theorem \ref{thm:main}  is applied to
the first five sequences introduced in \cite[\S 4]{EllMacOtt1}.
Located in the phase-coexistence region for all sufficiently large
$n$, the first two sequences converge to a second-order point, and
the last three sequences converge to the tricritical point. Possible
paths followed by these sequences are shown in Figure 2. For each of
the five sequences the quantities $\alpha_0$, $\theta$, and $\kappa$
appearing in Theorem \ref{thm:main} are specified in Table 1.1.

% RSE DELETED [Sequence] AND ADDED [h]
\begin{figure}[h]
\begin{center}
\epsfig{file=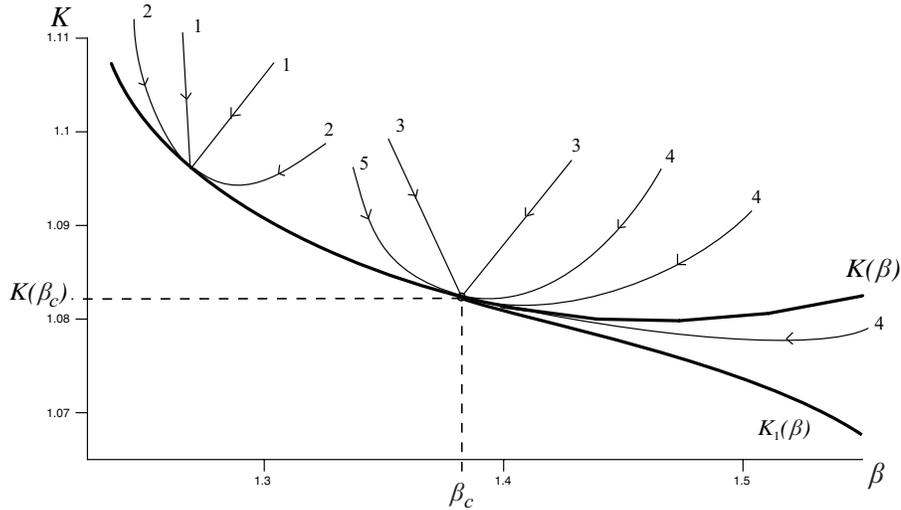,width=12cm}
\vspace{-.75in}
 \caption{\footnotesize Possible paths for the five
sequences converging to a second-order point and to the tricritical
point. In section \ref{section:verify6} 
and appendix A these sequences are defined
and are shown to satisfy the hypotheses of Theorem \ref{thm:main} and
Theorem \ref{thm:6.1}. The sequences labeled 1--5 in this figure
correspond to sequences 1a--5a in Table 1.1 and Table 5.1. }
\end{center}
\end{figure}

\iffalse
\begin{center}
% \label{table: SequenceInformation}
\begin{tabular}{lllll} \hline
{{\bf Seq.}} & {\bf Defn.} & \boldmath $\alpha_0$ \unboldmath &
\boldmath $\theta$ \unboldmath & \boldmath $\kappa$ \unboldmath \\
\hline 1 & (\ref{eqn:seq1}) & $\half$ & $\half$ & $\half
(1-\alpha)$ \\
2 & (\ref{eqn:seq2}) & $\frac{1}{2p}$ & $\frac{p}{2}$ &
$\half(1-p\alpha)$ \\ 3 & (\ref{eqn:seq3}) & $\frac{2}{3}$ &
$\frac{1}{4}$ & $\half(1-\alpha)$ \\ 4 & (\ref{eqn:seq4}) &
$\frac{1}{3}$ & $\half$ & $\half(1 - 2\alpha)$
\\ 5 & (\ref{eqn:seq5}) & $\frac{1}{3}$ & $\half$ & $\half(1- 2\alpha)$ \\ \hline
\end{tabular}
\end{center}
% \vspace{-.025in}
Table 1.1: {\small The equations where each of the five sequences is
defined and the values of $\alpha_0$, $\theta$, and $\kappa$ for
each sequence.} \fi

\begin{table}[h]
\begin{center}
% \label{table:13cases}
\begin{tabular}{||p{2cm}|p{2cm}|p{2cm}|p{2cm}|p{3cm}||}  \hline \hline
{{\bf Seq.}} & {\bf Defn.} & \boldmath $\alpha_0$ \unboldmath &
\boldmath $\theta$ \unboldmath & \boldmath $\kappa$ \unboldmath \\
\hline \hline 1a & (\ref{eqn:seq1}) & $\half$ & $\half$ & $\half
(1-\alpha)$ \\ \hline 2a & (\ref{eqn:seq2}) & $\frac{1}{2p}$ &
$\frac{p}{2}$ & $\half(1-p\alpha)$ \\ \hline 3a & (\ref{eqn:seq3}) &
$\frac{2}{3}$ & $\frac{1}{4}$ & $\half(1-\alpha)$ \\ \hline 4a &
(\ref{eqn:seq4}) & $\frac{1}{3}$ & $\half$ & $\half(1 - 2\alpha)$
\\ \hline
5a & (\ref{eqn:seq5}) &
$\frac{1}{3}$ & $\half$ & $\half(1- 2\alpha)$ \\
\hline \hline
\end{tabular}
\end{center}
% \vspace{-.25in}
Table 1.1: {\small The equations where each of the five sequences is
defined and the values of $\alpha_0$, $\theta$, and $\kappa$ for
each sequence.}
\end{table}

 The conditional limit (\ref{eqn:condlimit}) is closely related to
another result stated in part (a) of Theorem \ref{thm:6.1}. This
result is a new conditional central limit theorem for
$0<\alpha<\alpha_0$. As in (\ref{eqn:condlimit}), the conditioning is
on the event $\{S_n/n>\delta m(\beta_n,K_n)\}$, where $\delta \in
(0,1)$ is sufficiently close to 1. Under a set of hypotheses
satisfied by the first five sequences introduced in \cite[\S
4]{EllMacOtt1}, part (a) of Theorem \ref{thm:6.1} states that when
conditioned on $\{ S_n/n > \delta m(\beta_n, K_n)\}$, the $P_{n,
\beta_n,K_n}$-distributions of $n^\kappa(S_n/n-m(\beta_n,K_n))$
converge weakly to a normal random variable
$N(0,1/g^{(2)}(\bar{x}))$ with mean $0$ and variance
$1/g^{(2)}(\bar{x})$; in symbols,
\bea
\label{eqn:condclt}
\lefteqn{ \hspace{-1in}
P_{n,\beta_n,K_n}\{n^\kappa(S_n/n-m(\beta_n,K_n)) \in dx \  \big| \ S_n/n>\delta m(\beta_n,K_n)\} } \\
\hspace{-1in} && \nonumber \Longrightarrow
N(0,1/g^{(2)}(\bar{x})). \eea
Since $\kappa = \frac{1}{2}(1-\alpha/\alpha_0) + \theta\alpha$ is less than 1/2 [Thm.\
\ref{thm:6.1}(c)], the scaling in this result is non-classical.
An equivalent formulation is
that for any bounded, continuous function $f$
\bea
\label{eqn:condcltEquiv}
\lefteqn{ \lim_{n \goto \infty}
E_{n,\beta_n,K_n} \{ f(n^\kappa(S_n/n-m(\beta_n,K_n))) \ \big| \
S_n/n>\delta
m(\beta_n, K_n) \} } \\
&&=
\lim_{n \goto \infty}  E_{n,\beta_n,K_n} \{ f(S_n/n^{1-\kappa}-
n^{\kappa}m(\beta_n,K_n)) \ \big| \  S_n/n>\delta m(\beta_n, K_n)
\} \nonumber \\
&&=
E \{ f(N(0,1/g^{(2)}(\bar{x}))) \} \nonumber \\
&&=
\frac{1}{\int_{\mathbb{R}}\exp[-\frac{1}{2}
g^{(2)}(\bar{x})x^2]dx } \cdot \int_{\mathbb{R}} f(x) \exp[ \ts
-\frac{1}{2} g^{(2)}(\bar{x})x^2]dx.\nonumber \eea

Through the term $g^{(2)}(\bar{x})$ this conditional central limit theorem
and the asymptotic formula (\ref{eqn:refined}) exhibit
a sensitive dependence on the choice of the sequence $(\beta_n,K_n)$, which lies
in the phase coexistence region for all sufficiently large $n$ and converges to a second-order
point or the tricritical point. This contrasts sharply with the central limit theorem
that is valid for an arbitrary sequence $(\beta_n,K_n)$ that converges to a point $(\beta,K)$ in
the single-phase region defined by $\{(\beta, K) : 0 < \beta \leq \beta_c, 0 < K < K(\beta)\}$.
In this situation it is proved in Theorem 5.5 in \cite{CosEllOtt} that
\[
P_{n,\beta_n,K_n}\{S_n/n^{1/2} \in dx\} \Longrightarrow N(0,\sigma^2(\beta,K)),
\]
where the limiting variance $\sigma^2(\beta,K)$ depends only on $(\beta, K)$ and not on the sequence
$(\beta_n,K_n)$.

Formally, the conditional limit (\ref{eqn:condlimit}) follows from
the conditional central limit theorem (\ref{eqn:condcltEquiv}) if
one replaces the bounded, continuous function $f$ by the
absolute value function. Then (\ref{eqn:condcltEquiv}) would imply
\beas
\iffalse \label{eqn:condlimitEquiv} \fi
\lefteqn{ \lim_{n \goto \infty}
n^{\kappa}E_{n,\beta_n,K_n} \{ |S_n/n-m(\beta_n,K_n)| \ \big| \
S_n/n>\delta m(\beta_n, K_n) \} } \\
&&=
\lim_{n \goto \infty}  E_{n,\beta_n,K_n} \{ |S_n/n^{1-\kappa}-
n^{\kappa}m(\beta_n,K_n)| \ \big| \  S_n/n>\delta m(\beta_n,K_n)
\} \nonumber \\
&&=
\frac{1}{\int_{\mathbb{R}}\exp[-\frac{1}{2}
g^{(2)}(\bar{x})x^2]dx } \cdot \int_{\mathbb{R}} |x| \exp[\ts
-\frac{1}{2} g^{(2)}(\bar{x})x^2]dx = \bar{z}.\nonumber \eeas

\iffalse In order to justify this, one needs a uniform integrability
estimate. Unfortunately, we are unable to prove this estimate.
Because of this, we prove the conditional limit
(\ref{eqn:condlimit}) directly using a series of lemmas in section
\ref{section:LogicANDPreLemmas}. A parallel argument that proves the
conditional central limit theorem (\ref{eqn:condcltEquiv}) is given
in section 7 of \cite{LiEll}. \fi

In order to justify this formal derivation, one needs a uniform
integrability estimate. In fact, we can derive the conditional limit
(\ref{eqn:condlimit}) from a related weak convergence result via a
more circuitous route. The related weak convergence result, proved
in Lemma \ref{lem:Thm61aStep1}, involves two extra summands defined
in terms of a sequence of scaled normal random variables $W_n$. We
prove the conditional limit (\ref{eqn:condlimit}) by two steps: the
uniform integrability-type result in Proposition
\ref{prop:weakinteg} allows us to replace the bounded,
continuous function $f$ in Lemma \ref{lem:Thm61aStep1} by the
absolute value function; the calculations in Lemmas
\ref{lem:substep2aOFThm61b} and \ref{lem:Substep2bOFThm61b} show
that in the limit $n \goto \infty$ the extra summands involving the
normal random variables $W_n$ do not affect the limit. 
As we show at the end of section 7 in \cite{EllLi}, we also use the weak convergence result in Lemma
\ref{lem:Thm61aStep1} to prove the conditional central limit theorem
(\ref{eqn:condclt}) by an analogous but more straightforward argument. Again, a key step is
to show that in the limit $n \goto \infty$ the extra summands
involving the normal random variables $W_n$ do not affect the limit.

The conditional limit (\ref{eqn:condlimit}) is stated in part (b) of
Theorem \ref{thm:6.1}, the proof of which
\iffalse
We would like to motivate this limit by
giving some insight into the value of $\kappa$ and the formula for
$\bar{z}$. However, a straightforward motivation seems to be
unavailable because the proof of the conditional limit in part (b) of Theorem \ref{thm:6.1}
\fi
is subtle and complicated. In this proof Lemma
\ref{lem:G(xmn)miusG(mn)} is key. There we obtain two basic
estimates that allow us to apply the Dominated Convergence Theorem to
prove the weak convergence result in Lemma \ref{lem:Thm61aStep1}, from
which part (b) of Theorem \ref{thm:6.1} will be deduced. The value of
$\kappa$ can be motivated from
the calculation underlying the proof of part (a) of Lemma
\ref{lem:G(xmn)miusG(mn)}.

The contents of this paper are as follows. In section
\ref{section:BEGmodel} we define the mean-field B-C model and
summarize its phase-transition structure in Theorems
\ref{thm:secondorder} and \ref{thm:firstorder}. For a class of
sequences $(\beta_n,K_n)$ lying in the phase-coexistence region for
all sufficiently large $n$ and converging either to a second-order
point or to the tricritical point, Theorem \ref{thm:3.1} in section
\ref{section:resultPrevious} describes the asymptotic behavior of
$m(\beta_n, K_n) \goto 0$ as stated in (\ref{eqn:mbnkn}).  Theorem
\ref{thm:mainThm4.1Of2010} in section \ref{section:resultPrevious}
states one of the main results of \cite{EllMacOtt2010}, which is
that as $n \goto 0$, $m(\beta_n,K_n)$ is asymptotic to
$E_{n,\beta_n,K_n}\{|S_n/n|\}$ for $0<\alpha<\alpha_0$, proving that
for this range of $\alpha$ the thermodynamic magnetization
$m(\bn,\kn)$ is a physically relevant estimator of the finite-size
magnetization $E_{n,\beta_n,K_n}\{|S_n/n|\}$.

The main result in this paper is given in section
\ref{section:resultmain}. According to part (a) of Theorem
\ref{thm:main}, for $0 < \alpha < \alpha_0$
\[
E_{n,\beta_n,K_n}\{|S_n/n-m(\beta_n,K_n)|\} \sim \bar{z}/n^{\kappa},
\]
where $\bar{z} = (2 / [\pi g^{(2)}(\bar{x})] )^{1/2}$ and $\kappa =
\frac{1}{2}(1 - \alpha/\alpha_0) + \theta\alpha$. Part (a) of
Theorem \ref{thm:main} is applied in section \ref{section:verify6}
to five specific sequences $(\beta_n,K_n)$. The first two sequences
converge to a second-order point, and the last three sequences
converge to the tricritical point. In section
\ref{section:proofmain} part (a) of Theorem \ref{thm:6.1} states the
conditional central limit theorem (\ref{eqn:condclt}), and part (b)
of that theorem states the conditional limit (\ref{eqn:condlimit}).
\iffalse ------------------------------------------------------ In
section \ref{section:proofmain} Theorem \ref{thm:6.1} gives that
when $S_n/n$ is conditioned to lie in a suitable neighborhood of
$m(\beta_n,K_n)$, the $P_{n,\beta_n,K_n}$ -distributions of
$n^{\kappa}(S_n/n-m(\beta_n,K_n))$ converge in distribution to a
Gaussian and the corresponding conditioned expectation limit is a
explicit constant. Lemma \ref{lem:Sym1/2Small} uses MDP in Theorem
\ref{thm:MDP} and the asymptotic behavior of $m(\beta_n,K_n)$ to
give the probabilities of the event when $S_n/n$ are inside or
outside of a suitable neighborhood of $m(\beta_n,K_n)$. Then part
(a) of Theorem \ref{thm:main} is derived by combining Theorem
\ref{thm:MDP} with Lemma \ref{lem:Sym1/2Small} and using symmetry.
--------------------------------------\fi
\iffalse ------------In section \ref{section:Thm61a} we prove part
(a) of Theorem \ref{thm:6.1} using Lemmas
\ref{lem:RepresentF}--\ref{lem:Substep2bOFThm61a}. -------------\fi
In section \ref{section:LogicANDPreLemmas} we derive  a number of
lemmas that are applied in section \ref{section:Thm61b} to part (b)
of Theorem \ref{thm:6.1}. In section \ref{section:Thm61b} we prove
part (b) of Theorem \ref{thm:6.1} using these lemmas together with
Lemmas \ref{lem:Step1OFThm61b}, \ref{lem:substep2aOFThm61b}, and
\ref{lem:Substep2bOFThm61b} and the weaker form of the standard
uniform integrability estimate in Proposition \ref{prop:weakinteg}.
In appendix A we prove that sequences 1a--5a satisfy the limits in hypothesis
(iii$^\prime$) of Theorem \ref{thm:main}. In appendix B we prove the moderate
deviation principle in
part (a) of Theorem \ref{thm:MDP}. This result is used in the proof of one of our
main results in part (a) of Theorem \ref{thm:main}.

\iffalse In Figure 3, which appears at the end of the dissertation,
we give the logic graph of the proof of parts (a) and (b) of Theorem
\ref{thm:6.1}. In that figure the lemma at the tail of a arrow is
used to prove the lemma at the head of the same arrow. Figure 3
exhibits the complex web of relationships among the fourteen lemmas
in sections \ref{section:Thm61a} and \ref{section:Thm61b} and
Proposition \ref{prop:weakinteg}, indicating the subtlety and
complexity of the proof of Theorem \ref{thm:6.1}.

FIGURE 3 HERE. BIGGER ARROWS ON FIGURE 3.

\begin{figure}[h]
\begin{center}
\epsfig{file=logicgraph.eps,width=12cm} \iffalse \vspace{.5in} \fi
\caption{\footnotesize The relationships between the lemmas and a
proposition in section \ref{section:Thm61a} and \ref{section:Thm61b}
that are used to prove Theorem \ref{thm:6.1}.}
\end{center}
\end{figure}
\fi

\skp
\noi
{\bf Acknowledgement.}  The research of both authors was supported
in part by a grant from the National Science Foundation (NSF-DMS-0604071). We thank Peter T. Otto for permission
to use Figures 1 and 2.

%===================== new section =====================================
\section{Phase-Transition Structure of the Mean-Field B-C Model}
\label{section:BEGmodel}
\beginsec

For $N \in \mathbb{N}$ the mean-field Blume-Capel model is defined
on the complete graph on $N$ vertices $1,2,\ldots,N$. The spin at
site $j \in \{1,2,\ldots,N\}$ is denoted by $\omega_j$, a quantity
taking values in $\Lambda = \{-1,0,1\}$. The Hamiltonian for this
model is defined by
\[ H_{N,K}(\omega) = \sum_{j=1}^N \omega_j^2 - \frac{K}{N} \left( \sum_{j=1}^N
\omega_j \right) ^2,
\]
where $K>0$ is a positive parameter representing the interaction
strength and $\omega=(\omega_1,\ldots,\omega_N)\in \Lambda^N$. We
will refer to this model as the mean-field B-C model.

Let $P_N$ be the product measure on $\Lambda^N$ with identical
one-dimensional marginals $\rho=\frac{1}{3}(\delta_{-1} + \delta_0 +
\delta_1)$. Then $P_N$ assigns the probability $3^{-N}$ to each
$\omega \in \Lambda^N$. For inverse temperature $\beta >0$ and for
$K>0$, the canonical ensemble for the mean-field B-C model is the
sequence of probability measures that assign to each subset $B$ of
$\Lambda^N$ the probability
\bea \label{eqn:P} P_{N,\beta, K}(B)
&=&
\frac{1}{Z_N(\beta,K)} \cdot \int_B \exp[-\beta H_{N,K}] d P_N \\
&=& \frac{1}{Z_N(\beta,K)} \cdot \sum_{\omega \in B} \exp[-\beta
H_{N,K}(\omega)] \cdot 3^{-N}, \nonumber
\eea
where \be
Z_N(\beta,K)=\int_{\Lambda^N} \exp[-\beta H_{N,K}]dP_N =
\sum_{\omega \in \Lambda^N} \exp[-\beta H_{N,K(\omega)}]\cdot 3^{-N}
\nonumber. \ee

It is useful to rewrite this measure in a different form. Define
$S_N(\omega)= \sum_{j=1}^N \omega_j$ and let $P_{N,\beta}$ be the
product measure on $\Lambda_N$ with identical one-dimensional
marginals
\[ \rho_\beta(d \omega_j)= \frac{1}{Z(\beta)} \cdot \exp(-\beta
\omega_j^2) \rho(d \omega_j),
\]
where $Z(\beta) = \int_\Lambda \exp(-\beta
\omega_j^2) \rho (d \omega_j)= (1+2e^{-\beta})/3$. Define
\[ P_{N,\beta}(d\omega)=\prod_{j=1}^N \rho_\beta (d
\omega_j)=\frac{1}{[Z(\beta)]^N} \prod_{j=1}^N \exp(-\beta
\omega_j^2)\rho(d \omega_j)
\]
and
\[\tilde{Z}_N(\beta, K)= \int_{\Lambda^N} \exp[N\beta K(S_N/N)^2]d
P_{N,\beta} = \frac{Z_N(\beta,K)}{[Z(\beta)]^N}.
\]
Then we have
\be \label{eqn:PnbetaK}
P_{N,\beta,K}(d\omega)=\frac{1}{\tilde{Z}_N(\beta,K)} \exp[N\beta K
(S_N(\omega)/N)^2] P_{N,\beta}(d\omega) .
\ee

For $t \in \R$ and $x \in \R$ we also define the
cumulate generating function
\be \label{eqn:Cbeta}
c_\beta(t)=\log
\int_\Lambda \exp(t \omega_1) \rho_\beta (d \omega_1) = \log
\left[\frac{1+e^{-\beta}(e^t+e^{-t})}{1+2e^{-\beta}}\right]
\ee
and the free-energy function
\be \label{eqn:GbetaK}
G_{\beta,K}(x)=\beta K x^2 - c_\beta(2\beta Kx).
\ee
We denote by
$\mathcal{M}_{\beta,K}$ the set of equilibrium macrostates of the
mean-field B-C model. As shown in Proposition 3.4 in
\cite{EllOttTou}, $\mathcal{M}_{\beta,K}$ can be
characterized as the set of global minimum points of $G_{\beta,K}$:
\[ \mathcal{M}_{\beta,K}= \{ x \in [-1,1]: x \ \textrm{is the
global minimum points of }G_{\beta,K}(x)  \}. \] In
\cite{EllOttTou} $\mathcal{M}_{\beta,K}$ is denoted by
$\tilde {\mathcal {E}}_{\beta,K}$.

The critical inverse temperature for the mean-field B-C model is
$\beta_c=\log 4$. For $0<\beta \le \beta_c$, the next theorem states
that $\mathcal{M}_{\beta,K}$ exhibits a continuous bifurcation as
$K$ increases through a value $K(\beta)$. This bifurcation
corresponds to a second-order phase transition, and the curve $\{
(\beta, K(\beta)), 0<\beta<\beta_c \}$ is called the second-order
curve. The point $(\beta_c, K(\beta_c))$ is called the tricritical
point. Theorem \ref{thm:secondorder} is proved in Theorem 3.6 in
\cite{EllOttTou}, where $K(\beta)$ is denoted by
$K_c^{(2)}(\beta)$.

\begin{thm}
\label{thm:secondorder} For $0 < \beta \leq \beta_c$, we define
\[
\iffalse \label{eqn:kbeta} \fi
K(\beta) = {1}/[{2\beta c''_\beta(0)}] = ({e^\beta
+ 2})/({4\beta}).
\]
For these values of $\beta$,
$\mathcal{M}_{\beta,K}$ has the following structure.

{\em(a)} For $0 < K \leq K(\beta)$, ${\mathcal{M}}_{\beta,K} =
\{0\}$.

{\em(b)} For $K > K(\beta)$, there exists ${m}(\beta,K) > 0$ such
that ${\mathcal{M}}_{\beta,K} = \{\pm m(\beta,K) \}$.

{\em(c)} ${m}(\beta,K)$ is a positive, increasing, continuous
function for $K > K_c(\beta)$, and as $K \goto (K(\beta))^+$,
$m(\beta,K) \goto 0$. Therefore, ${\mathcal{M}}_{\beta,K}$ exhibits
a continuous bifurcation at $K(\beta)$.
\end{thm}

For $\beta > \beta_c$, the next theorem states that
$\mathcal{M}_{\beta,K}$ exhibits a discontinuous bifurcation as $K$
increases through a value $K_1(\beta)$. This bifurcation corresponds
to a first-order phase transition, and the curve $\{ (\beta,
K_1(\beta)), \beta>\beta_c \}$ is called the first-order curve.
Theorem \ref{thm:firstorder} is proved in Theorem 3.8 in
\cite{EllOttTou}, where $K_1(\beta)$ is denoted by
$K^{(1)}_c(\beta)$.

\begin{thm}
\label{thm:firstorder} For $\beta > \beta_c$,
$\mathcal{M}_{\beta,K}$ has the following structure in terms of the
quantity $K_1(\beta)$, denoted by $K_c^{(1)}(\beta)$ in {\em
\cite{EllOttTou}} and defined implicitly for $\beta > \beta_c$ on
page {\em 2231} of {\em \cite{EllOttTou}}.

{\em(a)} For $0 < K < K_1(\beta)$, ${\mathcal{M}}_{\beta,K} =
\{0\}$.

{\em(b)} For $K = K_1(\beta)$ there exists $m(\beta,K_1(\beta)) > 0$
such that ${\mathcal{M}}_{\beta,K_1(\beta)} = \{0,\pm
m(\beta,K_1(\beta))\}$.

{\em(c)} For $K > K_1(\beta)$ there exists $m(\beta,K) > 0$ such
that ${\mathcal{M}}_{\beta,K} = \{\pm m(\beta,K)\}$.

{\em(d)} $m(\beta,K)$ is a positive, increasing, continuous function
for $K \geq K_1(\beta)$, and as $K \goto K_1(\beta)^+$, $m(\beta,K)
\goto m(\beta,K_1(\beta)) > 0$.  Therefore,
${\mathcal{M}}_{\beta,K}$ exhibits a discontinuous bifurcation at
$K_1(\beta)$.
\end{thm}

The positive quantity $m(\beta, K)$ in Theorems
\ref{thm:secondorder} and \ref{thm:firstorder}  is called the
thermodynamic magnetization. In the next section we describe the
asymptotic behavior of the finite-size magnetization for suitable
sequences $(\beta_n, K_n)$ and relate this to the asymptotic
behavior of the thermodynamic magnetization $m(\beta_n, K_n)$.

%======================== new section =====================================
\section{Asymptotic Behavior of $E_{n, \beta_n, K_n}\{|S_n/n|\}$}
\label{section:resultPrevious}
\beginsec

For $\beta>0$ and $K>0$ the finite-size magnetization is defined as
\[ E_{N, \beta, K} \{ |S_N/N| \} = \int_{\Omega_N} |S_N/N| d P_{N, \beta,
K},
\]
where $P_{N, \beta, K}$ denotes the measure defined in (\ref{eqn:P})--(\ref{eqn:PnbetaK}).
In this section we describe the asymptotic behavior of
$E_{n,\beta_n,K_n}\{|S_n/n|\}$ for suitable sequences $(\beta_n,
K_n)$ lying in the phase-coexistence region. In this formulation we let $N = n$ in the
finite-size magnetization; i.e., we let the number of spins $N$
coincide with the index $n$ parametrizing the sequence $(\bn,\kn)$.

The phase-coexistence region is defined as the set of all points in
the positive $\beta$-$K$ quadrant for which
${\mathcal{M}}_{\beta,K}$ cosists of more than one value. According
to Theorems \ref{thm:secondorder} and \ref{thm:firstorder}, the
phase-coexistence region consists of all points above the
second-order curve, above the tricritical point, on the first-order
curve, and above the first-order curve; that is,
\[ \{(\beta,K): 0 < \beta \le \beta_c, K > K(\beta) \ \textrm{and} \ \beta > \beta_c, K \ge K_1(\beta)
\}.
\]
For a class of sequences $(\beta_n, K_n)$ lying in the
phase-coexistence region for all sufficiently large $n$ and
converging either to a second-order point or to the tricritical
point, Theorem \ref{thm:3.1} describes the asymptotic behavior of
the thermodynamic magnetization $ m(\beta_n,K_n) \goto 0$. The
asymptotic behavior is related to the unique positive, global
minimum point of the Ginzburg-Landau polynomial, which is defined in
hypothesis (iii) of the theorem.

Theorem \ref{thm:3.1} is a special case of the main theorem in
\cite{EllMacOtt1}, Theorem 3.2. In that paper we describe six
different sequences that satisfy the hypotheses of Theorem
\ref{thm:3.1}. The first five of these sequences are revisited in section
\ref{section:verify6} of this paper, where we show that they satisfy the
hypotheses of our main theorem, Theorem \ref{thm:main}.
These five sequences, labeled 1a--5a, are summarized in Table 5.1.
\iffalse With several modifications the hypotheses of the next
theorem are also the hypotheses of theorems related to the
asymptotic behavior of $E_{n, \beta_n, K_n}\{|S_n/n|\}$ and $E_{n,
\beta_n, K_n}\{||S_n/n|-m(\beta_n, K_n)|\}$. These are Theorem 4.1
in \cite{EllMacOtt2010} and the main theorem, Theorem 4.1 in this
paper, respectively. \fi The main conclusion of Theorem
\ref{thm:3.1} about the rate at which $ m(\beta_n,K_n) \goto 0$ will
be used in the proofs of a number of results in this paper. \iffalse
Theorem \ref{thm:3.1} is also true if hypothesis (iii)(b) is
replaced by the following hypothesis: the set of global minimum
points of the Ginzburg-Landau polynomial $g$ equals $\{\pm
\bar{x}\}$ and $\{0, \pm \bar{x}\}$. But we do not need it in this
paper. \fi

\iffalse
Let $a_n$ and $b_n$ be any positive sequences converging to 0. We
write
\[a_n \sim b_n    \  \  \mbox{if}   \  \ \lim_{n \goto \infty} a_n/b_n =1\ \ \ \mbox{and} \ \ \ a_n \gg b_n    \  \  \mbox{if}   \  \ \lim_{n \goto \infty} a_n/b_n =\infty.
\]
\fi

\begin{thm}
\label{thm:3.1} Let $(\bn,\kn)$ be a positive sequence that
converges either to a second-order point $(\beta,K(\beta))$, $0 <
\beta < \bc$, or to the tricritical point $(\beta,K(\beta)) =
(\bc,K(\bc))$. We assume that $(\bn,\kn)$ satisfies the following
four hypotheses.
\begin{itemize}
\item[\em (i)] $(\beta_n,K_n)$ lies in the phase-coexistence region for all sufficiently large $n$.

\item [\em (ii)] The sequence $(\bn,\kn)$ is parametrized by $\alpha > 0$. This parameter regulates the speed of approach of
$(\bn,\kn)$ to the second-order point or the tricritical point in
the following sense{\em :}
\[
b = \lim_{n \goto \infty} n^\alpha (\bn - \beta) \ \mbox{ and } \ k
= \lim_{n \goto \infty} n^\alpha (\kn - K(\beta))
\]
both exist, and $b$ and $k$ are not both $0$; if $b \neq 0$, then
$b$ equals $1$ or $-1$.

\item[\em (iii)] There exists an even
polynomial $g$ of degree $4$ or $6$ satisfying $g(x) \goto \infty$
as $|x| \goto \infty$ together with the following two properties{\em
;} $g$ is called the Ginzburg-Landau polynomial.

\begin{itemize}
 % \vspace{-.1in}
  \item[\em (a)]  There exist $\alpha_0 > 0$ and $\theta > 0$ such that for all $ \alpha > 0$
\[
\lim_{n \goto \infty} n^{\alpha/\alpha_0} G_{\beta_n,K_n}
(x/n^{\theta\alpha}) = g(x)
\]
uniformly for $x$ in compact subsets of $\R$.

\item [\em (b)] There exists $\bar{x}>0$ such that the set of global minimum points of $g$ equals $\{\pm
\barx\}$.
\end{itemize}

\item [\em (iv)] Consider $\alpha_0>0$ and $\theta >0$ in hypothesis {\em (iii)(a)}. 
There exists a polynomial $H$ satisfying $H(x) \goto \infty$ as $|x| \goto \infty$
together with the following property{\em :} for all $\alpha > 0$ 
there exists $R > 0$ such that for all $n \in \N$ sufficiently
large and for all $x \in \R$ satisfying $|x/n^{\theta\alpha}| < R$,
$n^{\alpha/\alpha_0} G_{\beta_n,K_n} (x/n^{\theta\alpha}) \geq
H(x)$.
\end{itemize}
Under hypotheses {\em (i)--(iv)}, for any $\alpha > 0$
\[
m(\beta_n,K_n) \sim {\bar{x}}/{n^{\theta\alpha}}; \ \mbox{ i.e., }
\lim_{n \goto \infty} n^{\theta\alpha} m(\beta_n,K_n) = \bar{x}.
\]
If $b \not = 0$, then this becomes $m(\beta_n,K_n) \sim \barx|\beta
- \bn|^{\theta}$.
\end{thm}

Theorem \ref{thm:mainThm4.1Of2010} restates Theorem 4.1 in
\cite{EllMacOtt2010}. The hypotheses are those of Theorem
\ref{thm:3.1} for all $0<\alpha<\alpha_0$ together with the
inequality $0<\theta\alpha_0<1/2$. These hypotheses are satisfied
by sequences 1a--5a in Table 5.1 as well as by a sixth sequence
described in Theorem 4.6 in \cite{EllMacOtt1}. Part (a) of the next theorem
gives the rate at which $E_{n, \beta_n, K_n}\{|S_n/n|\}\goto 0$ for
$0<\alpha<\alpha_0$, and part (b) states that for the same values of
$\alpha$, $E_{n, \beta_n, K_n}\{|S_n/n|\}\sim m(\beta_n, K_n)$. Thus
Theorem \ref{thm:mainThm4.1Of2010} shows that the asymptotic
behavior of $E_{n, \beta_n, K_n}\{|S_n/n|\}$ coincides with that of
$m(\beta_n,K_n)$ for $0<\alpha<\alpha_0$. Theorem 4.2 in
\cite{EllMacOtt2010} shows that for $\alpha>\alpha_0$, $m(\beta_n,
K_n)$ converges to 0 asymptotically faster than $E_{n, \beta_n, K_n}
\{|S_n/n|\}$.
\iffalse
; i.e., $E_{n, \beta_n, K_n}\{|S_n/n|\}\gg m(\beta_n,
K_n)$.
\fi

\begin{thm}
\label{thm:mainThm4.1Of2010} Let $(\beta_n, K_n)$ be a positive
sequence parametrized by $\alpha>0$ and converging either to a
second-order point $(\beta, K(\beta))$, $0<\beta<\beta_c$, or to the
tricritical point $(\beta_c, K(\beta_c))$. We assume that $(\beta_n,
K_n)$ satisfies the hypotheses of Theorem {\em \ref{thm:3.1}} for
all $0<\alpha<\alpha_0$. We also assume the inequality
$0<\theta\alpha_0<1/2$. The following conclusions hold.

{\em (a)} For all $0 < \alpha < \alpha_0$
\[
E_{n,\beta_n, K_n}\{|S_n/n|\} \sim {\bar{x}}/{n^{{\theta\alpha}}}; \
\mbox{i.e., } \lim_{n \goto \infty} n^{\theta\alpha} E_{n, \beta_n,
K_n} \left\{ \left| S_n/n \right| \right\} = \bar{x}.
\]

{\em (b)}  For all $0 < \alpha < \alpha_0$, $E_{n,\beta_n,
K_n}\{|S_n/n|\} \sim m(\beta_n, K_n)$.
\end{thm}

In Theorem \ref{thm:main} in the next section we state our main result on the rate at which $E_{n, \beta_n,K_n}\{
| |S_n/n| - m(\beta_n, K_n) |\}$ converges to 0 for $0<\alpha<\alpha_0$.
\iffalse
We also
compare this rate with the asymptotic behaviors of the finite-size
magnetization $E_{n, \beta_n,K_n}\{ |S_n/n|  \} \goto 0$ and of the
thermodynamic magnetization $m(\beta_n, K_n) \goto 0$. After the
statement, \fi
We then explain how Theorem \ref{thm:main}
generalizes Theorem \ref{thm:mainThm4.1Of2010}.

%======================== new section =====================================

\section{Asymptotic Behavior of $E_{n,\beta_n,K_n}\{| |S_n/n| - m(\beta_n,K_n) | \}$}
\label{section:resultmain}
\beginsec

We denote by $E_{n, \beta_n,K_n}$ expectation with respect to the
measure $P_{n,\beta_n,K_n}$. Theorem \ref{thm:main} is our main
result. In this theorem we investigate the asymptotic behavior of
the expectation $E_{n, \beta_n,K_n}\{ | |S_n/n| - m(\beta_n, K_n)  |
\}$ under the hypotheses of Theorem \ref{thm:3.1} and an additional
hypothesis (iii$^\prime$). Part (a) of Theorem \ref{thm:main} states that
the expected value of the fluctuations of $|S_n/n|$ around $m(\beta_n, K_n)$ is
asymptotic to $\bar{z}/n^\kappa$, where $\kappa = \frac{1}{2}
(1-\alpha/\alpha_0)+\theta\alpha$ and $\bar{z}
>0$ is given explicitly. Compared with the conclusion of Theorem \ref{thm:mainThm4.1Of2010},
part (a) of Theorem \ref{thm:main} is a
more refined statement. As we showed in the introduction, it yields the conclusion of Theorem
\ref{thm:mainThm4.1Of2010} as a corollary. The rate
$\bar{z}/n^\kappa$ at which $E_{n, \beta_n,K_n}\{ | |S_n/n| -
m(\beta_n, K_n)  | \}$ converges to 0 is much faster than the rate
$\bar{x}/n^{{\theta\alpha}}$ at which $E_{n, \beta_n,K_n}\{ |S_n/n|
\}$ and $m(\beta_n, K_n)$ converge to 0 separately. \iffalse The
result of Theorem \ref{thm:mainThm4.1Of2010} is also proved from
this statement. (change on April 25)\fi We comment on the hypotheses
of Theorem \ref{thm:main} at the end of this section.

%Add on April 25
Part (a) of Theorem \ref{thm:main} is proved in section
\ref{section:proofmain}. Part (b) of Theorem \ref{thm:main} asserts
that the hypotheses of this theorem are satisfied by sequences
1a--5a in Table 5.1. This is discussed in section
\ref{section:verify6}. For each of these sequences the Ginzburg-Landau polynomial
has degree 4 or 6.
%

\iffalse Let $a_n$ be a positive sequence converging to 0. In
stating the result on the rate at which $E_{n, \beta_n,K_n}\{ |
|S_n/n| - m(\beta_n, K_n)  | \} \goto 0$, we write
\[E_{n,
\beta_n,K_n}\{ | |S_n/n| - m(\beta_n, K_n)  | \} \sim a_n    \  \
\mbox{if}   \  \ \lim_{n \goto \infty} E_{n, \beta_n,K_n}\{ |
|S_n/n| - m(\beta_n, K_n)  | \}/a_n =1.\] \fi

\begin{thm}
\label{thm:main} Let $(\beta_n,K_n)$ be a positive sequence
converging either to a second-order point $(\beta, K(\beta))$,
$0<\beta<\beta_c$, or to the tricritical point $(\beta,
K(\beta))=(\beta_c, K(\beta_c))$. We assume that $(\beta_n,K_n)$
satisfies the hypotheses of Theorem {\em \ref{thm:3.1}} for all
$0<\alpha<\alpha_0$. We also assume the following additional
hypothesis on the Ginzburg-Landau polynomial $g$.
\begin{itemize}
\item[\em (iii$^\prime$)] Assume that $g$ has degree $4$. Then $\theta\alpha_0$
lies in the interval $[1/4, 1/2)$. In addition, for all $0 < \alpha < \alpha_0$
and for $j=2,3,4$
\[
\lim_{n \goto \infty} n^{\alpha/\alpha_0 - j \theta\alpha}
  G_{\beta_n,K_n}^{(j)}(m(\beta_n,K_n)) = g^{(j)}(\bar{x})>0.
\]
Assume that $g$ has degree $6$. Then $\theta\alpha_0$ lies in the
interval $[1/6, 1/2)$. In addition, for all $0 < \alpha < \alpha_0$ and for $j=2,3,4,5,6$
\[ \lim_{n \goto \infty}
n^{\alpha/\alpha_0 - j \theta\alpha}
  G_{\beta_n,K_n}^{(j)}(m(\beta_n,K_n)) = g^{(j)}(\bar{x})>0.
\]
\end{itemize}
For $\alpha \in (0, \alpha_0)$ we also define $\kappa = \frac{1}{2}
(1-\alpha/\alpha_0)+\theta\alpha$. Then for all $0<\alpha<\alpha_0$
the following conclusions hold.

{\em (a)} We have the asymptotic behavior
\[ E_{n,\beta_n,K_n} \{||S_n/n|-m(\beta_n,K_n)|\} \sim
\bar{z}/n^\kappa,
\]
where $\bar{z} = \left(2 /(\pi g^{(2)}(\bar{x}))\right)^{1/2}$;
i.e., $\lim_{n \goto \infty} n^\kappa E_{n,\beta_n,K_n}
\{||S_n/n|-m(\beta_n,K_n)|\}=\bar{z}$.

{\em (b)} The hypotheses of this theorem are satisfied by sequences {\em 1a--5a} in
Table {\em 5.1}.
\end{thm}

The hypotheses of Theorem \ref{thm:main} are those of Theorem
\ref{thm:3.1} together with the additional hypothesis
(iii$^{\prime}$) for all $0<\alpha<\alpha_0$. The latter hypothesis
takes two related forms depending on whether $g$ has degree 4 or
degree 6. In this hypothesis, the assumption on $\theta\alpha_0$
yields the inequality $0<\theta\alpha_0 <1/2$, which is required by
the moderate deviation principle stated in Theorem
\ref{thm:MDP}. \iffalse---------------------------------------\fi
Hypothesis (iii$^{\prime}$) also assumes both the asymptotic behavior of
certain derivatives of $n^{\alpha/\alpha_0}G_{\beta_n,K_n}$
evaluated at $m(\beta_n,K_n)$ and the positivity of the
corresponding derivatives of $g$ evaluated at the positive global
minimum point $\bar{x}$. These assumptions are needed in the proof
of Lemma \ref{lem:G(xmn)miusG(mn)}, a key result needed to prove
part (b) of Theorem \ref{thm:6.1}, which in turn
yields part (a) of Theorem \ref{thm:main}. The proof of that lemma
also requires the fact assumed in hypothesis (iii$^\prime$) that
$\theta\alpha_0$ lies in the interval [1/4, 1/2) or [1/6, 1/2)
depending on whether $g$ has degree 4 or degree 6.

In the next section we outline how to verify the hypotheses of
Theorem \ref{thm:main} for sequences 1a--5a in Table 5.1.
\iffalse
We will see that the asymptotic behavior of the derivatives of
$n^{\alpha/\alpha_0}G_{\beta_n,K_n}$ required in hypothesis
(iii$^\prime$) of Theorem \ref{thm:main} is a consequence of limits
given in Lemma \ref{lem:5.1}.
\fi

%===================== new section =====================================

\section{Verification of Hypotheses of Theorem \ref{thm:main} for Sequences 1a--5a}
\label{section:verify6}
\beginsec

Table 5.1 summarizes five sequences $(\beta_n, K_n)$ introduced in section 4 of
\cite{EllMacOtt1}. Depending on the inequalities on the
coefficients, sequences 1, 2, 3, and 5 each have two cases labeled a
and b, and sequence 4 has three cases labeled a, b, and c. All five
sequences 1a--5a lie in the phase-coexistence region for all
sufficiently large $n$ as required by hypothesis (i) of Theorem
\ref{thm:3.1}. \iffalse After defining the sequences, we then
restate the verification of the hypotheses of Theorem \ref{thm:3.1}
for sequences 1a--5a and verify the hypotheses of Theorem
\ref{thm:main} for sequences 1a--5a. Sequence 6 introduced in
\cite{EllMacOtt1} does not satisfy hypothesis (iii$^\prime$) in
Theorem \ref{thm:main}. \fi

\begin{table}[h]
\begin{center}
%\label{table:6sequences}
\begin{tabular}{||l|l|l|l|l|l|l||} \hline \hline
{{\bf Seq.}} & {\bf Defn.} &{\bf Case} &{\bf
Ineq.} &{\bf Region} & \boldmath $\mathcal {M}_g$ \unboldmath &{\bf Thm. in \cite{EllMacOtt1}} \\
\hline \hline 1 & (\ref{eqn:seq1}) & a & $K'(\beta)b-k <0$ & Ph-CR
& $\{\pm \bar{x}\}$ & Thm.4.1 \\
\cline{3-6} & & b & $K'(\beta)b-k >0$ & 1-PhR & $\{0\}$ & \\
\hline \hline 2 & (\ref{eqn:seq2}) & a & $(K^{(p)}(\beta)-\ell)b^p
<0$ & Ph-CR
& $\{\pm \bar{x}\}$ & Thm.4.2 \\
\cline{3-6} & & b & $(K^{(p)}(\beta)-\ell)b^p
<0$ & 1-PhR & $\{0\}$ & \\
\hline \hline 3 & (\ref{eqn:seq3}) & a & $K'(\beta_c)b-k <0$ & Ph-CR
& $\{\pm \bar{x}\}$ & Thm.4.3 \\
\cline{3-6} & & b & $K'(\beta_c)b-k >0$ & 1-PhR & $\{0\}$ & \\
\hline \hline 4 & (\ref{eqn:seq4}) & a & $\ell>\ell_c, \tilde{\ell}
\in \mathbb{R}$ & Ph-CR
& $\{\pm \bar{x}\}$ & Thm.4.4 \\
\cline{3-6} & & b & $\ell=\ell_c, \tilde{\ell} > K'''_1(\beta_c)$ & Ph-CR & $\{0, \pm \bar{x}\}$ & \\
\cline{3-6} & & c & $\ell<\ell_c, \tilde{\ell} \in
\mathbb{R}$ & 1-PhR & $\{0\}$ & \\
\hline \hline 5 & (\ref{eqn:seq5}) & a & $\ell > K''(\beta_c)$ &
Ph-CR
& $\{\pm \bar{x}\}$ & Thm.4.5 \\
\cline{3-6} & & b & $\ell<K''(\beta_c)$ & 1-PhR & $\{0\}$ & \\
\iffalse ************************************* \hline \hline 6 &
(\ref{eqn:seq6}) & a & $(K^{(p)}(\beta)-l)(-1)^p <0$ & Ph-CR
& $\{\pm \bar{x}\}$ & Thm.4.6 \\
\cline{3-6} & & b & $(K^{(p)}(\beta)-l)(-1)^p <0$ & 1-PhR & $\{0\}$
& \\    **************************************\fi
\hline \hline
\end{tabular}
\end{center}
% \vspace{-.025in}
%\begin{center}
\footnotesize Table 5.1: The equation where each of the 5 sequences is defined and the
inequalities on the coefficients guaranteeing that each sequence lies
in the phase-coexistence region (Ph-CR) or in the single-phase
region (1-PhR). The next-to-last column states the structure of the
set $\mathcal {M}_g$ of global minimum points of the Ginzburg-Landau
polynomial $g$ for each sequence in terms of a
positive number $\bar{x}$ that can be explicitly calculated. The theorems in
\cite{EllMacOtt1} where this information is verified are also given.
%\end{center}
\end{table}

The hypotheses of Theorem \ref{thm:main} consist of the hypotheses
of Theorem \ref{thm:3.1} for all $0<\alpha<\alpha_0$ and
hypothesis (iii$^\prime$). Hypothesis (iii$^\prime$) takes two forms
depending on the degree of the Ginzburg-Landau polynomial $g$. When
$g$ has degree 4, $\theta\alpha_0$ is assumed to lie in the interval
$[1/4, 1/2)$ and for all $\alpha \in (0, \alpha_0)$ and for $j=2,3,4$
\be
\label{eqn: hypiii'4}
\lim_{n \goto \infty}
n^{\alpha/\alpha_0 - j \theta\alpha}
  G_{\beta_n,K_n}^{(j)}(m(\beta_n,K_n)) = g^{(j)}(\bar{x}) > 0.
\ee
When $g$ has degree 6,
$\theta\alpha_0$ is assumed to lie in the interval $[1/6, 1/2)$ and for
all $\alpha \in (0, \alpha_0)$ and for $j=2,3,4, 5, 6$
\be
\label{eqn: hypiii'6}
\lim_{n \goto \infty}
n^{\alpha/\alpha_0 - j \theta\alpha}
  G_{\beta_n,K_n}^{(j)}(m(\beta_n,K_n)) = g^{(j)}(\bar{x}) > 0.
\ee
In this section we verify for
sequences 1a--5a that when $g$ has degree 4, we have $\theta\alpha_0 \in
[1/4, 1/2)$ and $g^{(j)}(\bar{x})>0$ for $j=2,3,4$ and that when $g$
has degree 6, we have $\theta\alpha_0 \in [1/6, 1/2)$ and
$g^{(j)}(\bar{x})>0$ for $j=2,3,4,5,6$. The verification of the
limits in (\ref{eqn: hypiii'4}) and (\ref{eqn: hypiii'6}) is carried
out in appendix A.

Sequence 6 introduced in Theorem 4.6 in \cite{EllMacOtt1} does not satisfy
hypothesis (iii$^\prime$) in Theorem \ref{thm:main}. In this case
$g$ has degree 4, but $\theta\alpha_0$ does not lie in the interval
$[1/4, 1/2)$.

The first two sequences converge to a second-order point $(\beta,
K(\beta))$, $0<\beta<\beta_c$, and the last three sequences converge
to the tricritical point $(\beta_c, K(\beta_c))$. For each sequence
1a--5a, the hypotheses of Theorem \ref{thm:3.1} are verified in
Theorems 4.1--4.5 in \cite{EllMacOtt1}. We follow the same method
used in that paper to verify hypothesis (iii$^{\prime}$) in Theorem
\ref{thm:main} for sequences 1a--5a.
\iffalse We first do some preparatory work for this. \fi
Hypothesis (iii$^{\prime}$) of Theorem \ref{thm:main}
takes two forms depending on whether the degree of the
Ginzburg-Landau polynomial $g$ is 4 or 6. We must
verify that $\theta\alpha_0$ lies in a certain interval and that
\be
\label{eqn:thislimit}
\lim_{n \goto \infty} n^{\alpha/\alpha_0 - j \theta\alpha}
  G_{\beta_n,K_n}^{(j)}(m(\beta_n,K_n)) = g^{(j)}(\bar{x})>0
\ee
for $j = 2,3,4$ when $g$ has degree 4 and for $j = 2,3,4,5,6$ when $g$ has degree 6.
The function $G_{\beta, K}$ is defined in
$(\ref{eqn:Cbeta})$--$(\ref{eqn:GbetaK})$.

\iffalse
Lemma \ref{lem:5.1} gives
a condition that implies the limit in the last display. This
condition is given in (\ref{eqn: 1assumeLem51}) when the degree of
$g$ is 4 and in (\ref{eqn: 2assumeLem51}) when the degree of $g$ is 6.
\fi

It is straightforward to show that the limit in (\ref{eqn:thislimit})
holds for a given $j$ provided the following limit holds uniformly
for $x$ in compact subsets of $\R$:
\be
\label{eqn:needthistoo}
\lim_{n \goto \infty} n^{\alpha/\alpha_0 - j \theta\alpha}
  G_{\beta_n,K_n}^{(j)}(x/n^{\theta\alpha}) = g^{(j)}(x).
\ee
The proof that the uniform convergence
in (\ref{eqn:needthistoo}) implies the limit
in (\ref{eqn:thislimit}) uses the fact that $n^{\theta\alpha} m(\beta_n,k_n) \goto \bar{x}$ [Thm.\ \ref{thm:3.1}].
\iffalse
Conditions (\ref{eqn: 1assumeLem51}) and (\ref{eqn: 2assumeLem51})
involve the uniform convergence of $n^{\alpha/\alpha_0 - j
\theta\alpha}  G_{\beta_n,K_n}^{(j)}(x/n^{\theta\alpha})$ to
$g^{(j)}(x)$ for $j=2,3,4$ when $g$ has degree 4 and for
$j=2,3,4,5,6$ when $g$ has degree 6. This
\fi
The uniform convergence in (\ref{eqn:needthistoo}) can be
obtained formally by taking the $j$-th derivative of the uniform
convergence limits in hypothesis (iii)(a) of Theorem \ref{thm:3.1}:
\[
\lim_{n \goto \infty} n^{\alpha/\alpha_0} G_{\beta_n,K_n}
(x/n^{\theta\alpha}) = g(x).
\]
The verification of the uniform convergence limits in (\ref{eqn:needthistoo}),
\iffalse
in (\ref{eqn: 1assumeLem51}) and (\ref{eqn: 2assumeLem51})
\fi
and thus the verification of the limits
(\ref{eqn: hypiii'4}) and (\ref{eqn: hypiii'6}) in hypothesis
(iii$^\prime$), depend on asymptotic properties of the Taylor
expansions of $G_{\beta_n,K_n}^{(j)}(x/n^{\theta\alpha})$. This
analysis closely parallels the proof of Theorem \ref{thm:3.1}, which
is based on a similar analysis of the Taylor expansions of
$G_{\beta_n,K_n}(x/n^{\theta\alpha})$ carried out in
\cite{EllMacOtt1}. The straightforward but tedious
calculations can be found in appendix A.

We now define the five sequences $(\beta_n,K_n)$ and summarize the
verification of the hypotheses of Theorem \ref{thm:main} for them.
\iffalse For every sequence we restate the verification of the
hypotheses of Theorem \ref{thm:3.1} proved in \cite{EllMacOtt2010}
because these hypotheses are needed to verify hypothesis
(iii$^\prime$) in Theorem \ref{thm:main}. \fi

\skp

%++++++++++++++++++++  sequence 1  ++++++++++++++++++++++++++++++++

\noindent {\bf Sequence 1a}

\noindent {\it Definition}. Given $0 < \beta < \bc$, $\alpha
> 0$, $b \in \{1,0,-1\}$, and $k \in \R$, $k \not = 0$, the sequence
is defined by \be \label{eqn:seq1} \beta_n = \beta + b/n^\alpha \
\mbox{ and } \ K_n = K(\beta) + k/n^\alpha. \ee This sequence
converges to the second-order point $(\beta,K(\beta))$ along a ray
with slope $k/b$ if $b \not = 0$. We assume that $K'(\beta)b-k <0$.
Under this assumption it is proved in Theorem 4.1 in \cite{EllMacOtt1}
that sequence 1 satisfies the hypotheses of Theorem
\ref{thm:3.1} with $\alpha_0=1/2$ and $\theta=1/2$. When $K'(\beta)b-k <0$, we refer to sequence 1 as sequence 1a.

\skp \noi {\it Hypothesis {\em (iii$^\prime$)} in Theorem {\em
\ref{thm:main}} for sequence {\em 1a}.} Since $\alpha_0=1/2$ and $\theta=1/2$,
$\theta\alpha_0$ lies in the interval $[1/4,1/2)$ as required by
hypothesis (iii$^\prime$). The limits in (\ref{eqn: hypiii'4}) for
$j=2,3,4$ are proved in appendix A. We now prove that
$g^{(j)}(\bar{x}) >0$ for $j=2, 3, 4$ using the formulas
for $g$ and $\bar{x}$ in Theorem 4.1 in \cite{EllMacOtt1}. Let
$c_4(\beta)=(e^\beta + 2)^2 (4-e^\beta)/8 \cdot 4 !$. Since $0<\beta<\beta_c=\log 4$, we have $e^\beta
< e^{\beta_c}=4$, which implies $c_4(\beta)>0$. Since
$K'(\beta)b-k <0$, these formulas yield
\iffalse
\[ g'(x)=2\beta (K'(\beta)b-k)x + 4 c_4(\beta) x^3
\]
and
\fi
\[
g^{(2)}(\bar{x}) = 2\beta (K'(\beta)b-k)+ 3 \cdot 4 c_4(\beta)
\bar{x}^2  = 4\beta (k-K'(\beta)b) > 0,
\]
\iffalse = 2\beta (K'(\beta)b-k)-6\beta (K'(\beta)b-k) \nonumber
\\ \fi
\[
g^{(3)}(\bar{x})= 4! c_4(\beta)\bar{x} >0, \ \mbox{ and } \  g^{(4)}(\bar{x})= 4! c_4(\beta) >0.
\]
Thus under the condition $K'(\beta)b-k<0$ sequence 1a  satisfies all
the hypotheses of Theorem \ref{thm:main}.

\skp

\noindent {\bf Sequence 2a}

\noindent {\it Definition}. Given $0 < \beta < \bc$, $\alpha
> 0$, $b \in \{1, -1\}$, an integer $p \ge 2$, and a real number $\ell \neq K^{(p)}(\beta)$, the sequence
is defined by \be \label{eqn:seq2} \beta_n = \beta + b/n^\alpha \
\mbox{ and } \ K_n = K(\beta) + \sum^{p-1}_{j=1}
K^{(j)}(\beta)b^j/(j! n^{j\alpha}) + \ell b^p / (p!n^{p\alpha}). \ee
This sequence converges to the second-order point $(\beta,
K(\beta))$ along a curve that coincides with the second-order curve
to order $n^{-(p-1)\alpha}$. We assume that
$(K^{(p)}(\beta)-\ell)b^p <0$. Under this assumption it is proved in Theorem 4.2 in \cite{EllMacOtt1}
that sequence 2 satisfies the hypotheses of Theorem \ref{thm:3.1} with
$\alpha_0=1/(2p)$ and $\theta=p/2$. When $(K^{(p)}(\beta)-\ell)b^p <0$,
we refer to sequence 2 as sequence 2a.

\skp \noi {\it Hypothesis {\em (iii$^\prime$)} in Theorem {\em
\ref{thm:main}} for sequence {\em 2a}.} Since $\alpha_0=1/(2p)$ and $\theta=p/2$,
$\theta\alpha_0$ lies in the interval $[1/4,1/2)$ as required by
hypothesis (iii$^\prime$). The limits in (\ref{eqn: hypiii'4}) for
$j=2,3,4$ are proved in appendix A. We now prove that
$g^{(j)}(\bar{x}) >0$ for $j=2, 3, 4$ using the formulas
for $g$ and $\bar{x}$ in Theorem 4.2 in \cite{EllMacOtt1}. Let
$c_4(\beta)=(e^\beta + 2)^2 (4-e^\beta)/8 \cdot 4!$. Since $0<\beta<\beta_c=\log 4$, we have $e^\beta < e^{\beta_c}=4$,
which implies $c_4(\beta)>0$. Since $(K^{(p)}(\beta)-\ell)b^p <0$, these formulas
yield
\iffalse
\[
g'(x)=\frac{2}{p!}\beta (K^{(p)}(\beta)-\ell)b^p x + 4
c_4(\beta) x^3
\]
and
\fi
\[g^{(2)}(\bar{x}) = \frac{2}{p!}\beta
(K^{(p)}(\beta)-\ell)b^p + 3 \cdot 4c_4(\beta)
\bar{x}^2
\iffalse = \frac{2}{p!}\beta
(K^{(p)}(\beta)-\ell)b^p-\frac{6}{p!}\beta (\ell-K^{(p)}(\beta))b^p
\nonumber \\ \fi
= \frac{4}{p!}\beta (\ell-K^{(p)}(\beta))b^p >0,
\]
\[
g^{(3)}(\bar{x})= 4! c_4(\beta)\bar{x} >0, \
\mbox{ and } \ g^{(4)}(\bar{x})= 4! c_4(\beta) >0.
\]
Thus under the condition $(K^{(p)}(\beta)-\ell)b^p<0$ sequence 2a
satisfies all the hypotheses of Theorem \ref{thm:main}.

\skp

\noindent {\bf Sequence 3a}

\noindent {\it Definition}. Given $\alpha
> 0$, $b \in \{1,0,-1\}$, and $k \in \R$, $k \not = 0$, the sequence
is defined by \be \label{eqn:seq3} \beta_n = \beta_c + b/n^\alpha \
\mbox{ and } \ K_n = K(\beta_c) + k/n^\alpha. \ee This sequence
converges to the tricritical point $(\beta_c, K(\beta_c))$ along a
ray with slope $k/b$ if $b \not = 0$. We assume that $K'(\beta_c)b-k
<0$. Under this assumption it is proved in Theorem 4.3 in \cite{EllMacOtt1}
that sequence 3 satisfies the hypotheses of
Theorem 3.1 with $\alpha_0=2/3$ and $\theta=1/4$. When $K'(\beta_c)b-k <0$, we refer to sequence 3 as sequence 3a.

\skp \noi {\it Hypothesis {\em (iii$^\prime$)} in Theorem {\em
\ref{thm:main}} for sequence {\em 3a}.} Since $\alpha_0=2/3$ and $\theta=1/4$,
$\theta\alpha_0$ lies in the interval $[1/6,1/2)$ as required by
hypothesis (iii$^\prime$). The limits in (\ref{eqn: hypiii'6}) for
$j=2,3,4,5,6$ are proved in appendix A. We now prove
that $g^{(j)}(\bar{x}) >0$ for $j=2, 3, 4,5,6$ using the
formulas for $g$ and $\bar{x}$ in Theorem 4.3 in \cite{EllMacOtt1}.
Let $c_6=9/40$. Since $K'(\beta_c)b-k <0$, these formulas yield
\iffalse
\[ g'(x)=2\beta_c (K'(\beta_c)b-k)x + 6 c_6 x^5,
\]
\fi
\[
g^{(2)}(\bar{x}) = 2\beta_c (K'(\beta_c)b-k)+ 5 \cdot 6 c_6
\bar{x}^4
\iffalse =
2\beta_c (K'(\beta_c)b-k)-10\beta_c (K'(\beta_c)b-k) \nonumber \\
\fi
= 8\beta_c (k-K'(\beta_c)b) >0,
\]
\[
g^{(3)}(\bar{x}) = 4 \cdot 5 \cdot 6 c_6 \bar{x}^3 >0, \  g^{(4)}(\bar{x})= 3 \cdot 4 \cdot 5 \cdot 6 c_6 \bar{x}^2 >0,
\]
\[g^{(5)}(\bar{x}) = 6! c_6 \bar{x} >0,  \ \mbox{ and } \
\ g^{(6)}(\bar{x}) = 6! c_6  >0.
\]
Thus under the condition $K'(\beta_c)b-k<0$ sequence 3a  satisfies
all the hypotheses of Theorem \ref{thm:main}.

\skp

\noindent {\bf Sequence 4a}

\noindent {\it Definition}. Given $\alpha
> 0$, a curvature parameter $\ell \in \mathbb{R}$, and another parameter $\tilde{\ell} \in \mathbb{R}$, the
sequence 4 is defined by \be \label{eqn:seq4} \beta_n = \beta_c +
1/n^\alpha \ \mbox{ and } \ K_n = K(\beta_c) + K'(\beta_c)/n^\alpha
+ \ell/(2n^\alpha) + \tilde{\ell}/(6n^{3\alpha}). \ee This sequence
converges from the right to the tricritical point $(\beta_c,
K(\beta_c))$ along the curve $(\beta, \tilde{K}(\beta))$, where for
$\beta > \beta_c$
\[\tilde{K}(\beta)=K(\beta_c)+ K'(\beta_c)(\beta-\beta_c) + \ell
(\beta-\beta_c)^2/2 + \tilde{\ell}(\beta-\beta_c)^3/6.
\]
The first-order curve $\{ (\beta, K_1(\beta)), \beta>\beta_c \}$ is shown in Figure 1 in the introduction.
In order to determine a condition on the coefficients
guaranteeing that sequence 4 satisfies the hypotheses of Theorem \ref{thm:3.1},
we must study $K_1(\beta)$ more closely.

Since $\lim_{\beta \goto \beta_c^+}
K_1(\beta)=K(\beta_c)$ \cite[Sects. 3.1, 3.3]{EllOttTou}, by
continuity we extend the definition of $K_1(\beta)$ from
$\beta>\beta_c$ to $\beta=\beta_c$ by define
$K_1(\beta_c)=K(\beta_c)$. In addition we must assume other properties of $K_1$ that
are stated in conjectures 1 and 2 on page 119 of \cite{EllMacOtt1}.
As a preliminary to stating these conjectures, we assume that the first three
right-hand derivatives of $K_1(\beta)$ exist at $\beta_c$ and
denote them by $K'_1(\beta_c)$, $K''_1(\beta_c)$, and
$K'''_1(\beta_c)$. We also define $\ell_c=
K''(\beta_c)-5/(4\beta_c)$. Conjectures 1 and 2 state the following:
(1) $K'_1(\beta_c)=K'(\beta_c)$ and (2) $K''_1(\beta_c) =
\ell_c < 0 < K''(\beta_c)$. These conjectures  are discussed in
detail in section 5 of \cite{EllMacOttUnpub} and are supported by
properties of the Ginzburg-Landau polynomials and numerical
calculations.

We assume that $\ell > \ell_c$, which by conjecture 1 equals $K''_1(\beta_c)$. Under this assumption
it is proved in Theorem 4.4 in \cite{EllMacOtt1}
that sequence 4 satisfies the hypotheses of Theorem \ref{thm:3.1}
with $\alpha_0=1/3$ and $\theta=1/2$. When $\ell > \ell_c$, we refer to sequence 4 as sequence 4a.

\skp \noi {\it Hypothesis {\em (iii$^\prime$)} in Theorem {\em
\ref{thm:main}} for sequence {\em 4a}.} Since $\alpha_0=1/3$ and $\theta=1/2$,
$\theta\alpha_0$ lies in the interval $[1/6,1/2)$ as required by
hypothesis (iii$^\prime$). The limits in (\ref{eqn: hypiii'6}) for
$j=2,3,4,5,6$ are proved in appendix A. Define
\be
\label{eqn:y}
y = \left(1 + \frac{3}{5}\beta_c (\ell - K''(\beta_c)) \right)^{1/2}.
\ee
Since $\ell > \ell_c = K''(\beta_c) - 5/(4\beta_c)$, we have $y > 1/2$.
We now prove that $g^{(j)}(\bar{x}) >0$ for $j=2, 3, 4,5,6$
using the formulas for $g$ and $\bar{x}$ in Theorem 4.4 in
\cite{EllMacOtt1}. Let $c_4=3/16$ and $c_6=9/40$. These formulas
yield
%%%%%%%%%%%%%%%%%%%%%%%%%%%%%%%%%%%%%%%%%%%%%%%%%%%%%%%%%%%%%%%%%%%%%%%%%%%%%%%%%%%%%%%%%%%%%%%%%%%%%%%%%%%%%%%%%%%%%%%%%%%
\iffalse
[ g'(x)=\beta_c (K''(\beta_c)-\ell)x - 16 c_4 x^3 + 6c_6 x^5,
\]
\bea g^{(2)}(\bar{x}) &=& \beta_c (K''(\beta_c)-l) - 48 \times \frac{3}{16} \times \frac{10}{9} [ 1+(1-\frac{3}{5}\beta_c (K''(\beta_c)-l))^{1/2} ] \nonumber \\
&&+ 30 \times \frac{9}{40} \times \frac{100}{81} [ 1+(1-\frac{3}{5}\beta_c(K''(\beta_c)-l))+2( 1- \frac{3}{5}\beta_c(K''(\beta_c)-l) )^{1/2} ]  \nonumber \\
&=& \beta_c (K''(\beta_c)-l) -10 -10( 1-\frac{3}{5}\beta_c (K''(\beta_c)-l) )^{1/2} \nonumber \\
&&+ \frac{50}{3} -5\beta_c( K''(\beta_c)-l ) + \frac{50}{3}( 1-\frac{3}{5}\beta_c (K''(\beta_c)-l)  )^{1/2} \nonumber \\
&=& -4\beta_c(K''(\beta_c)-l) + \frac{20}{3} + \frac{20}{3}
(1-\frac{3}{5}\beta_c (K''(\beta_c)-l))^{1/2} \nonumber. \eea  \fi
\iffalse Then we have \[ -4\beta_c(K''(\beta_c)-l)=
\frac{20}{3}y^2-\frac{20}{3}.
\] \fi
\[ g^{(2)}(\bar{x}) = \beta_c (K''(\beta_c)-\ell) - 3 \cdot 4 \cdot 4 c_4 \bar{x}^2 + 5 \cdot 6 c_6 \bar{x}^4
= \frac{20}{3}y^2+\frac{20}{3}y >0,
\]
\[
g^{(3)}(\bar{x}) = - 4! \cdot 4 c_4\bar{x} + 4 \cdot 5 \cdot 6 c_6\bar{x}^3 =
9\bar{x} \left( \frac{4}{3}+\frac{10}{3}y \right) >0,
\]
\[ g^{(4)}(\bar{x}) = - 4! \cdot 4 c_4 + 3 \cdot 4 \cdot 5 \cdot 6 c_6 \bar{x}^2 = -18 + 90(1+y)>0,
\]
\[ g^{(5)}(\bar{x}) = 6! c_6 \bar{x} >0, \ \mbox{ and } \ g^{(6)}(\bar{x}) = 6! c_6 >0.
\]
\iffalse
The details of the calculations are omitted.
\fi
Thus under the
condition $\ell > \ell_c = K''_1(\beta_c)$ sequence 4a satisfies all
the hypotheses of Theorem \ref{thm:main}.

\skp

\noindent {\bf Sequence 5a}

\noindent {\it Definition}. Given $\alpha
> 0$ and a real number $\ell \neq K''(\beta_c)$, the
sequence 5 is defined by \be \label{eqn:seq5} \beta_n = \beta_c -
1/n^\alpha \ \mbox{ and } \ K_n = K(\beta_c) - K'(\beta_c)/n^\alpha
+ \ell/(2n^\alpha). \ee This sequence converges to the tricritical
point $(\beta_c, K(\beta_c))$ from the left along the curve that
coincide with the second-order curve to order 2 in powers of
$\beta-\beta_c$. We assume that $\ell > K''(\beta_c)$. Under this
assumption it is proved in Theorem 4.5 in \cite{EllMacOtt1}
that sequence 5 satisfies the hypotheses of Theorem
\ref{thm:3.1} with $\alpha_0=1/3$ and $\theta=1/2$. When $\ell > K''(\beta_c)$, we refer to sequence 5 as sequence 5a.

\skp \noi {\it Hypothesis {\em (iii$^\prime$)} in Theorem {\em
\ref{thm:main}} for sequence {\em 5a}.} Since $\alpha_0=1/3$ and $\theta=1/2$,
$\theta\alpha_0$ lies in the interval $[1/6,1/2)$ as required by
hypothesis (iii$^\prime$). The limits in (\ref{eqn: hypiii'6}) for
$j=2,3,4,5,6$ are proved in appendix A. Define $y$ as in (\ref{eqn:y}).
\iffalse
\[
y = \left( 1+\frac{3}{5}\beta_c (\ell-K''(\beta_c)) \right)^{1/2}.
\]
\fi
Since $\ell > K''(\beta_c)$, we have $y>1$. We now prove that
$g^{(j)}(\bar{x}) >0$ for $j=2, 3, 4,5,6$ using the
formulas for $g$ and $\bar{x}$ in Theorem 4.5 in \cite{EllMacOtt1}.
Let $c_4=3/16$ and $c_6=9/40$. These formulas yield
%%%%%%%%%%%%%%%%%%%%%%%%%%%%%%%%%%%%%%%%%%%%%%%%%%%%%%%%%%%%%%%%%%%%%%%%%%%%%%%%%%%%%%%%%%%%%%%%%%%%%%%%%%%%%%%%%%%%%%%%%%%
\iffalse
\[ g'(x)=\beta_c (K''(\beta_c)-\ell)x - 16 c_4 x^3 + 6c_6 x^5,
\]
\bea g^{(2)}(\bar{x}) &=& \beta_c (K''(\beta_c)-l) + 48 \times \frac{3}{16} \times \frac{10}{9} [ -1+(1+\frac{3}{5}\beta_c (l-K''(\beta_c)))^{1/2} ] \nonumber \\
&&+ 30 \times \frac{9}{40} \times \frac{100}{81} [ 1+(1+\frac{3}{5}\beta_c(l-K''(\beta_c)))-2( 1+ \frac{3}{5}\beta_c(l-K''(\beta_c)) )^{1/2} ]  \nonumber \\
&=& \beta_c (K''(\beta_c)-l) -10 +10( 1+\frac{3}{5}\beta_c (l-K''(\beta_c)) )^{1/2} \nonumber \\
&&+ \frac{50}{3} +5\beta_c( l-K''(\beta_c) ) - \frac{50}{3}( 1+\frac{3}{5}\beta_c (l-K''(\beta_c))  )^{1/2} \nonumber \\
&=& 4\beta_c(l-K''(\beta_c)) + \frac{20}{3} - \frac{20}{3}
(1+\frac{3}{5}\beta_c (l-K''(\beta_c)))^{1/2} \nonumber. \eea Then
we have \[ 4\beta_c(l-K''(\beta_c))= \frac{20}{3}y^2-\frac{20}{3}.
\]
\fi
\[ g^{(2)}(\bar{x}) = \beta_c (K''(\beta_c)-\ell) + 3 \cdot 4 \cdot 4 c_4 \bar{x}^2 + 5 \cdot 6 c_6 \bar{x}^4 = \frac{20}{3}y(y-1) >0,
\]
\[ g^{(3)}(\bar{x}) = 4! \cdot 4 c_4\bar{x} + 4 \cdot 5 \cdot 6 c_6\bar{x}^3
\iffalse = 9\bar{x}(2+3\bar{x}^2) \fi >0, \ g^{(4)}(\bar{x}) =  4! \cdot 4 c_4 + 3 \cdot 4 \cdot 5 \cdot 6 c_6\cdot \bar{x}^2>0,
\]
\[ \ g^{(5)}(\bar{x}) = 6! c_6 \bar{x} >0,
\ \mbox{ and } \ g^{(6)}(\bar{x}) = 6! c_6  >0.
\]
\iffalse
The details of the calculations are omitted.
\fi
Thus under the
condition $\ell>K''(\beta_c)$ sequence 5a satisfies all the
hypotheses of Theorem \ref{thm:main}.

\skp We have completed the discussion of the verification of the hypotheses of Theorem
\ref{thm:main} for sequences 1a--5a in Table 5.1. This is the
content of part (b) of Theorem \ref{thm:main}. Part (a) of that
theorem is proved in the next section.

%===================== new section =====================================

\section{Proof of Part (a) of Theorem \ref{thm:main}}
\label{section:proofmain}
\beginsec

Theorem \ref{thm:6.1}, a new theorem stated in this section, has two parts.
\iffalse
The proof of part (a) of Theorem \ref{thm:main} depends on Theorem
\ref{thm:6.1}, a new result stated in this section. We start by
sketching the proof of part (a) of Theorem \ref{thm:main} from this new result.
\fi
Under the same hypotheses as Theorem \ref{thm:main},
part (a) of Theorem \ref{thm:6.1} states a conditional
central limit theorem: conditioned on the event $\{S_n/n>\delta
m(\beta_n,K_n)\}$ for $\delta \in (0,1)$ sufficiently close to 1,
the $P_{n, \beta_n, K_n}$-distributions of $n^\kappa
(S_n/n-m(\beta_n,K_n))$ converge weakly to an $N(0,1/g^{(2)}(\bar{x}))$-random variable
with mean 0 and variance $1/g^{(2)}(\bar{x})$. Under the same hypotheses as Theorem \ref{thm:main},
part (b) of Theorem \ref{thm:6.1} states the related conditional limit
\bea
\lefteqn{
\lim_{n \goto \infty} n^\kappa E_{n,\beta_n,K_n} \{
|S_n/n-m(\beta_n,K_n)| \ \big| \ S_n/n>\delta m(\beta_n,K_n) \} }
\nonumber \\
&&= E\{|N(0,1/g^{(2)}(\bar{x}))|\} = (2/(\pi g^{(2)}(\bar{x})))^{1/2} = \bar{z} \nonumber.\eea
\iffalse*********************************************************of
the expected value of $n^\kappa |S_n/n-m(\beta_n,K_n)|$. In Theorem
\ref{thm:6.1}(a) we prove the conditional limit of the expected
value of $f(n^\kappa (S_n/n-m(\beta_n,K_n)))$, where $f$ is any
bounded, continuous function. Proposition
\ref{prop:weakinteg} allows us to replace the bounded,
continuous $f$ by the absolute value function, yielding part (b) of
Theorem \ref{thm:6.1}. (change on April 25)*********************\fi
\iffalse EXPLAIN PHYSICAL INTERPRETATION OF CONDITIONING ON THE
EVENT $\{S_n/n>\delta m(\beta_n,K_n)\}$. \fi

We now sketch the proof of part (a) of Theorem \ref{thm:main} from part (b) of Theorem \ref{thm:6.1}.
In Lemma \ref{lem:Sym1/2Small}, we show that the moderate deviation principle
in Theorem \ref{thm:MDP} and the asymptotic behavior of
$m(\beta_n,K_n)$ in Theorem \ref{thm:3.1} imply that the event $\{
S_n/n > \delta m(\beta_n,K_n) \}$ and the symmetric event $\{ S_n/n
<- \delta m(\beta_n,K_n) \}$ have large probability and that the
event $ \{ \delta m(\beta_n,K_n)> S_n/n
> - \delta m(\beta_n,K_n) \}$ has an exponentially small
probability. As we show at the end of this section, combining part (b) of Theorem \ref{thm:6.1} with Lemma
\ref{lem:Sym1/2Small} and using symmetry yield
\[\lim_{n \goto \infty} n^\kappa E_{n,\beta_n,K_n}
\{||S_n/n|-m(\beta_n,K_n)|\}=\bar{z}.
\]
This is part (a) of Theorem \ref{thm:main}.

The proofs of parts (a) and (b) of Theorem \ref{thm:6.1} are
long and technical. 
\iffalse
--------------Part (a) is proved in subsections
\ref{subsection:step1OFpartA}, \ref{subsection:step2aOFpartA}, and
\ref{subsection:step2bOFpartA}, and
----------------------------------------------------------\fi Part
(b) is proved in subsections \ref{subsection:step1OFpartB},
\ref{subsection:step2aOFpartB}, and \ref{subsection:step2bOFpartB}
using a number of preparatory lemmas in section 7.
At the end of section 7 in \cite{EllLi} we outline 
the proof of part (a), which follows the pattern of proof of part (b) but
is more straightforward.
The weak convergence result proved in Lemma \ref{lem:Thm61aStep1}
is the seed that yields both the conditional central limit theorem in part (a)
of Theorem \ref{thm:6.1} and the conditional limit in part (b) of Theorem \ref{thm:6.1}.

The hypotheses of Theorem \ref{thm:6.1} coincide with the hypotheses of Theorem \ref{thm:main}.
Part (c) of Theorem \ref{thm:6.1} states that for $\alpha \in (0,
\alpha_0)$, $\kappa = \frac{1}{2}(1-\alpha/\alpha_0)+\theta\alpha$
lies in the interval $(\theta\alpha_0, 1/2)$. This fact is needed in
the proofs of Lemmas \ref{lem:tEgammaLowerValue},
\ref{lem:G(xmn)miusG(mn)}, and \ref{lem:Substep2bOFThm61b}. \iffalse
Part (c) of Theorem \ref{thm:6.1} is proved immediately after the
statement of the next theorem. \fi The proof that $\kappa \in (\theta\alpha_0, 1/2)$
is elementary. By hypothesis (iii$^\prime$) of Theorem \ref{thm:main}, we have
$\theta\alpha_0 < 1/2$, which gives $\theta < 1/(2\alpha_0)$.
Therefore
\[\kappa =
\frac{1}{2}(1-\alpha/\alpha_0)+\theta\alpha=\frac{1}{2}+
\alpha(\theta - 1/(2\alpha_0)) < 1/2.
\]
Since $0<\alpha<\alpha_0$ and $\theta < 1/(2\alpha_0)$, we have
$\kappa > \frac{1}{2}+\alpha_0(\theta - 1/(2\alpha_0))=\theta\alpha_0$.
This completes the proof of part (c) of Theorem \ref{thm:6.1}.

Concerning part (d) of Theorem \ref{thm:6.1},
the hypotheses of this theorem coincide with the hypotheses of
Theorem \ref{thm:main}. Thus, as shown in section 5 and appendix A, 
these hypotheses are satisfied by sequences
1a--5a in Table 5.1.
\iffalse
Since the hypotheses of Theorem \ref{thm:6.1} and Theorem \ref{thm:main} coincide,
these hypotheses are satisfied by sequences 1a--5a in Table 5.1,
as shown in section 5 and appendix A.
This is the assertion of part (d) of Theorem \ref{thm:6.1}.
\fi

\iffalse We then use part (b) of Theorem \ref{thm:6.1} and Lemma
\ref{lem:Sym1/2Small} to prove Theorem \ref{thm:main}. Theorem
\ref{thm:6.1} is also true if hypothesis (iii)(b) is replaced by the
following hypothesis: the set of global minimum points of the
Ginzburg-Landau polynomial $g$ equals $\{\pm \bar{x}\}$ and $\{0,
\pm \bar{x}\}$. But we do not need it in this dissertation. ADD MORE
ADD MORE ADD MORE. (change on April 25)\fi

\begin{thm}
\label{thm:6.1} Let $(\beta_n,K_n)$ be a positive sequence
converging either to a second-order point $(\beta, K(\beta))$,
$0<\beta<\beta_c$, or to the tricritical point $(\beta,
K(\beta))=(\beta_c, K(\beta_c))$. We assume that for all
$0<\alpha<\alpha_0$, $(\beta_n, K_n)$ satisfies the hypotheses of
Theorem {\em \ref{thm:main}}, which coincide with the hypotheses of
Theorem {\em \ref{thm:3.1}} together with hypothesis {\em
(iii$^\prime$)}. \iffalse************************** with the
exception that hypothesis {\em (iii)(b)} of Theorem {\em
\ref{thm:3.1}} is replaced by the following hypothesis, which
expends hypothesis {\em (iii)(b)} of Theorem {\em \ref{thm:3.1}}
\begin{itemize}
\item[(iii)(b$^\prime$)] The set of global minimum points of the
Ginzburg-Landau polynomial $g$ equals $\{\pm \bar{x}\}$ or $\{0, \pm
\bar{x}\}$.
\end{itemize}  ******************************************\fi
For $\alpha \in (0, \alpha_0)$ we define $\kappa =
\frac{1}{2}(1-\alpha/\alpha_0)+\theta\alpha$. Then for any
$0<\alpha<\alpha_0$ there exists $\Delta \in (0,1)$ such that for
any $\delta \in (\Delta, 1)$ the following conclusions hold.

{\em (a)} When conditioned on the event $\{ S_n/n > \delta m(\beta_n, K_n)
\}$, the $P_{n, \beta_n,K_n}$-distributions of
$n^\kappa(S_n/n-m(\beta_n,K_n))$ converge weakly to a normal random
variable $N(0,1/g^{(2)}(\bar{x}))$ with mean $0$ and
variance $1/g^{(2)}(\bar{x})$; in symbols,
\[ P_{n, \beta_n,K_n}\{n^\kappa(S_n/n-m(\beta_n,K_n)) \in dx \  \big| \ S_n/n>\delta m(\beta_n,K_n)  \}
\Longrightarrow N(0,1/g^{(2)}(\bar{x})) .
\]
\iffalse The limiting density is $
\exp[-\frac{1}{2}g^{(2)}(\bar{x})x^2]/
\int_{\mathbb{R}}\exp[-\frac{1}{2}g^{(2)}(\bar{x})x^2]dx$. \fi

{\em (b)} We have the conditional limit
\bea
\lefteqn{ \lim_{n \goto \infty} n^\kappa E_{n,\beta_n,K_n} \{
|S_n/n-m(\beta_n,K_n)|  \  \big| \  S_n/n>\delta m(\beta_n,K_n) \} }
\nonumber \\
&&=
\lim_{n \goto \infty}  E_{n,\beta_n,K_n} \{ |S_n/n^{1-\kappa} -
n^\kappa m(\beta_n,K_n)|  \  \big| \  S_n/n>\delta m(\beta_n,K_n) \}
=\bar{z} \nonumber,
\eea
where
\beas
\bar{z}
&=&
E \{|N(0,1/g^{(2)}(\bar{x}))|\} \\
&=&
\frac{1}{\int_{\mathbb{R}}\exp[-\frac{1}{2}
g^{(2)}(\bar{x})x^2]dx } \cdot \int_{\mathbb{R}} |x| \exp[\ts
-\frac{1}{2} g^{(2)}(\bar{x})x^2]dx = \ds \left(\frac{2}{ \pi g^{(2)}(\bar{x})}\right)^{1/2}.
\eeas

{\em (c)} For $\alpha \in (0, \alpha_0)$, $\kappa =
\frac{1}{2}(1-\alpha/\alpha_0)+\theta\alpha$ lies in the interval
$(\theta\alpha_0, 1/2)$.

{\em (d)} The hypotheses of this theorem are satisfied by sequences {\em{1a--5a}}
in Table {\em 5.1}.
\end{thm}

In part (a) of Theorem \ref{thm:MDP} we state a moderate deviation principle (MDP) for the mean-field B-C model.
This MDP will be used to prove Lemma \ref{lem:Sym1/2Small},
which in turn will be used to prove part (a) of Theorem
\ref{thm:main} from part (b) of Theorem \ref{thm:6.1}. The rate function in the MDP is the continuous
function $\Gamma(x)=g(x)-\inf_{y \in \mathbb{R}}g(y)$, where $g$ is the associated Ginzburg-Landau polynomial.
$\Gamma$ satisfies $\Gamma(x) \goto \infty$ as $|x| \goto \infty$. For $A$ a subset of $\R$
define $\Gamma(A)= \inf_{x \in A} \Gamma(x)$.

\begin{thm}
\label{thm:MDP} Let $(\beta_n,K_n)$ be a positive sequence
converging either to a second-order
point $(\beta, K(\beta))$, $0<\beta<\beta_c$, or to the tricritical
point $(\beta, K(\beta))=(\beta_c, K(\beta_c))$. We assume that
$(\beta_n, K_n)$ satisfies the hypotheses of Theorem
{\em\ref{thm:main}} for all $0<\alpha<\alpha_0$. The following
conclusions hold.

 % \vspace{-.1in}
{\em (a)} For all $0<\alpha<\alpha_0$,
  $S_n/n^{1-\theta\alpha}$ satisfies the MDP with respect to $P_{n,
  \beta_n,K_n}$ with exponential speed $n^{1-\alpha/\alpha_0}$ and
  rate function $\Gamma(x)=g(x)-\inf_{y \in \mathbb{R}}g(y)$;
\iffalse
  in symbols
  \[ P_{n, \beta_n, K_n} \{S_n/n^{1-\theta\alpha}\in dx\} \asymp
  \exp[-n^{1-\alpha/\alpha_0} \Gamma(x)]dx.
  \]
\fi
i.e.,
for any closed set $F$ in $\mathbb{R}$
\[
\iffalse \label{eqn:MDPupperbound}  \fi
\limsup_{n\goto\infty} \frac{1}{n^{1 - \alpha/\alpha_0}} \log
P_{n,\beta_n,K_n}\{S_n/n^{1-\theta\alpha} \in F\} \le -\Gamma(F)
\]
and for any open set $G$ in $\mathbb{R}$
\[
\iffalse \label{eqn:MDPlowerbound} \fi
\liminf_{n\goto\infty} \frac{1}{n^{1 - \alpha/\alpha_0}} \log
P_{n,\beta_n,K_n}\{S_n/n^{1-\theta\alpha} \in G\} \ge -\Gamma(G). \]

{\em (b)} The hypotheses of this theorem are satisfied by
  sequences {\em{1a--5a}} in Table {\em 5.1} .
\end{thm}

The MDP in part (a) of Theorem \ref{thm:MDP} is proved like
the MDP in part (a) of Theorem 8.1 in \cite{CosEllOtt} with only
changes in notation. Because of the importance of the MDP in part (a) of Theorem \ref{thm:MDP},
the proof is given in appendix B. Concerning part (b) of Theorem \ref{thm:MDP},
the hypotheses of this theorem coincide with the hypotheses of
Theorem \ref{thm:main}. Thus, as shown in section 5 of this paper and in appendix A, these hypotheses are satisfied by sequences
1a--5a in Table 5.1.

After proving the next lemma, we use it to derive part (a) of Theorem
\ref{thm:main} from part (b) of Theorem \ref{thm:6.1}.
%****add on April 28 end******

\begin{lem}
\label{lem:Sym1/2Small} We assume that $(\beta_n, K_n)$ satisfies
the hypotheses of Theorem {\em \ref{thm:main}} for all
$0<\alpha<\alpha_0$. Then for any $0<\alpha<\alpha_0$ and any
$\delta \in (0,1)$ there exists $c>0$ such that for all sufficiently
large n
\[P_{n,\beta_n,K_n} \{ \delta m(\beta_n,K_n) \ge S_n/n \ge -\delta m(\beta_n,K_n) \} \le \exp[-c n^{1-\alpha/\alpha_0}] \goto
0  \ \ \mbox{as} \ \ n \goto \infty.
\]
In addition,
\[\lim_{n \goto \infty} P_{n,\beta_n,K_n} \{S_n/n >\delta m(\beta_n,K_n)
\} = \lim_{n \goto \infty} P_{n,\beta_n,K_n} \{S_n/n <-\delta
m(\beta_n,K_n) \} =1/2.
\]
\end{lem}

\noindent{\bf Proof.}
\iffalse
By Theorem \ref{thm:MDP},
$S_n/n^{1-\theta\alpha}$ satisfies the MDP with exponential speed
$n^{1-\alpha/\alpha_0}$ and rate function $\Gamma(x)=g(x)-\inf_{y
\in \mathbb{R}} g(y)$, where $g$ is the associated Ginzburg-Landau
polynomial.
\fi
By hypothesis (iii)(b) of Theorem \ref{thm:3.1}
the global minimum points of $g$ are $\pm \bar{x}$, and by Theorem
\ref{thm:3.1}, $n^{\theta\alpha}m(\beta_n,K_n) \goto \bar{x}$ as $ n
\goto \infty$. Thus we can choose $\varepsilon
>0$ satisfying $(1+\varepsilon)\delta <1$ such that $n^{\theta\alpha}m(\beta_n,K_n) \le
(1+\varepsilon)\bar{x}$ for large $n$. Let $F$ be the closed set
$[-(1+\varepsilon) \delta \bar{x}, \ (1+\varepsilon) \delta \bar{x}]
$. Since $(1+\varepsilon) \delta \bar{x} < \bar{x}$ and
$-(1+\varepsilon) \delta \bar{x} > - \bar{x}$, we have
\[\inf_{y \in F} g(y) > \inf_{z \in \mathbb{R}} g(z) =
g(\bar{x}),
\]
which implies
\[ \Gamma(F) = \inf_{y \in F} \{ g(y)- \inf_{z \in
\mathbb{R}}g(z) \}
>0.
\]
\iffalse Choose the constant $c$ to be any number satisfying
$0<c<\inf_{|x|\le (1+\varepsilon) \delta \bar{x}} g(x)$. \fi We
write $m_n=m(\beta_n, K_n)$. The moderate deviation upper bound
in part (a) of Theorem \ref{thm:MDP} implies that for all sufficiently large
$n$
\bea
\lefteqn{ P_{n,\beta_n,K_n}\{ \delta m_n \ge S_n/n \ge
-\delta m_n \}
 }\nonumber \\
&&=
P_{n,\beta_n,K_n}\{ \delta n^{\theta\alpha} m_n \ge
S_n/n^{1-\theta\alpha} \ge
-\delta n^{\theta\alpha} m_n \} \nonumber \\
&&\le
P_{n,\beta_n,K_n}\{ (1+\varepsilon)\delta\bar{x} \ge
S_n/n^{1-\theta\alpha} \ge -(1+\varepsilon)\delta\bar{x} \}
\nonumber \\
&& \le
\exp[-n^{1-\alpha/\alpha_0}\Gamma(F)/2] \goto 0 \ \mbox{ as } n \goto \infty \nonumber.
\eea
\iffalse
It follows that
\[ P_{n,\beta_n,K_n}\{
\delta m_n \ge S_n/n \ge -\delta m_n \} \le
\exp[-n^{1-\alpha/\alpha_0}\Gamma(F)/2]  \goto 0   \ \ \mbox{as}  \
\ n \goto \infty.
\]
\fi
This yields the first assertion in the lemma.

To prove the second assertion, we write
\bea
1=
P_{n,\beta_n,K_n}\{ S_n/n \in \mathbb{R}
\}
&=&
P_{n,\beta_n,K_n} \{ \delta m_n \ge S_n/n \ge -\delta m_n \}
\nonumber \\
&&+
\ P_{n,\beta_n,K_n} \{ S_n/n > \delta m_n \} \nonumber
\\
&&+
\ P_{n,\beta_n,K_n} \{ S_n/n <- \delta m_n \}\nonumber,
\eea
Symmetry and the first assertion imply that
\[\lim_{n \goto \infty}P_{n,\beta_n,K_n} \{ S_n/n > \delta m_n \}= \lim_{n \goto \infty}P_{n,\beta_n,K_n} \{ S_n/n <- \delta m_n
\}=1/2.
\]
This completes the proof of the lemma. $\Box$

\skp Now we are ready to prove part (a) of Theorem \ref{thm:main}.

\skp \noindent{\bf Proof of part (a) of Theorem \ref{thm:main} from
part (b) of Theorem \ref{thm:6.1} and Lemma \ref{lem:Sym1/2Small}.}
We write $m_n=m(\beta_n, K_n)$. Define
\[ p_{n, \delta}^+ = P_{n,\beta_n,K_n} \{ S_n/n > \delta m_n \}, \
p_{n, \delta}^- = P_{n,\beta_n,K_n} \{ S_n/n <- \delta m_n \},
\]
\[ q_{n,\delta} = P_{n,\beta_n,K_n}\{ \delta m_n \ge S_n/n \ge -\delta m_n
\}.
\]
Since by symmetry $p_{n, \delta}^+ = p_{n, \delta}^-$  and
\bea
\lefteqn{ E_{n,\beta_n,K_n} \{ ||S_n/n|-m_n|
\ \big| \ S_n/n < -\delta m_n \}  }\nonumber \\
&&=
E_{n,\beta_n,K_n} \{ ||S_n/n|-m_n| \ \big| \ S_n/n > \delta m_n
\}  \nonumber \\
&&=
E_{n,\beta_n,K_n} \{ |S_n/n-m_n| \ \big| \ S_n/n > \delta m_n \}
\nonumber,
\eea
we have
\bea
E_{n,\beta_n,K_n} \{ ||S_n/n|-m_n| \}
&=&
E_{n,\beta_n,K_n} \{ ||S_n/n|-m_n|
\ \big| \ S_n/n > \delta m_n \} \cdot p_{n, \delta}^+  \nonumber \\
&& +
\ E_{n,\beta_n,K_n} \{ ||S_n/n|-m_n| \ \big| \ S_n/n < -\delta
m_n \} \cdot p_{n, \delta}^- \nonumber \\
&& +
\ E_{n,\beta_n,K_n} \{ ||S_n/n|-m_n| \ \big| \ \delta m_n \ge
S_n/n \ge -\delta m_n \} \cdot q_{n, \delta} \nonumber \\
&=&
2 \cdot E_{n,\beta_n,K_n} \{ |S_n/n-m_n|
\ \big| \ S_n/n > \delta m_n \} \cdot p_{n, \delta}^+  \nonumber \\
&& +
\ E_{n,\beta_n,K_n} \{ ||S_n/n|-m_n| \ \big| \ \delta m_n \ge
S_n/n \ge -\delta m_n \} \cdot q_{n, \delta} \nonumber.
\eea
By part
(b) of Theorem \ref{thm:6.1} and Lemma \ref{lem:Sym1/2Small}
\[ \lim_{n \goto \infty} n^\kappa E_{n,\beta_n,K_n} \{ |S_n/n-m_n|
\ \big| \ S_n/n > \delta m_n \} \cdot p_{n, \delta}^+ = \frac{1}{2}
\bar{z}.
\]
Since $|S_n/n| \le 1$ and $0 \le m_n \le 1$, Lemma
\ref{lem:Sym1/2Small} implies that there exists $c>0$ such that for
all sufficiently large $n$
\bea
\lefteqn{ n^\kappa E_{n,\beta_n,K_n}
\{ ||S_n/n|-m_n| \ \big| \ \delta m_n \ge S_n/n \ge -\delta m_n \}
\cdot q_{n, \delta}    } \nonumber \\
&&\le
2 n^\kappa q_{n, \delta}
\le
2n^\kappa
\exp[-cn^{1-\alpha/\alpha_0}] \goto 0  \  \  \mbox{as} \  \  n \goto
\infty \nonumber.
\eea
It follows that
\[\lim_{n \goto \infty} n^\kappa E_{n,\beta_n,K_n} \{ ||S_n/n|-m_n|
\} = \frac{1}{2}\bar{z} + \frac{1}{2}\bar{z} = \bar{z}.
\]
Part (a) of Theorem \ref{thm:main} is proved. $\Box$

%%%%%%%  %%%%%%%%%%%  %%%%%%%%%

\iffalse \skp In section \ref{section:verify6}, we verify the
hypotheses of Theorem \ref{thm:main} for sequences 1a--5a. This
gives part (b) of Theorem \ref{thm:main}. Hence the proof of Theorem
\ref{thm:main} is complete. \ \ $\Box$  (change on April 25)\fi

\skp In the next section we prove a number of lemmas that will be used in section 8 to prove part (b) of Theorem
\ref{thm:6.1}.

%===================== new section =====================================

\section{Preparatory Lemmas for Proof of Part (b) of Theorem \ref{thm:6.1}}
\label{section:LogicANDPreLemmas}
\beginsec

Let $(\beta_n, K_n)$ be a positive sequence. Throughout this section we work with $0 < \alpha < \alpha_0$ and denote
$m(\beta_n,K_n)$ by $m_n$. Let $W_n$ be a sequence of normal random
variables with mean 0 and variance $(2\beta_n K_n)^{-1}$ defined on
a probability space $(\Omega,\mathcal {F},Q)$. We denote by
$\tilde{E}_{n, \beta_n,K_n}$ expectation with respect to the product
measure $P_{n,\beta_n,K_n} \times Q$; $P_{n,\beta_n,K_n}$ is defined in (\ref{eqn:P})--(\ref{eqn:PnbetaK}). Because the proof of part (b)
of Theorem \ref{thm:6.1} is long and technical, we start by
explaining the logic. The hypotheses of this theorem coincide with those
of Theorem \ref{thm:main}.

Part (b) of Theorem \ref{thm:6.1} states that there exists $\Delta \in (0,1)$ such that for any $\delta
\in (\Delta, 1)$ \bea \label{eqn:Thm61b} \lefteqn{ \lim_{n \goto
\infty} n^{\kappa}E_{n,\beta_n,K_n} \{ |S_n/n-m_n| \ \big| \
S_n/n>\delta m_n \} } \\
&&= \lim_{n \goto \infty}  E_{n,\beta_n,K_n} \{ |S_n/n^{1-\kappa}-
n^{\kappa}m_n| \ \big| \  S_n/n>\delta m_n
\} \nonumber \\
&&= \frac{1}{\int_{\mathbb{R}}\exp[-\frac{1}{2}
g^{(2)}(\bar{x})x^2]dx } \cdot \int_{\mathbb{R}} |x| \exp[\ts
-\frac{1}{2} g^{(2)}(\bar{x})x^2]dx = \bar{z}.\nonumber \eea $\Delta
\in (0,1)$ is determined in Lemma \ref{lem:G(xmn)miusG(mn)}. The key
idea in proving (\ref{eqn:Thm61b}) is to show that adding suitably
scaled versions of the normal random variables $W_n$
yields a quantity with the following two properties: its limit equals the last line of (\ref{eqn:Thm61b})
and the second line of (\ref{eqn:Thm61b}) has the same limit; specifically,
\bea \label{eqn:Thm61bKey} \lefteqn{
\lim_{n \goto \infty} E_{n,\beta_n,K_n} \{ |S_n/n^{1-\kappa}-
n^{\kappa}m_n| \ \big| \ S_n/n>\delta m_n
\} } \\
&&= \lim_{n \goto \infty}  \tilde{E}_{n,\beta_n,K_n} \{
|S_n/n^{1-\kappa}+W_n/n^{1/2-\kappa}- n^{\kappa}m_n| \ \big| \
S_n/n+W_n/n^{1/2}>\delta m_n
\} \nonumber \\
&&= \frac{1}{\int_{\mathbb{R}}\exp[-\frac{1}{2}
g^{(2)}(\bar{x})x^2]dx } \cdot \int_{\mathbb{R}} |x| \exp[\ts
-\frac{1}{2} g^{(2)}(\bar{x})x^2]dx.\nonumber \eea

Formula (\ref{eqn:Thm61bKey}) is proved in two steps.

\skp \noi {\bf Step 1.} Prove the second limit in
(\ref{eqn:Thm61bKey}). This is done in part (b) of Lemma
\ref{lem:Step1OFThm61b} in subsection \ref{subsection:step1OFpartB}.

\skp \noi {\bf Step 2.} Prove the first limit in
(\ref{eqn:Thm61bKey}). This is done in two substeps, which we now
explain.

\skp \noi {\bf Substep 2a.} Define
\[ C_n = \tilde{E}_{n,\beta_n,K_n} \{ |S_n/n^{1-\kappa}- n^{\kappa}m_n| \
\big| \  S_n/n>\delta m_n \}
\]
and
\[ D_n = \tilde{E}_{n,\beta_n,K_n} \{
|S_n/n^{1-\kappa}+W_n/n^{1/2-\kappa}- n^{\kappa}m_n| \ \big| \
S_n/n>\delta m_n \}.
\]
Thus $D_n$ is obtained from $C_n$ by replacing $S_n/n^{1-\kappa}$ by
$S_n/n^{1-\kappa}+W_n/n^{1/2-\kappa}$. Substep 2a is to prove that
$\lim_{n \goto \infty} | C_n - D_n |=0$.
This is done in Lemma \ref{lem:substep2aOFThm61b} in subsection
\ref{subsection:step2aOFpartB}.

\skp \noi {\bf Substep 2b.} Define
\[ F_n = \tilde{E}_{n,\beta_n,K_n} \{
|S_n/n^{1-\kappa}+W_n/n^{1/2-\kappa}- n^{\kappa}m_n| \ \big| \
S_n/n+W_n/n^{1/2}>\delta m_n \}.
\]
Thus $F_n$ is obtained from $D_n$ by replacing $S_n/n$ in the
conditioned event $\{ S_n/n > \delta m_n\}$ by $S_n/n+W_n/n^{1/2}$.
Substep 2b is to prove that
\[
\lim_{n \goto \infty} D_n = \lim_{n \goto \infty} F_n =
\frac{1}{\int_{\mathbb{R}}\exp[-\frac{1}{2} g^{(2)}(\bar{x})x^2]dx }
\cdot \int_{\mathbb{R}} |x| \exp[ \ts -\frac{1}{2}
g^{(2)}(\bar{x})x^2]dx.
\]
The limit of $F_n$ as $n \goto \infty$
is calculated in Step 1. Substep 2b is proved in part (b) of Lemma
\ref{lem:Substep2bOFThm61b} in subsection
\ref{subsection:step2bOFpartB}. The explanation of the logic of the proof of part (b)
of Theorem \ref{thm:6.1} is now complete.

\skp
We next state and prove the preparatory lemmas needed to carry out
Step 1, Substep 2a, and Substep 2b in the proof of part (b) of Theorem \ref{thm:6.1}.

\iffalse
\skp The hypotheses of Theorem \ref{thm:6.1} coincide with the
hypotheses of Theorem \ref{thm:main}. Step 1 in the proof of part
(b) of Theorem \ref{thm:6.1} is to prove that
\bea \label{eqn:Step1}
\lefteqn {\lim_{n \goto \infty} F_n} \\
&=& \lim_{n \goto \infty} \tilde{E}_{n,\beta_n,K_n} \{
|S_n/n^{1-\kappa}+W_n/n^{1/2-\kappa}- n^{\kappa}m_n| \ \big| \
S_n/n+W_n/n^{1/2}>\delta m_n
\} \nonumber \\
&=& \frac{1}{\int_{\mathbb{R}}\exp[ -\frac{1}{2}
g^{(2)}(\bar{x})x^2]dx } \cdot \int_{\mathbb{R}} |x| \exp [\ts
-\frac{1}{2} g^{(2)}(\bar{x})x^2 ] dx.\nonumber \eea
This is proved
in part (b) of Lemma \ref{lem:Step1OFThm61b}. In order to do this
proof, we need Lemmas \ref{lem:RepresentF}--\ref{lem:intRnkexpIneq}
and \ref{lem:Thm61aStep1}(a), several of which will also be used
elsewhere. \iffalse in other steps. \fi
\fi

\iffalse \skp LEMMAS 7.1-7.7 HERE ARE THE SAME CONTENT AND ORDER AS
LEMMA 7.1-7.7 IN THESIS. OF COURSE I MAKE SOME SMALL CHANGES INSIDE
SEVERAL LEMMAS. \fi

Lemma \ref{lem:RepresentF} is a representation formula that
will be used to study the limit of the conditional expectation in the second line of
(\ref{eqn:Thm61bKey}). This lemma can be proved like Lemma 3.3 in
\cite{EllNew}, which applies to the Curie-Weiss model, or like Lemma
3.2 in \cite{EllWan}, which applies to the Curie-Weiss-Potts model.
It will also be used to prove Lemma
\ref{lem:tEgammaLowerValue} and Lemma \ref{lem:intRnkexpIneq}.

\begin{lem}
\label{lem:RepresentF} Given a positive sequence $(\beta_n,K_n)$,
let $W_n$ be a sequence of normal random variables with mean 0 and
variance $(2\beta_n K_n)^{-1}$ defined on a probability space
$(\Omega,\mathcal {F},Q)$. Then for any $\bar{\gamma} \in [0,1)$ and
any bounded, measurable function $\varphi$
\bea
\lefteqn{ \int_{\Lambda^n \times \Omega} \varphi(S_n/n^{1-\bar{\gamma}}+W_n/n^{1/2-\bar{\gamma}})d(P_{n,\beta_n,K_n}\times Q) } \nonumber \\
&&=
\frac{1}{\int_{\mathbb{R}} \exp[-n
G_{\beta_n,K_n}(x/n^{\bar{\gamma}})]dx} \cdot \int_{\mathbb{R}}
\varphi(x)\exp[-n G_{\beta_n,K_n}(x/n^{\bar{\gamma}})]dx \nonumber.
\eea
In this formula $G_{\beta_n,K_n}$ is the free energy function defined in {\em (\ref{eqn:GbetaK})}.
\end{lem}

Lemma \ref{lem:tEgammaLowerValue} uses the representation
formula in the preceding lemma to rewrite the conditional
expectation in the second line of (\ref{eqn:Thm61bKey}). \iffalse Given a positive sequence
$(\beta_n,K_n)$, let $W_n$ be a sequence of normal random variables
with mean 0 and variance $(2\beta_n K_n)^{-1}$ defined on a
probability space $(\Omega,\mathcal {F},Q)$. We denote by
$\tilde{E}_{n, \beta_n,K_n}$ the expectation with respect to the
product measure $P_{n,\beta_n,K_n} \times Q$. \fi

\iffalse Lemma \ref{lem:tEgammaLowerValue} is very useful. We apply
it to $h=f$ and to $\bar{\gamma} = \kappa$ in the proof of part (b)
of Lemma \ref{lem:Thm61aStep1}.  We also apply part (b) of this
lemma to $h(x)=|x|$ and to $\bar{\gamma} = \kappa$ in the proof of
Lemma \ref{lem:Step1OFThm61b}. We apply (\ref{eqn:tE1An}) in this
lemma with $\bar{\gamma} = \kappa$ to prove substep 2b of Theorem
\ref{thm:6.1}(a)(i.e. Lemma \ref{lem:Substep2bOFThm61a}) and substep
2b of Theorem \ref{thm:6.1}(b)(i.e. Lemma
\ref{lem:Substep2bOFThm61b}).  \fi

\iffalse \skp CHANGE IN PART (B): APPLYING TO $h(x) = |x|$, WHICH
WILL BE USED IN THE PROOF OF PART (B) OF LEMMA 8.1. \fi

\begin{lem}
\label{lem:tEgammaLowerValue} We assume that $(\beta_n, K_n)$
satisfies the hypotheses of Theorem {\em{\ref{thm:main}}} for all
$0<\alpha<\alpha_0$. For any
$\bar{\delta} \in (0,1)$ define
\[ A_n(\bar{\delta}) = \{ S_n/n + W_n/n^{1/2} > \bar{\delta} m_n \},
\]
where $m_n=m(\beta_n, K_n)$. Given any $\alpha \in (0, \alpha_0)$, define
$\kappa = \frac{1}{2}(1-\alpha/\alpha_0)+\theta\alpha$. The following conclusions hold.

{\em (a)} For any bounded, measurable function $h$
\bea
\label{eqn:tEhmeasurable}
\lefteqn{
\tilde{E}_{n,\beta_n,K_n} \{ h(S_n/n^{1-\kappa}+ W_n/n^{1/2-\kappa}
-n^{\kappa} m_n )
\cdot 1_{A_n(\bar{\delta})} \}   } \\
&&=
\frac{1}{Z_{n, \kappa}}
\int_{-n^{\kappa}(1-\bar{\delta})m_n}^{\infty} h(x) \exp[-n
G_{\beta_n,K_n}(x/n^{\kappa}+m_n)]dx \nonumber,
\eea
where $Z_{n,\kappa} = \int_{\mathbb{R}} \exp[-n
G_{\beta_n,K_n}(x/n^{\kappa})]dx$.
In particular, if $h \equiv 1$, then
\bea \label{eqn:tE1An}
\tilde{E}_{n,\beta_n,K_n} \{ 1_{A_n(\bar{\delta})} \} &=&
(P_{n,\beta_n,K_n} \times Q) \{ A_n(\bar{\delta})  \}  \\
&=&
\frac{1}{Z_{n, \kappa}}
\int_{-n^{\kappa}(1-\bar{\delta})m_n}^{\infty} \exp[-n
G_{\beta_n,K_n}(x/n^{\kappa}+m_n)]dx .\nonumber
\eea

{\em (b)} We have the representation
\bea \label{eqn:tEh|x|} \lefteqn{ \tilde{E}_{n,\beta_n,K_n} \{
|S_n/n^{1-\kappa}+ W_n/n^{1/2-\kappa} -n^{\kappa} m_n |
\big| A_n(\bar{\delta})\}   } \\
&=&
\frac{1}{\int_{-n^{\kappa}(1-\bar{\delta})m_n}^{\infty} \exp[-n
G_{\beta_n,K_n}(x/n^{\kappa}+m_n)+nG_{\beta_n, K_n}(m_n)]dx}
\nonumber
\\
&& \cdot \int_{-n^{\kappa}(1-\bar{\delta})m_n}^{\infty} |x| \exp[-n
G_{\beta_n,K_n}(x/n^{\kappa}+m_n)+nG_{\beta_n, K_n}(m_n)]dx
\nonumber. \eea
\end{lem}

{\noi}{\bf Proof.} (a) \iffalse According to part (c) of Theorem
\ref{thm:6.1}, $0<\kappa<1/2<1$. If $h$ is bounded and measurable,
then we \fi
By part (c) of Theorem \ref{thm:6.1}, $\kappa \in (\theta\alpha_0, 1/2)$.
We apply Lemma \ref{lem:RepresentF} with
$\varphi(x)=h(x-n^{\kappa}m_n) \cdot 1_{(n^{\kappa} \bar{\delta}m_n,
\infty)}(x)$ and $\bar{\gamma} = \kappa$, obtaining
\bea
\lefteqn{ \tilde{E}_{n,\beta_n,K_n}\{
h(S_n/n^{1-\kappa}+W_n/n^{1/2-\kappa} - n^{\kappa}
 m_n) \cdot 1_{ \{
S_n/n^{1-\kappa}+W_n/n^{1/2-\kappa} > n^{\kappa}
\bar{\delta} m_n  \}  } \} }  \nonumber \\
&&=
\int_{\Lambda^n \times \Omega}
h(S_n/n^{1-\kappa}+W_n/n^{1/2-\kappa} - n^{\kappa}
 m_n) \cdot 1_{ \{
S_n/n^{1-\kappa}+W_n/n^{1/2-\kappa} > n^{\kappa}
\bar{\delta} m_n  \}  } d(P_{n,\beta_n,K_n} \times Q) \nonumber \\
\iffalse
&&=
\int_{\Lambda^n \times \Omega}
\varphi (S_n/n^{1-\kappa}+W_n/n^{1/2-\kappa} ) d(P_{n,\beta_n,K_n} \times Q) \nonumber \\
&&=
\frac{1}{\int_{\mathbb{R}} \exp[-n
G_{\beta_n,K_n}(x/n^{\kappa})]dx} \cdot \int_{\mathbb{R}}
\varphi(x)\exp[-n G_{\beta_n,K_n}(x/n^{\kappa})]dx \nonumber \\
\fi
&&=
\frac{1}{\int_{\mathbb{R}} \exp[-n
G_{\beta_n,K_n}(x/n^{\kappa})]dx} \cdot \int_{\mathbb{R}}
h(x-n^{\kappa}m_n) \cdot 1_{( n^{\kappa}\bar{\delta}m_n ,\infty)}(x)
\exp[-n G_{\beta_n,K_n}(x/n^{\kappa})]dx \nonumber \\
&&=
\frac{1}{Z_{n, \kappa}} \cdot
\int_{-n^{\kappa}(1-\bar{\delta})m_n}^{\infty} h(x)  \exp[-n
G_{\beta_n,K_n}(x/n^{\kappa} + m_n)]dx \nonumber.
\eea
\iffalse
The change
of variables $y=x-n^{\kappa}m_n$ in the last integral gives
\bea
\lefteqn{ \hspace{-.5in}
  \frac{1}{\int_{\mathbb{R}} \exp[-n
G_{\beta_n,K_n}(y/n^{\kappa})]dx} \cdot
\int_{-n^{\kappa}(1-\bar{\delta})m_n}^{\infty} h(y)  \exp[-n
G_{\beta_n,K_n}(y/n^{\kappa} + m_n)]dy     } \nonumber \\
&=&
\frac{1}{Z_{n, \kappa}} \cdot
\int_{-n^{\kappa}(1-\bar{\delta})m_n}^{\infty} h(y)  \exp[-n
G_{\beta_n,K_n}(y/n^{\kappa} + m_n)]dy \nonumber.
\eea
\fi
This yields (\ref{eqn:tEhmeasurable}). Formula (\ref{eqn:tE1An}) follows by taking $h \equiv 1$.

(b) We apply part (a) to the sequence of bounded, measurable functions $h_j(x) = |x| \wedge j$, $j \in \N$.
By the monotone convergence theorem we obtain (\ref{eqn:tEhmeasurable}) with $h(x)$ replaced by $|x|$. Part (b)
now follows by using the definition of conditional expectation and multiplying the numerator and denominator
of the resulting fraction by $\exp[nG_{\beta_n, K_n}(m_n)]$.  The proof of Lemma \ref{lem:tEgammaLowerValue}
is complete. $\Box$
\iffalse By part
(a) applied to $h=f$ \bea \lefteqn{ \tilde{E}_{n,\beta_n,K_n} \{
f(S_n/n^{1-\kappa}+ W_n/n^{1/2-\kappa} -n^{\kappa} m_n )
\big| A_n(\bar{\delta})\}   } \nonumber \\
&=&
\frac{1}{P_{n,\beta_n,K_n}\times Q \{ A_n(\bar{\delta}) \}}
\cdot \tilde{E}_{n,\beta_n,K_n} \{ f(S_n/n^{1-\kappa}+
W_n/n^{1/2-\kappa} -n^{\kappa} m_n ) \cdot
1_{A_n(\bar{\delta})} \}   \nonumber \\
&=&
\frac{1}{\int_{-n^{\kappa}(1-\bar{\delta})m_n}^{\infty} \exp[-
\hspace{0.04 in} n G_{\beta_n,K_n}(x/n^{\kappa}+m_n)]dx} \nonumber
\\
&& \cdot \int_{-n^{\kappa}(1-\bar{\delta})m_n}^{\infty} f(x) \exp[-
\hspace{0.04 in} n G_{\beta_n,K_n}(x/n^{\kappa}+m_n)]dx \nonumber.
\eea Multiplying numerator and denominator by $\exp [n G_{\beta_n,
K_n}(m_n)]$ completes the proof for $h=f$. Define $\Pi(x)=|x|$ and
apply previous proof to bounded, continuous function $f_j =
\min(\Pi,j)$ for $j \in \mathbb{N}$. Since $f_j \ge 0$ and $f_j
\uparrow \Pi$, the monotone convergence theorem implies \bea
\lefteqn{ \tilde{E}_{n,\beta_n,K_n} \{ |S_n/n^{1-\kappa}+
W_n/n^{1/2-\kappa} -n^{\kappa} m_n | \ \big| \ A_n(\bar{\delta})\}
        } \nonumber \\
&&=  \frac{\int_{-n^\kappa(1-\bar{\delta})m_n}^{\infty}
|x|\exp[-nG_{\beta_n,K_n}(x/n^\kappa+m_n)+ nG_{\beta_n,K_n}(m_n)]dx
          }
          {\int_{-n^\kappa(1-\bar{\delta})m_n}^{\infty}
\exp[-nG_{\beta_n,K_n}(x/n^\kappa+m_n)+ nG_{\beta_n,K_n}(m_n)]dx
          }
\nonumber. \eea  \fi

\skp Lemma \ref{lem:AsyOFGmn} gives the asymptotic behavior of
$G_{\beta_n,K_n}(m_n)$. This lemma is used to prove
Lemma \ref{lem:AsyOFZnk} and part (a) of Lemma
\ref{lem:Step1OFThm61b}.

\begin{lem}
\label{lem:AsyOFGmn} We assume that $(\beta_n, K_n)$ satisfies the
hypotheses of Theorem {\em \ref{thm:main}} for all
$0<\alpha<\alpha_0$. Let $m_n=m(\beta_n,K_n)$. Then for all
$0<\alpha<\alpha_0$,
\[\lim_{n \goto \infty}
n^{\alpha/\alpha_0}G_{\beta_n,K_n}(m_n) = g(\bar{x}) < 0 .\]
\end{lem}

{\noindent}{\bf Proof.} We have
\bea
\label{eqn:needthis}
| n^{\alpha/\alpha_0}
G_{\beta_n, K_n}(m_n) - g(\bar{x}) |
&\le&
|n^{\alpha/\alpha_0}
G_{\beta_n, K_n}(n^{\theta\alpha}m_n/n^{\theta\alpha}) -
g(n^{\theta\alpha}m_n)|
\\
&&+ \,
|g(n^{\theta \alpha}m_n)-g(\bar{x})|. \nonumber
\eea
The hypotheses of Theorem \ref{thm:3.1} for all $0 < \alpha < \alpha_0$ consist of a subset of the
hypotheses of Theorem \ref{thm:main}. By hypothesis (iii)(a) of
Theorem \ref{thm:3.1}
\[
\lim_{n \goto \infty} n^{\alpha/\alpha_0} G_{\beta_n,K_n}
(x/n^{\theta\alpha}) = g(x)
\]
uniformly for $x$ in compact subsets of $\R$. According to Theorem \ref{thm:3.1}, $n^{\theta\alpha}m_n \goto
\bar{x}$, and so for any $\varepsilon>0$ the sequence
$n^{\theta\alpha}m_n$ lies in the compact set $[\bar{x}-\varepsilon,
\bar{x}+\varepsilon ]$ for all sufficiently large $n$. Setting $x=n^{\theta
\alpha}m_n$, we see that the first term on the right-hand side of (\ref{eqn:needthis})
has the limit 0. Because of the limit $n^{\theta\alpha}m_n
\goto \bar{x}$ and the continuity of  $g$, the second term on the
right-hand side of (\ref{eqn:needthis}) also converges to 0 as $n \goto
\infty$. It follows that
\[
\lim_{n \goto \infty} | n^{\alpha/\alpha_0} G_{\beta_n, K_n}(m_n) - g(\bar{x}) | = 0.
\]
By hypothesis (iii)(b) of Theorem \ref{thm:3.1}, $\bar{x} > 0$ is the unique
nonnegative, global minimum point of $g$. Thus
\[ g(\bar{x}) < g(0) = \lim_{n \goto \infty} n^{\alpha/\alpha_0}G_{\beta_n,K_n}(0)
= 0.
\]
The proof of lemma is complete.
$\Box$

\skp Lemma \ref{lem:AsyOFZnk}
\iffalse
states the asymptotic
behavior of $Z_{n, \kappa}$ defined in Lemma
\ref{lem:tEgammaLowerValue}. Part (b) of Lemma \ref{lem:AsyOFZnk}
\fi
gives an inequality involving $nG_{\beta_n,K_n}(m_n)$ and the quantity
$Z_{n,\kappa}$ defined in part (a) of Lemma \ref{lem:tEgammaLowerValue}.
This inequality is used in the proof of Lemma
\ref{lem:intRnkexpIneq} and the proof of part (a) of Lemma
\ref{lem:Substep2bOFThm61b}.

\begin{lem}
\label{lem:AsyOFZnk} We assume that $(\beta_n, K_n)$ satisfies the
hypotheses of Theorem {\em \ref{thm:main}} for all
$0<\alpha<\alpha_0$. For any $0 < \alpha < \alpha_0$ define
$\kappa=\frac{1}{2}(1-\alpha/\alpha_0)+\theta\alpha$ and
\[
Z_{n,\kappa}=\int_{\mathbb{R}} \exp[ -n G_{\beta_n,K_n} (x/n^{\kappa})]dx.
\]
Let $m_n = m(\beta_n,K_n)$.
\iffalse
The following conclusions hold.
\begin{itemize}
\item [\em (a)] $Z_{n,\kappa}=\int_{\mathbb{R}} \exp[ -n G_{\beta_n,K_n} (x/n^\kappa)
]dx$ has the asymptotic behavior
\[ \lim_{n \goto \infty} \frac{1}{n^{1-\alpha/\alpha_0}} \log Z_{n,\kappa}= -
g(\bar{x}).
\]

\item[\em (b)] The asymptotic behaviors of $\log Z_{n,\kappa}$ and
$nG_{\beta_n,K_n}$ are related by
\[ \lim_{n \goto \infty} \frac{1}{n^{1-\alpha/\alpha_0}} \log Z_{n,\kappa}= - \lim_{n \goto \infty}
\frac{1}{n^{1-\alpha/\alpha_0}}nG_{\beta_n,K_n}(m_n).
\]
\end{itemize}
\fi
Then for any $\varepsilon >0$ and all sufficiently large
$n$
\[ \exp[ nG_{\beta_n,K_n}(m_n) ] \cdot Z_{n, \kappa} \le \exp[ \varepsilon
n^{1-\alpha/\alpha_0} ].
\]
\end{lem}

{\noindent}{\bf Proof.} For any $0<\alpha< \alpha_0$ define
\[ Z_{n,\theta\alpha}=\int_{\mathbb{R}} \exp[ -n G_{\beta_n,K_n} (x/n^{\theta\alpha})
]dx.
\]
Changing variables shows that $Z_{n,\kappa} = n^{\kappa - \theta\alpha} Z_{n,\theta\alpha}$.
The MDP stated in Theorem \ref{thm:MDP} is proved in Theorem 8.1 in
\cite{CosEllOtt} via an associated Laplace principle. A key step in
this proof is the limit
\[\lim_{n \goto \infty} \frac{1}{n^{1-\alpha/\alpha_0}}\log \int_{\mathbb{R}} \exp[ n^{1-\alpha/\alpha_0} \psi(x)
- n G_{\beta_n, K_n}(x/n^{\theta\alpha}) ]dx = \sup_{x \in
\mathbb{R}}\{\psi(x)-g(x)\},
\]
where $\psi$ is any bounded, continuous function mapping $\mathbb{R}$
to $\mathbb{R}$. This is proved on page 546 of \cite{CosEllOtt} with $v = -(1 - \alpha/\alpha_0)$ and $\gamma = \theta\alpha$.
Setting $\psi=0$ gives the limit
\[ \lim_{n \goto \infty} \frac{1}{n^{1-\alpha/\alpha_0}}
\log Z_{n, \theta\alpha} = \lim_{n \goto \infty}
\frac{1}{n^{1-\alpha/\alpha_0}} \log \int_{\mathbb{R}} \exp[ -n
G_{\beta_n, K_n} (x/n^{\theta\alpha}) ]dx     \\
= -\inf_{y \in \mathbb{R}} g(y).
\]
\iffalse
Since $\kappa = \frac{1}{2}(1-\alpha/\alpha_0) + \theta\alpha$, we
\fi
Since $Z_{n, \kappa} = n^{\kappa-\theta\alpha} Z_{n, \theta\alpha}$ and
$g$ has a unique positive, global minimum point at $\bar{x}$,
\beas
\lim_{n \goto \infty}
\frac{1}{n^{1-\alpha/\alpha_0}}\log Z_{n,\kappa} & = & \lim_{n \goto
\infty} \frac{1}{n^{1-\alpha/\alpha_0}}\log (n^{\kappa-\theta\alpha}
Z_{n,\theta\alpha}) \\
\iffalse
& =  & \lim_{n \goto \infty}
\frac{1}{n^{1-\alpha/\alpha_0}}\log Z_{n,\theta\alpha} \\
\fi
& = & -\inf_{y \in \mathbb{R}} g(y) = - g(\bar{x}).
\eeas
\iffalse************************** \bea \lim_{n \goto \infty}
\frac{1}{n^{1-\alpha/\alpha_0}}\log Z_{n,\kappa} &=& \lim_{n \goto
\infty} \frac{1}{n^{1-\alpha/\alpha_0}}\log
n^{\kappa-\theta\alpha} Z_{n,\theta\alpha} \nonumber \\
&=&  \lim_{n \goto \infty} \frac{1}{n^{1-\alpha/\alpha_0}}[\log
n^{\kappa-\theta\alpha} + \log Z_{n,\theta\alpha}] \nonumber \\
&=&  \lim_{n \goto \infty} \frac{1}{n^{1-\alpha/\alpha_0}} \cdot
\log n^{\kappa-\theta\alpha} + \lim_{n \goto \infty}
\frac{1}{n^{1-\alpha/\alpha_0}} \cdot \log Z_{n,\theta\alpha}
\nonumber. \eea The first term on the right-hand side is 0.
************************************************************\fi
\iffalse
Since $g$ has a unique positive, global minimum point $\bar{x}$, we
have
\bea
\lim_{n \goto \infty} \frac{1}{n^{1-\alpha/\alpha_0}}\log
Z_{n,\kappa}
&=&
\lim_{n \goto \infty}
\frac{1}{n^{1-\alpha/\alpha_0}} \cdot \log Z_{n,\theta\alpha}
\nonumber \\
&=&
-\inf_{y \in \mathbb{R}} g(y)
=
-g(\bar{x})
\nonumber .
\eea
This proves part (a) of Lemma \ref{lem:AsyOFZnk}.

\skp (b) By part (a) of this lemma and
\fi
By Lemma \ref{lem:AsyOFGmn}
$\lim_{n \goto \infty} n^{\alpha/\alpha_0}G_{\beta_n,K_n}(m_n) = g(\bar{x})$.
\iffalse
\[\lim_{n \goto \infty} \frac{1}{n^{1-\alpha/\alpha_0}}\log
Z_{n,\kappa} = -g(\bar{x}) = - \lim_{n \goto \infty}
n^{\alpha/\alpha_0}G_{\beta_n,K_n}(m_n).
\]
This implies that
\fi
Hence the asymptotic behaviors of $\log Z_{n,\kappa}$
and $nG_{\beta_n,K_n}(m_n)$ are related by
\[ \lim_{n \goto \infty} \frac{1}{n^{1-\alpha/\alpha_0}}\log
Z_{n,\kappa} =  - \lim_{n \goto \infty}
\frac{1}{n^{1-\alpha/\alpha_0}} nG_{\beta_n,K_n}(m_n).
\]
Thus for any $\varepsilon >0$ and all sufficiently large $n$
\[ \frac{1}{n^{1-\alpha/\alpha_0}}\log
Z_{n,\kappa} + \frac{1}{n^{1-\alpha/\alpha_0}} nG_{\beta_n,K_n}(m_n)
\le \varepsilon,
\]
\iffalse
It follows that
\[ \log Z_{n,\kappa} + nG_{\beta_n,K_n}(m_n)\le \varepsilon n^{1-\alpha/\alpha_0}
\]
\fi
or equivalently
\iffalse
\[\exp \left[ \log Z_{n,\kappa} + nG_{\beta_n,K_n}(m_n) \right] =
\fi
\[
\exp[ nG_{\beta_n,K_n}(m_n) ] \cdot Z_{n, \kappa}
 \le \exp(\varepsilon n^{1-\alpha/\alpha_0}).
\]
The proof of Lemma \ref{lem:AsyOFZnk} is complete. $\Box$

\iffalse*********************************************Under the
hypotheses of Theorem \ref{thm:6.1}, part (a) of Lemma
\ref{lem:G(xmn)miusG(mn)} calculates
\[\lim_{n \goto \infty} ( nG_{\beta_n,K_n}(x/n^\kappa + m_n)-
nG_{\beta_n,K_n}(m_n))
\]
Part (b) of Lemma
\ref{lem:G(xmn)miusG(mn)} gives the lower bound of this term with
another restriction on $x/n^\kappa$. Part (a) and part (b) of Lemma
\ref{lem:G(xmn)miusG(mn)} provide the reason why Dominated
Convergence Theorem (DCT) can be used in the proof of Lemma
\ref{lem:Thm61aStep1}. Part (b) of Lemma \ref{lem:G(xmn)miusG(mn)}
is also used for proving part (a) of Lemma \ref{lem:Step1OFThm61b},
and the related term is $A_j$.

MORE EXPLANATION NEEDED    *******************************\fi

\iffalse The proof of Lemma \ref{lem:G(xmn)miusG(mn)} is long and
technical. Taylor expansion will be used. The information about
intervals in which $\theta\alpha_0$ lies, the positivity of
$g^{(j)}(\bar{x})$ for corresponding $j$ and condition $g(x) \goto
\infty$ as $|x| \goto \infty$ also play a role in the proof. The
asymptotic behavior of the derivatives $G^{(j)}_{\beta_n,K_n}(m_n)$
given by hypothesis (iii$^\prime$) of Theorem \ref{thm:main} is also
required. And as we mentioned at the beginning of this section, the
formula $m_n \sim \bar{x}/n^{\theta\alpha}$ in Theorem \ref{thm:3.1}
will be very useful. We separate the whole proof for $x >0$ and
$x<0$ in either case $g$ has degree 4 or $g$ has degree 6
respectively.         \fi

\skp We recall that Step 1 in the proof of part (b) of Theorem
\ref{thm:6.1} is to prove the second limit in (\ref{eqn:Thm61bKey}):
\bea \lefteqn{\hspace{-.5in}  \lim_{n \goto \infty} \tilde{E}_{n,\beta_n,K_n} \{
|S_n/n^{1-\kappa}+W_n/n^{1/2-\kappa}- n^{\kappa}m_n| \ \big| \
S_n/n+W_n/n^{1/2}>\delta m_n
\} } \nonumber \\
&&= \frac{1}{\int_{\mathbb{R}}\exp[-\frac{1}{2}
g^{(2)}(\bar{x})x^2]dx } \cdot \int_{\mathbb{R}} |x| \exp[\ts
-\frac{1}{2} g^{(2)}(\bar{x})x^2]dx.\nonumber \eea
By part (b) of
Lemma \ref{lem:tEgammaLowerValue} the limit of the
conditional expectation equals the limit of the product in the last
two lines of (\ref{eqn:tEh|x|}) with $\bar{\delta} = \delta$. For $\delta \in (0,1)$ this
product has the form
\bea \lefteqn {
\frac{1}{\int_{-n^{\kappa}(1-{\delta})m_n}^{\infty} \exp[-n
G_{\beta_n,K_n}(x/n^{\kappa}+m_n)+nG_{\beta_n, K_n}(m_n)]dx} }
\nonumber\\
&& \cdot \int_{-n^{\kappa}(1-{\delta})m_n}^{\infty} |x| \exp[-n
G_{\beta_n,K_n}(x/n^{\kappa}+m_n)+nG_{\beta_n, K_n}(m_n)]dx
\nonumber.
\eea
The calculation of the limit of this product depends in part on Lemma \ref{lem:Thm61aStep1},
\iffalse
Step 1 in the proof of part (b) of Theorem
\ref{thm:6.1} requires Lemma \ref{lem:Thm61aStep1},
\fi
which will be proved via the Dominated Convergence Theorem (DCT). Two key estimates
are given in the next lemma. Part (b) of the next lemma also removes
an error term \iffalse ($A_j$) \fi that arises in the proof of part
(a) of Lemma \ref{lem:Step1OFThm61b}. The proof of Lemma
\ref{lem:G(xmn)miusG(mn)} is postponed until the end of this
section.

\begin{lem}
\label{lem:G(xmn)miusG(mn)} We assume that $(\beta_n, K_n)$
satisfies the hypotheses of Theorem {\em \ref{thm:main}} for all
$0<\alpha<\alpha_0$. For any $0 < \alpha < \alpha_0$ define
\iffalse
$u=1-\alpha/\alpha_0$, $\gamma = \theta\alpha$, and
\fi
$\kappa=\frac{1}{2}(1-\alpha/\alpha_0) + \theta\alpha$, and let
$m_n=m(\beta_n,K_n)$. The following conclusions hold.

{\em (a)} For all $x \in \mathbb{R}$
\[\lim_{n \goto \infty}
\left(nG_{\beta_n,K_n}( x/n^\kappa + m_n) -nG_{\beta_n,K_n}(m_n)\right) =
\frac{1}{2} g^{(2)}(\bar{x}) x^2. \]

{\em (b)} There exists $\Delta \in (0,1)$ such that for any $\bar{\delta} \in (\Delta, 1)$
there exists $R>0$ such that for all sufficiently large n and all $x
\in \mathbb{R}$ satisfying $|x/n^\kappa| < R$ and $x/n^\kappa >
-(1-\bar{\delta})m_n$
\[ nG_{\beta_n,K_n}( x/n^\kappa + m_n) -nG_{\beta_n,K_n}(m_n) \ge
\frac{1}{8} g^{(2)}(\bar{x})x^2.
\]
\end{lem}

\iffalse \skp Lemma \ref{lem:intRnkexpIneq} is used in the proof of
Lemma \ref{lem:Thm61aStep1}(a). It is also needed in the proof of
Lemma \ref{lem:Step1OFThm61b}(a), and the related term is $B_n$.
\iffalse The proof of Lemma \ref{lem:intRnkexpIneq} requires Lemma
\ref{lem:RepresentF}, the inequality in Lemma \ref{lem:AsyOFZnk}(b)
and the Large Deviation Principle (LDP) of $S_n/n + W_n/n^{1/2}$ in
Lemma 4.4(b) in \cite{CosEllOtt}. \fi  \fi

The next lemma removes an error term that arises in applying the DCT to
prove Lemma \ref{lem:Thm61aStep1}. The next lemma also removes an
error term \iffalse ($B_n$) \fi that arises in the proof of part (a)
of Lemma \ref{lem:Step1OFThm61b}.

\begin{lem}
\label{lem:intRnkexpIneq} We assume that $(\beta_n, K_n)$ satisfies
the hypotheses of Theorem {\em \ref{thm:main}} for all
$0<\alpha<\alpha_0$. Then there exist a constant $c_2>0$ such that
for all sufficiently large $n$
\[ \int_{Rn^\kappa}^\infty  \exp[-nG_{\beta_n,K_n}( x/n^\kappa + m_n) +
nG_{\beta_n,K_n}(m_n)]dx \le \exp[-c_2 n] \goto 0  \ \ \mbox{as} \ \
n \goto \infty,
\]
where $R$ is chosen as in part {\em (b)} of Lemma {\em
\ref{lem:G(xmn)miusG(mn)}} and $m_n = m(\beta_n,K_n)$.
\end{lem}

{\noindent}{\bf Proof.}
We start by applying Lemma \ref{lem:RepresentF} with
$\varphi(x)= 1_{(Rn^\kappa + n^\kappa m_n, \infty)}(x)$ and $\bar{\gamma} = \kappa$, obtaining
\bea
\lefteqn{ (P_{n,\beta_n,K_n} \times Q) \{ S_n/n+W_n/n^{1/2} \ge
R+m_n \}      } \nonumber  \\
&=&
(P_{n,\beta_n,K_n} \times Q) \{ S_n/n^{1-\kappa}
+W_n/n^{1/2-\kappa} \ge Rn^\kappa +m_n n^\kappa \} \nonumber \\
&=&
\int_{\Lambda^n \times \Omega} 1_{  \{ S_n/n^{1-\kappa}
+W_n/n^{1/2-\kappa} \ge Rn^\kappa +m_n n^\kappa  \}  }  d
(P_{n,\beta_n,K_n} \times Q) \nonumber \\
&=&
\frac{1}{Z_{n,\kappa}} \cdot \int_{\mathbb{R}} 1_{[Rn^\kappa+m_n
n^\kappa ,\infty)}(x) \exp[-n G_{\beta_n,K_n}(x/n^{\kappa})]dx
\nonumber \\
& = & \frac{1}{Z_{n,\kappa}} \cdot \int_{Rn^\kappa}^{\infty} \exp[-n G_{\beta_n,K_n}(x/n^{\kappa}
+m_n)]dx, \nonumber
\eea
where $Z_{n,\kappa} = \int_{\R} \exp[-nG_{\beta_n,K_n}(x/n^\kappa)] dx$.
\iffalse
The change of variables $y=x-n^{\kappa}m_n$ in the last line of the display gives
\[
\frac{1}{Z_{n,\kappa}} \cdot \int_{Rn^\kappa}^{\infty} \exp[-n G_{\beta_n,K_n}(y/n^{\kappa}
+m_n)]dy,
\]
\fi
Thus we have
\bea
\lefteqn{   \int_{Rn^\kappa}^\infty
\exp[-nG_{\beta_n,K_n}( x/n^\kappa + m_n) + nG_{\beta_n,K_n}(m_n)]dx     } \nonumber \\
&&=
\exp[nG_{\beta_n,K_n}(m_n)] \cdot Z_{n,\kappa} \cdot
(P_{n,\beta_n,K_n} \times Q) \{ S_n/n+W_n/n^{1/2} \ge R+m_n \}
\nonumber \\
&&\le
\exp[nG_{\beta_n,K_n}(m_n)] \cdot Z_{n,\kappa} \cdot
(P_{n,\beta_n,K_n} \times Q) \{ S_n/n+W_n/n^{1/2} \ge R \}
\nonumber.
\eea
By part (b) of Lemma 4.4 in \cite{CosEllOtt}, with respect to
$P_{n,\beta_n,K_n} \times Q$, $S_n/n+W_n/n^{1/2}$ satisfies the large deviation principle
on $\mathbb{R}$ with exponential speed $n$ and rate function
$G_{\beta,K(\beta)}$. In particular, for the closed set $[  R,
\infty  )$ we have the large deviation upper bound
\[ \limsup_{n \goto \infty} \frac{1}{n} \log (P_{n,\beta_n,K_n} \times Q) \{ S_n/n+W_n/n^{1/2}
\ge R \} \le - \inf_{x \ge R} G_{\beta,K(\beta)}(x).
\]
By part (a) of Theorem \ref{thm:secondorder}, since $0<\beta \le
\beta_c$, we have ${\mathcal{M}}_{\beta,K(\beta)} = \{0\}$. Thus
$G_{\beta,K(\beta)}$ has a unique global minimum point at 0. Since
$R>0$, it follows that
\[ \inf_{x \ge R} G_{\beta,K(\beta)}(x) >
\inf_{x \in \mathbb{R}} G_{\beta,K(\beta)}(x) =0.
\]
Therefore for all sufficiently large $n$
\[(P_{n,\beta_n,K_n} \times Q) \{ S_n/n+W_n/n^{1/2} \ge R
\} \le \exp[-c_1 n],
\]
where $c_1 = \inf_{x \ge R} G_{\beta,K(\beta)}(x)/2 > 0$.
\iffalse
Therefore there exists $c_1>0$ such that for all sufficiently large $n$
\[(P_{n,\beta_n,K_n} \times Q) \{ S_n/n+W_n/n^{1/2} \ge R
\} \le \exp[-2 c_1 n].
\]
\fi
We now appeal to Lemma \ref{lem:AsyOFZnk}, which states that for any $\varepsilon
>0$ and all sufficiently large $n$
\[ \exp[nG_{\beta_n,K_n}(m_n)] \cdot Z_{n,\kappa} \le \exp[\varepsilon
n^{1-\alpha/\alpha_0}].
\]
Since $0 < 1 - \alpha/\alpha_0 < 1$, it follows that for all sufficiently large $n$
\bea
\lefteqn{
\int_{Rn^\kappa}^\infty
\exp[-nG_{\beta_n,K_n}( x/n^\kappa + m_n) + nG_{\beta_n,K_n}(m_n)]dx     } \nonumber \\
&&\le
\exp[nG_{\beta_n,K_n}(m_n)] \cdot Z_{n,\kappa} \cdot
(P_{n,\beta_n,K_n} \times Q) \{ S_n/n+W_n/n^{1/2} \ge R \} \nonumber \\
&&\le
\exp[\varepsilon n^{1-\alpha/\alpha_0}] \cdot \exp[-c_1 n] \nonumber \\
&&\le
\exp[- c_1 n/2]  \nonumber.
\eea
\iffalse
The last inequality follows
from the fact $0 < 1-\alpha/\alpha_0 <1$.
\fi
This gives the conclusion of Lemma \ref{lem:intRnkexpIneq} with $c_2 = c_1/2$. The proof of the lemma
is complete. $\Box$

\iffalse \skp  Part (a) of Lemma \ref{lem:Thm61aStep1}, i.e.,
formula (\ref{eqn:intLowerfexp}) is used in the proof of substep 2b
of both parts of Theorem \ref{thm:6.1} (in Lemma
\ref{lem:Substep2bOFThm61a} and Lemma \ref{lem:Substep2bOFThm61b})
and the proof of both parts of Lemma \ref{lem:Step1OFThm61b}. The
part (b) of Lemma \ref{lem:Thm61aStep1} is just the step 1 of
Theorem \ref{thm:6.1}(a). \fi

\skp Lemma \ref{lem:Thm61aStep1} is a key result in the proof of the conditional
limit stated in part (b) of Theorem \ref{thm:6.1}. The lemma deals with the weak convergence of
certain measures needed in the proof of part (a) of Lemma
\ref{lem:Step1OFThm61b}. Lemma \ref{lem:Thm61aStep1} is also used
with $f \equiv 1$ in the proof of part (b) of Lemma
\ref{lem:Step1OFThm61b} and part (a) of Lemma
\ref{lem:Substep2bOFThm61b}.

\begin{lem}
\label{lem:Thm61aStep1} We assume that $(\beta_n, K_n)$ satisfies
the hypotheses of Theorem {\em \ref{thm:main}} for all
$0<\alpha<\alpha_0$. Given any $\alpha \in (0, \alpha_0)$, define
$\kappa = \frac{1}{2}(1-\alpha/\alpha_0)+\theta\alpha$ and
\[
Z_{n,\kappa}= \int_{\mathbb{R}} \exp[-n
G_{\beta_n,K_n}(x/n^{\kappa})]dx.
\]
 For
$\bar{\delta} \in (0,1)$ define
\[ A_n(\bar{\delta})=\{ S_n/n + W_n/n^{1/2} > \bar{\delta}m_n \},
\]
where $m_n=m(\beta_n,K_n)$. Let $f$ be any bounded,
continuous function and let $\Delta \in (0,1)$ be the number
determined in part {\em (b)} of Lemma {\em \ref{lem:G(xmn)miusG(mn)}}. Then for any
$0<\alpha<\alpha_0$ and any $\bar{\delta} \in (\Delta, 1)$ we have
the limit
\bea \lefteqn{
\lim_{n \goto \infty} \exp[
nG_{\beta_n,K_n}(m_n) ] \cdot Z_{n, \kappa} \cdot
\tilde{E}_{n,\beta_n,K_n} \{ f(S_n/n^{1-\kappa}+ W_n/n^{1/2-\kappa}
-n^{\kappa} m_n ) \cdot 1_{A_n(\bar{\delta})} \}
} \nonumber \\
&&= \lim_{n \goto \infty}
\int_{-n^\kappa(1-\bar{\delta})m_n}^{\infty}
f(x)\exp[-nG_{\beta_n,K_n}(x/n^\kappa+m_n)+ nG_{\beta_n,K_n}(m_n)]dx
\nonumber
\\
&&= \int_{\mathbb{R}} f(x)\exp[\ts -\frac{1}{2}g^{(2)}(\bar{x})
x^2]dx. \label{eqn:anotherone}
\eea
\end{lem}

{\noindent}{\bf Proof.} The first equality follows by applying
part (a) of Lemma \ref{lem:tEgammaLowerValue} to $h=f$. Concerning the second
equality, we denote by $I_n$ the integral in the second line of
(\ref{eqn:anotherone}). We write
$I_n=I_{n_1} + I_{n_2}$, where
\[ I_{n_1} = \int_{-n^\kappa(1-\bar{\delta})m_n}^{Rn^\kappa}
f(x)\exp[-nG_{\beta_n,K_n}(x/n^\kappa+m_n)+ nG_{\beta_n,K_n}(m_n)]dx
\]
and
\[ I_{n_2}= \int_{Rn^\kappa}^{\infty}
f(x)\exp[-nG_{\beta_n,K_n}(x/n^\kappa+m_n)+
nG_{\beta_n,K_n}(m_n)]dx.
\]
The number $R$ is chosen as in part (b) of Lemma
\ref{lem:G(xmn)miusG(mn)}. Since $f$ is bounded,
Lemma \ref{lem:intRnkexpIneq} implies that there exists $c_2 > 0$ such that for all sufficiently large $n$
\bea
I_{n_2}
&\le&
\|f\|_\infty \int_{Rn^\kappa}^\infty \exp[-nG_{\beta_n,K_n}(
x/n^\kappa + m_n) +
nG_{\beta_n,K_n}(m_n)]dx \nonumber \\
&\le&
\|f\|_\infty \exp[-c_2 n] \goto 0  \ \ \mbox{as} \ \ n \goto \infty
\nonumber. \eea Thus $I_{n_2} \goto 0$ as $n
\goto \infty.$

Define
\[ h_n(x) = f(x) \exp[-nG_{\beta_n,K_n}(x/n^\kappa+m_n) +
nG_{\beta_n,K_n}(m_n)]
\]
and
\[ h(x) = f(x) \exp[\ts -\frac{1}{2}g^{(2)}(\bar{x})x^2].
\]
By part (a) of Lemma \ref{lem:G(xmn)miusG(mn)}, $h_n(x) \goto h(x)$
for all $x \in \mathbb{R}$. In addition, by part (b) of Lemma
\ref{lem:G(xmn)miusG(mn)}, if $x \in ( -(1-\bar{\delta})m_n
n^\kappa, Rn^\kappa )$, then for all sufficiently large $n$
\[nG_{\beta_n,K_n}( x/n^\kappa + m_n) -nG_{\beta_n,K_n}(m_n) \ge
H(x)= \frac{1}{8} g^{(2)}(\bar{x})x^2.
\]
\iffalse
It follows that $|h_n(x)| \leq \|f\|_\infty \exp[\ts -\frac{1}{8} g^{(2)}(\bar{x})x^2]$.
\bea
|h_n(x)|
&=&
|f(x)|
\exp[-nG_{\beta_n,K_n}(x/n^\kappa+m_n) +
nG_{\beta_n,K_n}(x/n^\kappa)] \nonumber \\
&\le&
\|f\|_\infty
\exp[\ts -\frac{1}{8} g^{(2)}(\bar{x})x^2] \nonumber. \eea
\fi
Since $\exp[-H(x)]$ is integrable, the Dominated Convergence Theorem implies that
\[
\lim_{n \goto \infty} I_{n_1} = \lim_{n \goto \infty}\int_{-n^\kappa(1-\bar{\delta})m_n}^{Rn^\kappa} h_n(x)dx
= \int_{\mathbb{R}} h(x) dx = \int_{\mathbb{R}} f(x)
\exp[\ts -\frac{1}{2}g^{(2)}(\bar{x})x^2]dx.
\]
\iffalse
i.e.,
\[ \lim_{n \goto \infty} I_{n_1} = \int_{\mathbb{R}} f(x)
\exp[\ts -\frac{1}{2}g^{(2)}(\bar{x})x^2]dx.
\]
\fi
We conclude that
\[ \lim_{n \goto \infty} I_n =
\lim_{n \goto \infty} I_{n_1}
+ \lim_{n \goto \infty} I_{n_2}
= \int_{\mathbb{R}} f(x) \exp[\ts
-\frac{1}{2}g^{(2)}(\bar{x})x^2]dx.
\]
This completes the proof of Lemma \ref{lem:Thm61aStep1}. $\Box$

\iffalse \skp LEMMAS 7.1-7.7 HERE ARE THE SAME CONTENT AND ORDER AS
LEMMA 7.1-7.7 IN THESIS. OF COURSE I MAKE SOME SMALL CHANGES INSIDE
EACH LEMMA.

\skp Part (b) of the next lemma is used in the proof of part (a) of
Lemma \ref{lem:Substep2bOFThm61b}. I DELETE LEMMA 7.8 IN THESIS,
WHICH WAS SUBSTEP 2A IN THM 6.1(A). SO THE LEMMA 7.9 IN THESIS
BECOMES LEMMA 7.8 NOW. \fi

\skp
The next lemma collects several elementary but useful facts concerning the normal random variables $W_n$.

\begin{lem}
\label{lem:normalrvs}
Let $(\bn,\kn)$ be a positive sequence that
converges either to a second-order point $(\beta,K(\beta))$, $0 <
\beta < \bc$, or to the tricritical point $(\beta,K(\beta)) =
(\bc,K(\bc))$. Let $W_n$ be a sequence of normal random
variables with mean 0 and variance $\sigma_n^2 = (2\beta_n K_n)^{-1}$ defined on
a probability space $(\Omega,\mathcal {F},Q)$. The following conclusions hold.

{\em (a)} For $b > 0$ and $\zeta > 0$ there exists a constant $c > 0$ such that for all $n$,
$Q\{|W_n| > b n^\zeta\} \leq \exp[-c n^{2\zeta}]$.

{\em (b)} There exist a constant $c_1 > 0$ such that for all $n$
\[
\int_\Omega |W_n|^2 dQ \leq c_1 \ \mbox{ and } \ \int_\Omega |W_n| dQ \leq \sqrt{c_1}.
\]
\end{lem}

\noi
{\bf Proof.} (a) We have the bound
\[
Q\{|W_n| > b n^\zeta\} = {\frac{\sqrt{2}}{\sqrt{\pi}\sigma_n}} \int_{b n^\zeta}^\infty \exp[-x^2/2\sigma_n^2] dx
\leq \frac{\sqrt{2}\sigma_n}{\sqrt{\pi} b n^\zeta} \exp[-b^2 n^{2\zeta}/2\sigma_n^2].
\]
Part (a) now follows from the fact that
since $(\beta_n,K_n)$ is a positive sequence converging to $(\beta,K(\beta))$ for $0 < \beta \leq \beta_c$,
the positive sequences $\sigma_n$ and $\sigma_n^2$ are bounded.

(b) Since $\int_\Omega |W_n|^2 dQ = \sigma_n^2$ and $\int_\Omega |W_n| dQ \leq (\int_\Omega |W_n|^2 dQ)^{1/2} = \sigma_n$, this follows from
the fact that the positive sequences $\sigma_n^2$ and $\sigma_n$
are bounded. The proof of the lemma is complete. $\Box$

\skp The next lemma is used in the proof of part (a) of Lemma
\ref{lem:Substep2bOFThm61b}. Under the hypotheses of Theorem \ref{thm:main}, for any $0 < \alpha < \alpha_0$ the interval
$(0,\frac{1}{2}-\theta\alpha)$ appearing in the next lemma is nonempty because
by hypothesis (iii$'$) $\half - \theta\alpha > \half - \theta\alpha_0 > 0$.

\begin{lem}
\label{lem:IneqOFtwoEvent} We assume that $(\beta_n, K_n)$ satisfies
the hypotheses of Theorem {\em \ref{thm:main}} for all
$0<\alpha<\alpha_0$. For $\bar{\delta} \in (0,1)$ define
\[ A_n(\bar{\delta})=\{ S_n/n + W_n/n^{1/2} > \bar{\delta}m_n \},
\]
where $m_n=m(\beta_n,K_n)$.
\iffalse
Let $f$ be any nonnegative, bounded, uniformly
continuous function and
\fi
Let $\Delta \in (0,1)$ be the number
determined in part {\em (b)} of Lemma {\em \ref{lem:G(xmn)miusG(mn)}}. Assume that
$0<\alpha < \alpha_0$ and choose any numbers $\delta_1, \delta,
\delta_2$ and $\zeta$ satisfying $\Delta <
\delta_1<\delta<\delta_2<1$ and $\zeta \in (0,
\frac{1}{2}-\theta\alpha)$. Then there exist constants $c>0$ and $c_2>0$ such that the following conclusions hold.
\iffalse
\item [\em (a)] For all sufficiently large $n$
\bea \lefteqn{
\tilde{E}_{n,\beta_n,K_n} \{
f(S_n/n^{1-\kappa}+W_n/n^{1/2-\kappa}-n^\kappa m_n) \cdot
1_{A_n(\delta_1)} \} +\|f\|_\infty e^{-cn^{2\zeta}}    } \nonumber \\
&&\ge
\tilde{E}_{n,\beta_n,K_n} \{
f(S_n/n^{1-\kappa}+W_n/n^{1/2-\kappa}-n^\kappa m_n) \cdot
1_{ \{ S_n/n>\delta m_n \} } \}  \nonumber \\
&&\ge \tilde{E}_{n,\beta_n,K_n} \{
f(S_n/n^{1-\kappa}+W_n/n^{1/2-\kappa}-n^\kappa m_n) \cdot
1_{A_n(\delta_2)} \} - \|f\|_\infty e^{-cn^{2\zeta}} \nonumber. \eea
\fi

{\em (a)} For all sufficiently large $n$
\bea
(P_{n,\beta_n,K_n}\times Q) \{ A_n(\delta_1)  \} + e^{-cn^{2\zeta}}  & \geq &
P_{n,\beta_n,K_n}  \{ S_n/n > \delta m_n \}  \nonumber \\
& \ge &
(P_{n,\beta_n,K_n}\times Q) \{ A_n(\delta_2)  \} -
e^{-cn^{2\zeta}} \nonumber. \eea

{\em (b)} For all sufficiently large $n$
\bea \lefteqn{ \tilde{E}_{n,\beta_n,K_n} \{
|S_n/n^{1-\kappa}+W_n/n^{1/2-\kappa}-n^\kappa m_n| \cdot
1_{A_n(\delta_1)} \}     }\nonumber \\
&&+ \, 2n^{\kappa} e^{-cn^{2\zeta}} + c_2n^{\kappa-1/2}e^{-c n^{2\zeta}/2}     \nonumber \\
&\ge& \tilde{E}_{n,\beta_n,K_n} \{
|S_n/n^{1-\kappa}+W_n/n^{1/2-\kappa}-n^\kappa m_n| \cdot
1_{ \{ S_n/n>\delta m_n \} } \}  \nonumber \\
&\ge& \tilde{E}_{n,\beta_n,K_n} \{
|S_n/n^{1-\kappa}+W_n/n^{1/2-\kappa}-n^\kappa m_n| \cdot
1_{A_n(\delta_2)} \} \nonumber \\
&&- \, 2n^{\kappa} e^{-cn^{2\zeta}} -
c_2n^{\kappa-1/2}e^{-c n^{2\zeta}/2} \nonumber. \eea
\end{lem}

{\noindent}{\bf Proof of part (a) of Lemma \ref{lem:IneqOFtwoEvent}.} We choose $\zeta \in (0,
\frac{1}{2}-\theta\alpha)$. The proof is based on the following two claims, which are proved later.
\iffalse
The proof of Claim 1 appears at end of the proof of this lemma; the analogous
proof of Claim 2 is omitted.
\fi

\skp {\bf Claim 1.}  \ \ \ For all sufficiently large $n$, $ \{S_n/n
> \delta m_n \} \subset A_n(\delta_1) \cup \{ |W_n|>\half n^{\zeta} \}. $

\skp {\bf Claim 2.}  \ \ \ For all sufficiently large $n$, $ \{S_n/n
> \delta m_n \}  \supset A_n(\delta_2)  \backslash  \{
|W_n|> \half n^{\zeta} \}.$

\skp
\iffalse
Since
\[ P_{n,\beta,K_n} \{ S_n/n > \delta m_n \} = E_{n,\beta_n,K_n} \{ 1_{\{ S_n/n > \delta m_n \}}
\} = \tilde{E}_{n,\beta_n,K_n} \{ 1_{\{ S_n/n > \delta m_n \}} \},
\]
\[ P_{n,\beta,K_n} \times Q \{ A_n(\delta_1) \} = \tilde{E}_{n,\beta_n,K_n} \{ 1_{A_n(\delta_1)}
\},
\]
and
\[ P_{n,\beta,K_n} \times Q \{ A_n(\delta_2) \} = \tilde{E}_{n,\beta_n,K_n} \{ 1_{A_n(\delta_2)}
\},
\]
we can prove
\fi
\iffalse
We prove part (a) by showing that
\bea
\tilde{E}_{n,\beta_n,K_n} \{
1_{A_n(\delta_1)} \} + e^{-cn^{2\zeta}}  & \ge &
E_{n,\beta_n,K_n} \{
1_{ \{ S_n/n>\delta m_n \} } \}  \nonumber \\
& \ge & \tilde{E}_{n,\beta_n,K_n} \{  1_{A_n(\delta_2)} \} -
e^{-cn^{2\zeta}} \nonumber. \eea
We use Claim 1 to prove the first inequality in this display.
For all sufficiently large $n$
\bea \tilde{E}_{n,\beta_n,K_n} \{1_{ \{
S_n/n>\delta m_n \} } \} &=& \int_{\mathbb{R}} 1_{ \{
S_n/n>\delta m_n \} }  d (P_{n, \beta_n, K_n} \times Q)  \\
&\le& \int_{\mathbb{R}} 1_{
A_n(\delta_1) \cup \{ |W_n|>n^{\zeta} \} }  d (P_{n, \beta_n, K_n} \times Q)  \nonumber \\
&\le& \int_{\mathbb{R}} 1_{
A_n(\delta_1) }  d (P_{n, \beta_n, K_n} \times Q ) \nonumber \\
&&+ \int_{\mathbb{R}} 1_{
 \{ |W_n|>n^{\zeta} \} }  d (P_{n, \beta_n, K_n} \times Q)
\nonumber. \eea
\fi
By Claims 1 and 2, for all sufficiently large $n$
\beas
\iffalse \label{eqn:TwoEventf1} \fi
\lefteqn{
(P_{n, \beta_n, K_n} \times Q)\{A_n(\delta_1)\} + Q\{\{ |W_n|>\half n^{\zeta} \} } \\
& & = (P_{n, \beta_n, K_n} \times Q)\{A_n(\delta_1)\} + (P_{n, \beta_n, K_n} \times Q)\{\{ |W_n|> \half n^{\zeta} \} \nonumber \\
& & \geq  (P_{n, \beta_n, K_n} \times Q)\{S_n/n>\delta m_n \} = P\{S_n/n>\delta m_n \} \nonumber \\
\nonumber && \geq (P_{n, \beta_n, K_n} \times Q)\{A_n(\delta_2)\} - Q\{ |W_n|>\half n^{\zeta} \}.
\eeas
Part (a) of Lemma \ref{lem:normalrvs} completes the proof. Thus, given Claims
1 and 2, the proof of part (a) is complete.

\skp{\noindent}{\bf Proof of part (b) of Lemma \ref{lem:IneqOFtwoEvent}.} We use Claim 1 to
prove the first inequality in part (b). For all
sufficiently large $n$
\bea \label{eqn:TwoEventAbsolute1} \lefteqn{
\tilde{E}_{n,\beta_n,K_n} \{
|S_n/n^{1-\kappa}+W_n/n^{1/2-\kappa}-n^\kappa m_n| \cdot 1_{ \{
S_n/n>\delta m_n \} }
\}      }  \\
&=& \int_{\Lambda^n \times \Omega} |S_n/n^{1-\kappa}+W_n/n^{1/2-\kappa}-n^\kappa
m_n| \cdot 1_{ \{
S_n/n>\delta m_n \} }  d (P_{n, \beta_n, K_n} \times Q ) \nonumber \\
\iffalse
&\le& \int_{\Lambda^n \times \Omega}
|S_n/n^{1-\kappa}+W_n/n^{1/2-\kappa}-n^\kappa m_n| \cdot 1_{
A_n(\delta_1) \cup \{ |W_n|>n^\zeta \} }  d (P_{n, \beta_n, K_n} \times Q)  \nonumber \\
\fi
&\le& \int_{\Lambda^n \times \Omega}
|S_n/n^{1-\kappa}+W_n/n^{1/2-\kappa}-n^\kappa m_n| \cdot 1_{
A_n(\delta_1) }  d (P_{n, \beta_n, K_n} \times Q ) \nonumber \\
&&+ \int_{\Lambda^n \times \Omega} |S_n/n^{1-\kappa}+W_n/n^{1/2-\kappa}-n^\kappa
m_n| \cdot 1_{
 \{ |W_n|> \half n^\zeta \} }  d (P_{n, \beta_n, K_n} \times Q)
\nonumber \\
&\le& \int_{\Lambda^n \times \Omega}
|S_n/n^{1-\kappa}+W_n/n^{1/2-\kappa}-n^\kappa m_n| \cdot 1_{
A_n(\delta_1) }  d (P_{n, \beta_n, K_n} \times Q)  \nonumber \\
&&+ \int_{\Lambda^n \times \Omega} |S_n/n^{1-\kappa}-n^\kappa m_n| \cdot 1_{
 \{ |W_n|> \half n^\zeta \} }  d (P_{n, \beta_n, K_n} \times Q ) \nonumber \\
&&+ \int_{\Omega} |W_n/n^{1/2-\kappa}| \cdot 1_{
 \{ |W_n|> \half n^\zeta \} }  d Q
\nonumber. \eea Since $|S_n/n| \le 1$ and $m_n \in (0,1)$, we have
$|S_n/n^{1-\kappa}-n^\kappa m_n| \le 2n^\kappa$.
\iffalse
Using the fact that
\[ Q\{ |W_n|>n^\zeta \} \le \frac{\sqrt{2} \sigma_n}{\sqrt{\pi}n^\zeta} \cdot \exp( -n^{2\zeta} /2 \sigma_n^2
),
\]
\fi
Using part (a) of Lemma \ref{lem:normalrvs}, for all sufficiently large $n$ we bound the next to last integral
in (\ref{eqn:TwoEventAbsolute1}) by
\[
\iffalse \label{eqn:AbsoluteQ} \fi
2n^\kappa \cdot Q\{ |W_n|> \half n^\zeta \} \le
2n^\kappa \exp(-c n^{2 \zeta}),
\]
where $c > 0$ is a constant. The next step is to apply the Cauchy-Schwartz inequality to the last
integral in (\ref{eqn:TwoEventAbsolute1}) and use parts (a)
and (b) of Lemma \ref{lem:normalrvs}. There exist constants $c > 0$ and $c_2 = \sqrt{c_1} > 0$ such that for all $n$
\beas
\lefteqn{ \int_{\Omega} |W_n/n^{1/2-\kappa}| \cdot 1_{
 \{ |W_n|> \half n^\zeta \} }  d Q    }
\nonumber \\
&& \le \left( \int_{\Omega} |W_n/n^{1/2-\kappa}|^2 d Q
\right)^{1/2}  \cdot  \left( Q\{ |W_n|> \half n^\zeta \}\right)^{1/2}
\leq c_2 n^{\kappa - 1/2} \exp[-cn^{2\zeta}/2] \nonumber.
\eeas
\iffalse
Since $W_n$ is a sequence of $N(0, \sigma_n^2)$
random variables defined on a probability space $(\Omega, \mathcal
{F}, Q)$, where $\sigma_n^2=(2\beta_n K_n)^{-1}$, we have \iffalse
\[ \tilde{E}_{n, \beta_n, K_n}\{W_n\}=0.
\]
Then \fi
\[\tilde{E}_{n, \beta_n, K_n}\{(W_n)^2\}= \iffalse  \sigma_n^2 + (\tilde{E}_{n, \beta_n,
K_n}\{W_n\})^2 = \sigma_n^2= \fi (2\beta_n K_n)^{-1}.
\]
Since $(\beta_n, K_n) \goto (\beta, K(\beta))$ with
$0<\beta<\beta_c$, there exists a constant $c_2>0$ such that for all
$n \in \N$
\bea
\int_{\mathbb{R}} |W_n/n^{1/2-\kappa}|^2 d Q
&=&
\tilde{E}_{n, \beta_n, K_n}\{|W_n/n^{1/2-\kappa}|^2\} \nonumber \\
&=& \frac{1}{(n^{1/2-\kappa})^2} \cdot \tilde{E}_{n, \beta_n,
K_n}\{(W_n)^2\} \le \frac{(c_2)^2}{(n^{1/2-\kappa})^2} \nonumber.
\eea
Because
\[ \int_{\mathbb{R}} 1_{ \{ |W_n|>n^\zeta \} } d Q =
Q\{|W_n|>n^\zeta\} \le \exp(-c n^{2\zeta})\nonumber,
\]
the last integral in
(\ref{eqn:TwoEventAbsolute1}) is bounded by
\be \label{eqn:WnBound}
\frac{c_2}{n^{1/2-\kappa}} \cdot \exp(-cn^{2\zeta}/2) \le c_2
n^{\kappa-1/2} \cdot \exp(-c_3n^{2\zeta}),
\ee
where $c_3>0$ is a constant.
\fi
It follows that for all sufficiently large $n$
\bea
\label{eqn:again}
\lefteqn{ \tilde{E}_{n,\beta_n,K_n} \{
|S_n/n^{1-\kappa}+W_n/n^{1/2-\kappa}-n^\kappa m_n| \cdot 1_{ \{
S_n/n>\delta m_n \} } \}   }  \\
&&\le \tilde{E}_{n,\beta_n,K_n} \{
|S_n/n^{1-\kappa}+W_n/n^{1/2-\kappa}-n^\kappa m_n| \cdot
1_{A_n(\delta_1)} \} \nonumber \\
&& \hspace{.15in} + \, 2n^\kappa e^{-cn^{2\zeta}} +
c_2n^{\kappa-1/2}e^{-c n^{2\zeta}/2} \nonumber.
\eea
This completes
the proof of the first inequality in part (b).

We now use Claim 2 to prove the second inequality in part (b).
For all sufficiently large $n$ \bea
\label{eqn:TwoEventAbsolute2} \lefteqn{ \tilde{E}_{n,\beta_n,K_n} \{
|S_n/n^{1-\kappa}+W_n/n^{1/2-\kappa}-n^\kappa m_n| \cdot 1_{ \{
S_n/n>\delta m_n \} }
\}      }  \\
&=& \int_{\Lambda^n \times \Omega} |S_n/n^{1-\kappa}+W_n/n^{1/2-\kappa}-n^\kappa
m_n| \cdot 1_{ \{
S_n/n>\delta m_n \} }  d P_{n, \beta_n, K_n}  \nonumber \\
\iffalse
&\ge& \int_{\Lambda^n \times \Omega}
|S_n/n^{1-\kappa}+W_n/n^{1/2-\kappa}-n^\kappa m_n| \cdot 1_{
A_n(\delta_2) \backslash \{ |W_n| > \half n^\zeta \} }  d (P_{n, \beta_n, K_n} \times Q)  \nonumber \\
\fi
&\ge& \int_{\Lambda^n \times \Omega}
|S_n/n^{1-\kappa}+W_n/n^{1/2-\kappa}-n^\kappa m_n| \cdot 1_{
A_n(\delta_2) }  d (P_{n, \beta_n, K_n} \times Q)  \nonumber \\
&&- \int_{\Lambda^n \times \Omega} |S_n/n^{1-\kappa}+W_n/n^{1/2-\kappa}-n^\kappa
m_n| \cdot 1_{
 \{ |W_n| > \half n^\zeta \} }  d (P_{n, \beta_n, K_n} \times Q)
\nonumber \\
&\ge& \int_{\Lambda^n \times \Omega}
|S_n/n^{1-\kappa}+W_n/n^{1/2-\kappa}-n^\kappa m_n| \cdot 1_{
A_n(\delta_2) }  d (P_{n, \beta_n, K_n} \times Q ) \nonumber \\
&&- \int_{\Lambda^n \times \Omega} |S_n/n^{1-\kappa}-n^\kappa m_n| \cdot 1_{
 \{ |W_n| > \half n^\zeta\} }  d (P_{n, \beta_n, K_n} \times Q ) \nonumber \\
&&- \int_{\Omega} |W_n/n^{1/2-\kappa}| \cdot 1_{
 \{ |W_n| > \half n^\zeta \} }  d Q
\nonumber. \eea
\iffalse
The next to last inequality follows from the fact
that for any sets $C$ and $D$, $ 1_{C \backslash D} \ge 1_C-1_D$.
The second term after the last inequality in
(\ref{eqn:TwoEventAbsolute2}) coincides with the second term after
the last inequality  in (\ref{eqn:TwoEventAbsolute1}) and so can be
bounded using (\ref{eqn:AbsoluteQ}). The last term after the last
inequality in (\ref{eqn:TwoEventAbsolute2}) coincides with the last
term after the last inequality in (\ref{eqn:TwoEventAbsolute1}) and
so can be bounded using (\ref{eqn:WnBound}). Hence from
(\ref{eqn:TwoEventAbsolute2}) we obtain for all sufficiently large
$n$ \bea \lefteqn{ \tilde{E}_{n,\beta_n,K_n} \{
|S_n/n^{1-\kappa}+W_n/n^{1/2-\kappa}-n^\kappa m_n| \cdot 1_{ \{
S_n/n>\delta m_n \} }
\}      } \nonumber \\
&&\ge \tilde{E}_{n,\beta_n,K_n} \{
|S_n/n^{1-\kappa}+W_n/n^{1/2-\kappa}-n^\kappa m_n| \cdot
1_{A_n(\delta_2)} \} - 2n^\kappa e^{-cn^{2\zeta}} -
c_2n^{\kappa-1/2}e^{-c_3n^{2\zeta}}  \nonumber. \eea Thus, we have
for all sufficiently large $n$ \bea \lefteqn{
\tilde{E}_{n,\beta_n,K_n}
\{|S_n/n^{1-\kappa}+W_n/n^{1/2-\kappa}-n^\kappa m_n| \cdot
1_{A_n(\delta_1)} \} +2n^\kappa e^{-cn^{2\zeta}} +c_2n^{\kappa-1/2}e^{-c_3n^{2\zeta}}   } \nonumber \\
&&\ge \tilde{E}_{n,\beta_n,K_n} \{
|S_n/n^{1-\kappa}+W_n/n^{1/2-\kappa}-n^\kappa m_n| \cdot
1_{ \{ S_n/n>\delta m_n \} } \}  \nonumber \\
&&\ge \tilde{E}_{n,\beta_n,K_n} \{
|S_n/n^{1-\kappa}+W_n/n^{1/2-\kappa}-n^\kappa m_n| \cdot
1_{A_n(\delta_2)} \} - 2n^\kappa e^{-cn^{2\zeta}} -
c_2n^{\kappa-1/2}e^{-c_3n^{2\zeta}}. \nonumber \eea
\fi
The last two integrals in (\ref{eqn:TwoEventAbsolute2}) coincide with the last two integrals in
(\ref{eqn:TwoEventAbsolute1}) and hence can be bounded the same way. For all sufficiently large $n$ this yields
\bea \lefteqn{ \tilde{E}_{n,\beta_n,K_n} \{
|S_n/n^{1-\kappa}+W_n/n^{1/2-\kappa}-n^\kappa m_n| \cdot 1_{ \{
S_n/n>\delta m_n \} }
\}      } \nonumber \\
&&\ge \tilde{E}_{n,\beta_n,K_n} \{
|S_n/n^{1-\kappa}+W_n/n^{1/2-\kappa}-n^\kappa m_n| \cdot
1_{A_n(\delta_2)} \} \\
&& \hspace{.15in} - \, 2n^\kappa e^{-cn^{2\zeta}} -
c_2n^{\kappa-1/2}e^{-c n^{2\zeta}/2}  \nonumber, \eea
where $c > 0$ and $c_2 > 0$ are constants.
In combination with (\ref{eqn:again}), the last inequality yields part (b).

In order to complete the proofs of parts (a) and (b) we now turn to the proofs of Claims 1 and 2.

\skp{\noindent}{\bf Proof of Claim 1.} We write
\bea
\{ S_n/n
>\delta m_n \}
& = & ( \{ S_n/n >\delta m_n \} \cap  \{ |W_n| \le \half n^\zeta \} )
                 \cup ( \{ S_n/n >\delta m_n \} \cap  \{ |W_n| > \half n^\zeta \} ) \nonumber \\
& \subset & ( \{ S_n/n >\delta m_n \} \cap  \{ |W_n| \le \half n^\zeta \} )
                 \cup \{ |W_n| > \half n^\zeta \} . \nonumber
\eea
Claim 1 follows if we prove for all sufficiently large $n$
\bea
\label{eqn:claim1}
\{ S_n/n
>\delta m_n \} \cap  \{ |W_n| \le \half n^\zeta \}  \subset A_n(\delta_1)
= \{S_n/n+W_n/n^{1/2}
> \delta_1 m_n \}. \eea
We have
\bea
\{ S_n/n >\delta m_n \} \cap \{
|W_n| \le \half n^\zeta \}
&=&
(\{
S_n/n >\delta m_n \} \cap \{ 0 \le W_n \le \half n^\zeta \} ) \nonumber \\
&&\cup
(\{ S_n/n >\delta m_n \} \cap \{ -\half n^\zeta \le W_n <0 \} )
\nonumber.\eea
If $S_n/n >\delta m_n$ and $0 \le W_n \le \half n^\zeta$,
then $S_n/n+W_n/n^{1/2} \ge S_n/n >\delta m_n >\delta_1
 m_n$. Thus
\be
\label{eqn:claim1A}
\{ S_n/n >\delta m_n \} \cap \{ 0 \le W_n
\le \half n^\zeta \} \subset A_n(\delta_1). \ee
Now assume that $S_n/n
>\delta m_n$ and $-\half n^\zeta \le W_n <0$. Since $\zeta < \frac{1}{2}-
\theta\alpha$, we have for all sufficiently large $n$
\[ (\delta-\delta_1)\bar{x} > n^{\zeta - 1/2 + \theta\alpha}.
\]
Since $\lim_{n \goto \infty} n^{\theta\alpha}m_n = \bar{x}$ [Thm.\ \ref{thm:3.1}], it
follows that for all sufficiently large $n$
\[ (\delta-\delta_1)m_n > \half n^{\zeta - 1/2}.
\]
Thus $\delta m_n - \half n^{\zeta-1/2} > \delta_1 m_n$ for all
sufficiently large $n$. Hence, if $S_n/n > \delta m_n$ and $W_n/n^{1/2}
\ge - \half n^{\zeta-1/2}$, then for all sufficiently large $n$
\[S_n/n + W_n/n^{1/2} > \delta m_n + W_n/n^{1/2} \ge
\delta m_n - \half n^{\zeta-1/2} > \delta_1 m_n.
\]
It follows that for all sufficiently large $n$
\[
\iffalse \label{eqn:claim1B} \fi
\{ S_n/n
>\delta m_n \} \cap \{ -n^{\zeta} \le W_n <0 \} \subset
A_n(\delta_1). \] Therefore (\ref{eqn:claim1}) follows from
(\ref{eqn:claim1A}) and the last display. This completes the
proof of Claim 1.

\skp{\noindent}{\bf Proof of Claim 2.} It suffices to prove that $A_n(\delta_2)
\subset \{S_n/n >\delta m_n \} \cup \{ |W_n| > \half n^\zeta \}$. We
write
\bea
A_n(\delta_2)
&=&
\{S_n/n+W_n/n^{1/2} > \delta_2 m_n \} \nonumber \\
&=&
(\{S_n/n+W_n/n^{1/2} > \delta_2 m_n \} \cap \{ |W_n| > \half n^\zeta
\}) \nonumber \\
&&\cup
( \{S_n/n+W_n/n^{1/2} > \delta_2 m_n \}
\cap \{ |W_n| \le \half n^\zeta\})
 \nonumber. \\
& \subset & \{ |W_n| > \half n^\zeta\} \cup
( \{S_n/n+W_n/n^{1/2} > \delta_2 m_n \}
\cap \{ |W_n| \le \half n^\zeta\}. \nonumber
\eea
Hence Claim 2 follows if we prove for all sufficiently large $n$
\[
\iffalse \label{eqn:claim2} \fi
\{S_n/n+W_n/n^{1/2} > \delta_2 m_n \} \cap \{
|W_n| \le \half n^\zeta \} \subset   \{ S_n/n >\delta m_n \}.
\]
We omit the proof, which is similar to the proof of (\ref{eqn:claim1}). The proof of Lemma \ref{lem:IneqOFtwoEvent}
is complete. $\Box$
\skp We now prove Lemma \ref{lem:G(xmn)miusG(mn)}, completing the preparatory lemmas that will be used
in the next section to prove part (b) of Theorem \ref{thm:6.1}.

\skp

{\noindent}{\bf Proof of Lemma \ref{lem:G(xmn)miusG(mn)}.} This is
done when $g$ has degree 4. We omit the analogous but more complicated proof when $g$ has
degree 6.
\iffalse
We recall that for $\alpha \in (0, \alpha_0)$
$u=1-\alpha/\alpha_0$, $\gamma=\theta\alpha$, and
$\kappa=\frac{1}{2}(1-\alpha/\alpha_0)+\theta\alpha$.
\fi
\iffalse \skp \noi {\it $g$ has degree 4.} \fi

\skp
\noi
{\bf Proof of part (a) of Lemma \ref{lem:G(xmn)miusG(mn)} when \boldmath $g$ \unboldmath has degree 4.}
By Taylor's theorem,
for any $R>0$, all $n \in \mathbb{N}$, and all $x \in \mathbb{R}$
satisfying $|x/n^{\kappa}| < R$ \bea \lefteqn{
nG_{\beta_n,K_n}(x/n^\kappa+m_n) -nG_{\beta_n,K_n}(m_n) }
\nonumber \\
&=& \sum_{j=1}^4 \frac{1}{n^{j\kappa-1}} \cdot
\frac{G_{\beta_n,K_n}^{(j)}(m_n)}{j!} x^j + \frac{1}{n^{5\kappa-1}}
\cdot \frac{G_{\beta_n,K_n}^{(5)}(m_n + \tau x/n^\kappa)}{5!} x^5
\nonumber, \eea where $\tau$ is a number in [0,1]. The quantity $m_n
+ \tau x/n^\kappa$ lies in the interval $[ \ m_n-|x/n^\kappa| ,
m_n+|x/n^\kappa| \ ]$. Since $m_n \in (0,1)$, $m_n \goto 0$ and $|
x/n^\kappa |<R$, we have $m_n+\tau x/n^\kappa \in (-R,R+1)$ for all
$n$. Since the sequence $(\beta_n,K_n)$ is bounded and
positive, there exists $a \in (0, \infty)$ such that $0 \le \beta_n
\le a$ and $0 \leq K_n \le a$ for all $n$. As a continuous function of
$(\beta, K, y)$ on the compact set $[  0, a  ]\times[ 0,a  ]\times[
 -R, R+1  ]$, it follows that $G_{\beta,K}^{(5)}(y)$ is uniformly
bounded. Since $m_n+\tau x/n^\kappa \in (-R,R+1)$ for all $n \in
\mathbb{N}$ and $x \in \mathbb{R}$ satisfying $|x/n^\kappa|<R$,
$G_{\beta_n,K_n}^{(5)}(m_n + \tau x/n^\kappa)$ is uniformly bounded
for $n \in \N$ and $x \in (-Rn^\kappa, Rn^\kappa)$. We summarize the last display
by writing \bea \lefteqn{ nG_{\beta_n,K_n}(x/n^\kappa+m_n)
-nG_{\beta_n,K_n}(m_n) }
\nonumber \\
&=& \sum_{j=1}^4 \frac{1}{n^{j\kappa-1}} \cdot
\frac{G_{\beta_n,K_n}^{(j)}(m_n)}{j!} x^j + \mbox{O} \!
\left(\frac{1}{n^{5\kappa-1}}\right) x^5 \nonumber, \eea where the
big-oh term is uniform for $x \in (-Rn^\kappa, Rn^\kappa)$.

Let $\varepsilon_n$ denote a sequence that converges to 0 and that
represents the various error terms arising in the proof. The
same notation $\varepsilon_n$ will be used to represent different error terms. To simplify the arithmetic,
we introduce $u = 1 - \alpha/\alpha_0 > 0$. We
have the following three properties:

\begin{itemize}
\item [ (1)] Since $m_n$ is the unique positive, global minimum point of
$G_{\beta_n,K_n}$, $G_{\beta_n,K_n}^{(1)}(m_n)=0$.

\item [(2)] By hypothesis (iii$^\prime$) of Theorem \ref{thm:main}, for $j=2,3,4$, we have
\[ G_{\beta_n,K_n}^{(j)}(m_n)= (g^{(j)}(\bar{x})+\varepsilon_n)/n^{\alpha/\alpha_0-j\theta\alpha} =
(g^{(j)}(\bar{x})+\varepsilon_n)/n^{1-u-j\theta\alpha}.
\]

\item [(3)] Since $\kappa=\frac{1}{2}u+\theta\alpha$, we have $j\kappa -u -j\theta\alpha = \left(\frac{j}{2}-1\right)\!u$
\iffalse
\[ j\kappa -u -j\theta\alpha =
\frac{j}{2}u+j\gamma-u-j\theta\alpha=\left(\frac{j}{2}-1\right)u
\]
\fi
for $j= 2, 3, 4$.
\end{itemize}

Using these properties, we obtain the following asymptotic formula,
which is valid for any $R>0$, all $n \in \mathbb{N}$, and all $x \in
\mathbb{R}$ satisfying $|x/n^\kappa|<R$: \bea
nG_{\beta_n,K_n}(x/n^\kappa + m_n) -nG_{\beta_n,K_n}(m_n)
&=&  \sum_{j=2}^4  \frac{1}{n^{j\kappa -1}} \cdot  \frac{g^{(j)}(\bar{x})+\varepsilon_n}{n^{1-u-j\theta\alpha}
\cdot j!}  \cdot x^j  + \mbox{O}\!\left(\frac{1}{n^{5\kappa-1}}\right) x^5         \nonumber \\
&=&  \sum_{j=2}^4  \frac{1}{j!} \cdot  \frac{g^{(j)}(\bar{x})+\varepsilon_n}{n^{(j/2-1)u}}  \cdot x^j  + \mbox{O}\!\left(\frac{1}{n^{5\kappa-1}}\right) x^5         \nonumber \\
&=&  \frac{1}{2!} \cdot (g^{(2)}(\bar{x})+\varepsilon_n) x^2     +
     \frac{1}{3!} \cdot \frac{(g^{(3)}(\bar{x})+\varepsilon_n)}{n^{u/2}} x^3
     \nonumber \\
&&+  \  \frac{1}{4!} \cdot
\frac{(g^{(4)}(\bar{x})+\varepsilon_n)}{n^u} x^4  +
\mbox{O}\!\left(\frac{1}{n^{5\kappa-1}}\right) x^5    \nonumber.
\eea

\iffalse    {\noindent}{\it Proof of part (a) of Lemma {\em
\ref{lem:G(xmn)miusG(mn)}} when $g$ has degree 4.}           \fi
By hypothesis (iii$^\prime$) of Theorem \ref{thm:main} and part (c)
of Theorem \ref{thm:6.1}, we have $1/4 \le
\theta\alpha_0<\kappa<1/2$. Therefore $5\kappa-1 > 5 \theta\alpha_0
-1 >0$. Since $u>0$ and $\varepsilon_n \goto 0$, we have for all $x
\in \mathbb{R}$
\[ \lim_{n \goto \infty} (nG_{\beta_n,K_n}(x/n^\kappa + m_n)
-nG_{\beta_n,K_n}(m_n)) =\frac{1}{2} g^{(2)}(\bar{x})x^2.
\]
This completes the proof of part (a) of Lemma
\ref{lem:G(xmn)miusG(mn)} when $g$ has degree 4.

\skp

\iffalse   {\noindent}{\it Proof of part (b) of Lemma {\em
\ref{lem:G(xmn)miusG(mn)}} when \boldmath $g$ \unboldmath has degree 4.}      \fi

\noi
{\bf Proof of part (b) of Lemma \ref{lem:G(xmn)miusG(mn)} when $g$ has degree 4.}
Hypothesis (iii$^\prime$) of Theorem \ref{thm:main} states that
$g^{(j)}(\bar{x}) >0$ for  $j=2, 3, 4$. It follows that for all
sufficiently large $n$ and all $x \in \mathbb{R}$
\[\frac{1}{2!} \cdot (g^{(2)}(\bar{x})+\varepsilon_n) x^2
\ge \frac{1}{2 \cdot 2!} \cdot g^{(2)}(\bar{x}) x^2
\]
and
\[\frac{1}{4!} \cdot \frac{(g^{(4)}(\bar{x})+\varepsilon_n)}{n^u} x^4
\ge \frac{1}{2 \cdot 4!} \cdot \frac{g^{(4)}(\bar{x})}{n^u} x^4
\]
and that for all sufficiently large $n$
\[\frac{1}{3!} \cdot \frac{g^{(3)}(\bar{x})+\varepsilon_n}{n^{u/2}} x^3
\ge \frac{1}{2 \cdot 3!} \cdot \frac{g^{(3)}(\bar{x})}{n^{u/2}} x^3
\ \ \mbox{for all} \ \ x \ge 0
\]
and
\[\frac{1}{3!} \cdot \frac{g^{(3)}(\bar{x})+\varepsilon_n}{n^{u/2}} x^3
\ge \frac{2}{3!} \cdot \frac{g^{(3)}(\bar{x})}{n^{u/2}} x^3 \ \
\mbox{for all} \ \ x<0.
\]
\iffalse and
\[O(\frac{1}{n^{5\kappa-1}}) x^5 = \frac{1}{n^u}  \cdot O(\frac{1}{n^{5\kappa-1-u}})
x^5.
\]\fi

We first consider $x \in [0, Rn^\kappa)$. Since
$g^{(4)}(\bar{x})>0$, for all sufficiently large $n$ and all such
$x$ we have
\bea \label{eqn:x+GGg4}
\lefteqn{ nG_{\beta_n,K_n}(x/n^\kappa +
m_n) -nG_{\beta_n,K_n}(m_n) } \\
&\ge& \frac{1}{2 \cdot 2!} \cdot
g^{(2)}(\bar{x}) x^2 +
     \frac{1}{2 \cdot 3!} \cdot \frac{g^{(3)}(\bar{x})}{n^{u/2}} x^3  \nonumber  \\
&&+ \  \frac{1}{2 \cdot 4!} \cdot \frac{g^{(4)}(\bar{x})}{n^u} x^4 +
\frac{1}{n^u} \mbox{O} \! \left( \frac{x}{n^{5\kappa-1-u}}\right)x^4  \nonumber \\
&\ge& \frac{1}{2 \cdot 2!} \cdot g^{(2)}(\bar{x}) x^2 +
     \frac{1}{2 \cdot 3!} \cdot \frac{g^{(3)}(\bar{x})}{n^{u/2}} x^3     \nonumber \\
&&+ \  \frac{1}{2 \cdot 4!} \cdot \frac{g^{(4)}(\bar{x})}{n^u} x^4
(1  + \mbox{O}(x/n^{5\kappa-1-u})). \nonumber \eea %%%%%%%%%%%%%%%%%%%%%%%%%%
By hypothesis (iii$'$) of Theorem \ref{thm:main}, $\theta\alpha_0 \in [1/4,1/2)$. Hence
$4\theta-1/\alpha_0
\ge 0$, and so
\iffalse
\bea 5\kappa-1-u-\kappa &=& 4\kappa-1-u =
4(u/2+\gamma)-1-u = 2u+4\gamma-1-u \nonumber \\
&=& 4\gamma +u-1 = 4\theta\alpha-\alpha/\alpha_0 = (4\theta
-1/\alpha_0)\alpha \ge 0. \nonumber \eea
\fi
\[
5\kappa-1-u = \kappa + (4\theta
-1/\alpha_0)\alpha \ge \kappa.
\]
\iffalse
Thus $5\kappa-1-u \ge
\kappa$, and so
\fi
Hence for all $0<x<Rn^\kappa$ we have $0 \le
x/n^{5\kappa-1-u} \le x/n^\kappa <R$. Thus the term
$\mbox{O}(x/n^{5\kappa-1-u})$ appearing in (\ref{eqn:x+GGg4}) can be
made larger than $-1$ for all $0 \le x/n^\kappa <R$ by choosing $R$
to be sufficiently small. Since $g^{(3)}(\bar{x}) >0$,
$g^{(4)}(\bar{x})>0$, and $1+\mbox{O}(x/n^{5\kappa-1-u}) >0$, we
have that for all sufficiently large $n$ and all $x
\in[0,Rn^\kappa)$
\[nG_{\beta_n,K_n}( x/n^\kappa + m_n) -nG_{\beta_n,K_n}(m_n)
\ge \frac{1}{2 \cdot 2!} g^{(2)}(\bar{x})x^2 \ge \frac{1}{8}
g^{(2)}(\bar{x})x^2.
\]
This is the conclusion of part (b) of Lemma
\ref{lem:G(xmn)miusG(mn)} for all $0 \le x <Rn^\kappa$ when $g$ has
degree 4.

We now consider $x \in (-Rn^\kappa,0]$. Since
$g^{(4)}(\bar{x})>0$, for all sufficiently large $n$ and all such
$x$ we have \bea \label{eqn:x-GGg4} nG_{\beta_n,K_n}(x/n^\kappa +
m_n) -nG_{\beta_n,K_n}(m_n) &\ge& \frac{1}{2 \cdot 2!} \cdot
g^{(2)}(\bar{x}) x^2 +
     \frac{2}{ 3!} \cdot \frac{g^{(3)}(\bar{x})}{n^{u/2}} x^3      \\
&&+  \  \frac{1}{2 \cdot 4!} \cdot \frac{g^{(4)}(\bar{x})}{n^u} x^4
+
\frac{1}{n^u} \mbox{O} \! \left( \frac{x}{n^{5\kappa-1-u}} \right)x^4 \nonumber \\
&\ge& \frac{1}{2 \cdot 2!} \cdot g^{(2)}(\bar{x}) x^2 +
     \frac{2}{ 3!} \cdot \frac{g^{(3)}(\bar{x})}{n^{u/2}} x^3     \nonumber \\
&&+ \  \frac{1}{2 \cdot 4!} \cdot \frac{g^{(4)}(\bar{x})}{n^u} x^4
(1  + \mbox{O}(x/n^{5\kappa-1-u})). \nonumber \eea %%%%%%%%%%%%%%%%%%%%%%%%%%
Since $5\kappa-1-u \ge  \kappa$, for all $-Rn^\kappa<x<0$ we have
$-R < x/n^\kappa \le x/n^{5\kappa-1-u} <0$. Thus the term
$\mbox{O}(x/n^{5\kappa-1-u})$ appearing in (\ref{eqn:x-GGg4}) can be
made larger than $-1$ for all $-R<x/n^\kappa <0$ by choosing $R$ to
be sufficiently small. Since $g^{(4)}(\bar{x})>0$ and
$1+\mbox{O}(x/n^{5\kappa-1-u}) >0$, we have that for all
sufficiently large $n$ and all $x \in (-Rn^\kappa,0)$
\[nG_{\beta_n,K_n}( x/n^\kappa + m_n) -nG_{\beta_n,K_n}(m_n)
\ge \frac{1}{2 \cdot 2!} g^{(2)}(\bar{x})x^2 + \frac{2}{ 3!} \cdot
\frac{g^{(3)}(\bar{x})}{n^{u/2}} x^3.
\]
By Theorem \ref{thm:3.1} we have $m_n  \sim \bar{x}/n^{\theta\alpha}$.
Thus $n^{\theta\alpha} m_n =
\bar{x}+\varepsilon_n$, and
\[
n^\kappa m_n = n^{u/2} \cdot n^{\theta\alpha}
m_n = n^{u/2}(\bar{x}+\varepsilon_n).
\]
In part (b) of Lemma
\ref{lem:G(xmn)miusG(mn)} we assume that $x/n^\kappa >
-(1-\bar{\delta})m_n$ and $0<\bar{\delta}<1$. Thus for all
sufficiently large $n$ and all such $x$
\[ x > -(1-\bar{\delta})n^\kappa m_n
= -(1-\bar{\delta})n^{u/2}(\bar{x}+\varepsilon_n) \ge
-(1-\bar{\delta})n^{u/2} \cdot 2\bar{x}.
\]
Since $g^{(3)}(\bar{x})>0$, we see that for all sufficiently large
$n$, all $x \in (-Rn^\kappa,0)$, and all $x/n^\kappa >
-(1-\bar{\delta})m_n$ there exists $\Delta \in (0,1)$ such that for
any $\bar{\delta} \in (\Delta, 1)$ the following inequalities hold:
\bea \lefteqn{
nG_{\beta_n,K_n}(x/n^\kappa + m_n) -nG_{\beta_n,K_n}(m_n) } \nonumber \\
&\ge& \frac{1}{2 \cdot 2!} \cdot g^{(2)}(\bar{x}) x^2
+ \frac{2}{ 3!} \cdot \frac{g^{(3)}(\bar{x})}{n^{u/2}}
     x^2 \cdot [-(1-\bar{\delta})n^{u/2} \cdot 2\bar{x}] \nonumber \\
&=& \left( \ \frac{1}{2 \cdot 2!} \cdot g^{(2)}(\bar{x}) - \frac{2
\cdot 2}{ 3!} \cdot g^{(3)}(\bar{x})(1-\bar{\delta})\bar{x}  \
\right) x^2  \nonumber \\
&\ge& \frac{1}{2} \cdot  \frac{1}{2 \cdot 2!}
g^{(2)}(\bar{x})x^2 = \frac{1}{8} g^{(2)}(\bar{x})x^2\nonumber. \eea This is the
conclusion of part (b) of Lemma \ref{lem:G(xmn)miusG(mn)} for
$-Rn^\kappa <x<0$ and $x/n^\kappa >-(1-\bar{\delta})m_n$ when $g$
has degree 4.

\iffalse
Since $m_n \goto 0$, for all sufficiently large $n$ and all
$0<\bar{\delta}<1$ we have $-Rn^\kappa < -(1-\bar{\delta})n^\kappa
m_n$, and thus for $x<0$ the condition $x/n^\kappa >
-(1-\bar{\delta})m_n$ is more restrictive than $x>-Rn^\kappa$.
However, it is useful to use the latter condition to prove that
$1+\mbox{O}(x/n^{5\kappa-1-u}) >0$.
\fi

We have shown that for any $\bar{\delta} \in (\Delta, 1)$ there
exists $R>0$ such that for all sufficiently large $n$ and all $x \in
\mathbb{R}$ satisfying $|x/n^\kappa|<R$ and $x/n^\kappa > -(1-
\bar{\delta})m_n$
\[nG_{\beta_n,K_n}( x/n^\kappa + m_n) -nG_{\beta_n,K_n}(m_n)
\ge \frac{1}{8} g^{(2)}(\bar{x})x^2.
\]
This completes the proof of part (b) of Lemma
\ref{lem:G(xmn)miusG(mn)} when $g$ has degree 4. Because we are omitting the proof of part (b) when $g$ has degree 6,
the proof of Lemma
\ref{lem:G(xmn)miusG(mn)} is complete. $\Box$

\skp \iffalse Having completed the logic of proof of part (b) of
Theorem \ref{thm:6.1} and the statement and proof of preparatory
lemmas, \fi We have completed the statements and proofs of the
preparatory lemmas. We now turn to the proof of part (b) of Theorem
\ref{thm:6.1}.

%===================== new section =====================================

\section{Proof of Part (b) of Theorem \ref{thm:6.1}}
\label{section:Thm61b}
\beginsec

As we saw at the start of section \ref{section:LogicANDPreLemmas}, the proof
of part (b) of Theorem \ref{thm:6.1} involves two steps. Step 1 is
proved in the next subsection. Step 2 is subdivided into two
substeps. Substep 2a is proved in subsection
\ref{subsection:step2aOFpartB}, and Substep 2b is proved in
subsection \ref{subsection:step2bOFpartB}.

%****************************************************************************************************************************************

\subsection{Proof of Step 1 in Proof of Theorem \ref{thm:6.1} (b)}
\label{subsection:step1OFpartB}

Part (b) of the following lemma states Step 1 in the proof of
part (b) of Theorem \ref{thm:6.1}. \iffalse The proof of part (a) of
Lemma \ref{lem:Step1OFThm61b} requires Lemma
\ref{lem:Thm61aStep1}(a), \ref{lem:AsyOFGmn},
\ref{lem:G(xmn)miusG(mn)}, \ref{lem:intRnkexpIneq}, Proposition
\ref{prop:weakinteg} and Lemma 4.4(a) in \cite{CosEllOtt}. The proof
of part (b) of Lemma \ref{lem:Step1OFThm61b} needs lemma
\ref{lem:Thm61aStep1}(a), Lemma \ref{lem:tEgammaLowerValue}(b) with
$\bar{\gamma}= \kappa$ and part (a) of this lemma. \fi We recall
that $W_n$ is a sequence of normal random variables with mean 0 and
variance $(2\beta_n K_n)^{-1}$ defined on a probability space
$(\Omega,\mathcal {F},Q)$. We denote by $\tilde{E}_{n, \beta_n,K_n}$
expectation with respect to the product measure $P_{n,\beta_n,K_n}
\times Q$; $P_{n,\beta_n,K_n}$ is defined in (\ref{eqn:P})--(\ref{eqn:PnbetaK}).

\begin{lem}
\label{lem:Step1OFThm61b} We assume that $(\beta_n, K_n)$ satisfies
the hypotheses of Theorem {\em \ref{thm:main}} for all
$0<\alpha<\alpha_0$. For $\bar{\delta} \in (0,1)$ define
\[ A_n(\bar{\delta})=\{ S_n/n + W_n/n^{1/2} > \bar{\delta}m_n \}
\]
where $m_n=m(\beta_n, K_n)$. Let $\Delta \in (0,1)$ be the number
determined in part {\em (b)} of Lemma {\em \ref{lem:G(xmn)miusG(mn)}}. Then for any
$0<\alpha<\alpha_0$ and any $\bar{\delta} \in (\Delta, 1)$, the
following conclusions hold.

\begin{itemize}
\item [\em (a)] We have the limit
\bea
\lefteqn{ \lim_{n \goto \infty}
\int_{-n^\kappa(1-\bar{\delta})m_n}^{\infty} |x|
\exp[-nG_{\beta_n,K_n}(x/n^\kappa+m_n)+ nG_{\beta_n,K_n}(m_n)]dx }
\nonumber \\
&&=
\int_{\mathbb{R}} |x|\exp[\ts
-\frac{1}{2}g^{(2)}(\bar{x}) x^2]dx. \nonumber \ \ \ \ \ \ \ \ \ \ \
\ \ \  \ \ \ \ \ \ \ \ \ \ \ \ \ \ \ \ \ \ \ \ \ \ \ \ \ \ \ \ \ \
\eea

\item [\em (b)] We have the limit
\bea
\lefteqn{ \lim_{n \goto \infty} \tilde{E}_{n,\beta_n,K_n} \{
|S_n/n^{1-\kappa}+ W_n/n^{1/2-\kappa} -n^{\kappa} m_n | \
\big| \  A_n(\bar{\delta})\}   } \nonumber  \\
&=&
\frac{1}{\int_{\mathbb{R}}
\exp[-\frac{1}{2}g^{(2)}(\bar{x})x^2]dx} \nonumber \cdot
\int_{\mathbb{R}} |x| \exp[\ts -\frac{1}{2}g^{(2)}(\bar{x})x^2]dx
\nonumber \\
&=& \bar{z} \nonumber. \eea
\end{itemize}

\end{lem}

{\noindent}{\bf Proof of part (a) of Lemma \ref{lem:Step1OFThm61b}.}
Let $\Psi_n$ and $\Psi$ denote the measures on $\mathbb{R}$ defined
by
\[ \Psi_n(dx) = 1_{(-n^\kappa(1-\bar{\delta})m_n,\infty)}(x) \cdot
\exp[-nG_{\beta_n,K_n}(x/n^\kappa+m_n)+nG_{\beta_n,K_n}(m_n)]dx
\]
and
\[ \Psi(dx) = \exp[\ts -\frac{1}{2} g^{(2)}(\bar{x})x^2]dx.
\]
According to Lemma \ref{lem:Thm61aStep1}, $\Psi_n$ converges weakly to $\Psi$.
\iffalse
We work with these measures rather than with the normalized
probability measures because in terms of $\Psi_n$ and $\Psi$ the
limit in Lemma \ref{lem:Thm61aStep1} can be directly expressed;
namely, for any bounded, continuous function $f$
\[ \lim_{n \goto \infty} \int_{\mathbb{R}}f d\Psi_n=\int_{\mathbb{R}}f
d\Psi.
\]
\fi
The limit in part (a) of Lemma \ref{lem:Step1OFThm61b} can
be expressed as
\[ \lim_{n \goto \infty} \int_{\mathbb{R}} |x| d\Psi_n = \int_{\mathbb{R}} |x|
d\Psi.
\]
As discussed in Theorem 4 in \S II.6 of \cite{Shiryaev}, this limit would follow from the weak convergence of $\Psi_n$ to
$\Psi$ if one could prove the
uniform integrability estimate
\[ \lim_{j \goto \infty} \sup_{n \in \mathbb{R}} \int_{\{ |x|>j \}}
|x| d \Psi_n =0.
\]
\iffalse Under the assumption that $\int_{\mathbb{R}} |x|d \Psi <
\infty$, \fi The next proposition shows that the limit $\lim_{n
\goto \infty} \int_{\mathbb{R}} |x| d\Psi_n  = \int_{\mathbb{R}} |x|
d\Psi$ is a consequence of a condition that is weaker than uniform
integrability.
\iffalse
It is proved in Proposition 8.3 in \cite{EllMacOtt2010}.
\fi

\begin{prop}
\label{prop:weakinteg} Let $\Psi_n$ be a sequence of measures on
$\mathbb{R}$ that converges weakly to a measure $\Psi$ on
$\mathbb{R}$.
\iffalse
; i.e., for any bounded, continuous function $\tilde{f}$
\[ \lim_{n \goto \infty} \int_{\mathbb{R}} \tilde{f} d \Psi_n = \int_{\mathbb{R}} \tilde{f} d
\Psi.
\]
\fi
Assume in addition that $\int_{\mathbb{R}} |x| d \Psi < \infty$ and
that
\[ \lim_{j \goto \infty} \limsup_{n \goto \infty} \int_{ \{ |x|>j \} } |x| d
\Psi_n = 0.
\]
It then follows that
\[ \lim_{n \goto \infty} \int_{\mathbb{R}} |x| d \Psi_n = \int_{\mathbb{R}} |x| d
\Psi.
\]
\end{prop}

\noi
{\bf Proof.} Since $\Psi_n \Longrightarrow \Psi$, we have $\Psi_n(\R) \goto \Psi(\R)$. Hence the proposition
is a consequence of Proposition 8.3 in \cite{EllMacOtt2010} applied to the sequence of probability
measures
\[
\frac{1}{\Psi_n(\R)} \cdot \Psi_n(dx) \Longrightarrow \frac{1}{\Psi(\R)} \cdot \Psi(dx).
\]
This completes the proof. $\Box$.

\iffalse {\noindent}{\bf Proof.} For $j \in \mathbb{N}$, $f_j$
denotes the bounded, continuous function that equals $|x|$ for
$|x|\le j$ and
equals $j$ for $|x|>j$. Then
\bea
\lefteqn{ \left|\int_{\mathbb{R}}|x|d \Psi_n  - \int_{\mathbb{R}}|x|d\Psi  \right|  } \nonumber \\
&&\le
\int_{\mathbb{R}} ||x|-f_j| d \Psi_n
+ \left| \int_{\mathbb{R}}f_j d\Psi_n - \int_{\mathbb{R}}f_j d\Psi  \right|
+  \int_{\mathbb{R}} ||x|-f_j| d \Psi \nonumber \\
&&\le
2\int_{\{|x|>j\}} |x| d \Psi_n
+
\left| \int_{\mathbb{R}}f_j
d\Psi_n - \int_{\mathbb{R}}f_j d\Psi  \right|
+
2\int_{\{|x|>j\}}
|x| d \Psi \nonumber. \eea
Since $\Psi_n  \Longrightarrow \Psi$, we
have $\int_{\mathbb{R}}f_j d\Psi_n  \goto \int_{\mathbb{R}}f_j
d\Psi$, and therefore
\bea
\lefteqn{ \limsup_{n \goto \infty}
\left| \int_{\mathbb{R}}|x|d \Psi_n  - \int_{\mathbb{R}}|x|d\Psi
\right| }  \nonumber  \\
&&\le
2\limsup_{n \goto \infty}\int_{\{|x|>j\}} |x| d \Psi_n +
2\int_{\{|x|>j\}} |x| d \Psi \nonumber. \eea
By sending $j \goto
\infty$, both terms on the right side of above inequality have
limits 0.  $\Box$ \fi

\skp We now verify the following hypotheses of Proposition
\ref{prop:weakinteg} for the measures $\Psi_n$ and $\Psi$:
\begin{itemize}
\item [ (1)]  $\Psi_n  \Longrightarrow \Psi$.
\item [  (2)]  $ \int_{\mathbb{R}} |x| d
\Psi < \infty$.
\item [ (3)]  $\lim_{j \goto \infty} \limsup_{n \goto \infty} \int_{ \{ |x|>j \} } |x| d
\Psi_n = 0$.
\end{itemize}

\iffalse
\skp $\Psi_n$ and $\Psi$ are not probability measures. However, by
Lemma \ref{lem:Thm61aStep1} with $f \equiv 1$ \bea \lim_{n \goto
\infty} \Psi_n(\mathbb{R}) &=& \lim_{n \goto \infty}
\int_{-n^\kappa(1-\bar{\delta})m_n}^{\infty}
\exp[-nG_{\beta_n,K_n}(x/n^\kappa+m_n)+ nG_{\beta_n,K_n}(m_n)]dx
\nonumber\\
&=& \int_{\mathbb{R}} \exp[\ts -\frac{1}{2}g^{(2)}(\bar{x}) x^2]dx =
\Psi(\mathbb{R}). \nonumber
\eea
\iffalse Because of this limit, it
is elementary to show that
\fi
It follows that if $\lim_{n \goto
\infty} \int_{\mathbb{R}}f d\Psi_n=\int_{\mathbb{R}}f d\Psi$ for all
bounded, continuous functions $f$, then $\Psi_n
\Longrightarrow \Psi$ \iffalse       \cite[Remark 1, p.
313]{Shiryaev}            \fi. This is the weak convergence needed
in item (1). Item (2) follows from
\[ \int_{\mathbb{R}} |x| d \Psi = \int_{\mathbb{R}} |x| \exp[\ts -\frac{1}{2}
g^{(2)}(\bar{x})x^2]dx < \infty.
\]
\fi

Item (1) is proved in Lemma \ref{lem:Thm61aStep1}, and item (2) is immediate
from the definition of $\Psi$.
We now prove item (3).
Since
\[ \int_{\{|x| > j\}} |x| d \Psi_n = \int_{\{|x| > j\} \cap  \{ x >
-n^\kappa (1- \bar{\delta})m_n \}  } |x| d \Psi_n,
\]
we can prove item (3) by showing that
\[ \lim_{j \goto \infty} \limsup_{n \goto \infty} \int_{ \{ |x|>j \}  \cap \{ x >
-n^\kappa (1- \bar{\delta})m_n \}  } |x| d \Psi_n = 0.
\]
In order to do this we find, for any $j \in \mathbb{N}$ and all
sufficiently large $n$, quantities $A_j$, $B_n$ and $C_n$ with the
properties that
\[ \int_{ \{ |x|>j \}  \cap \{ x >
-n^\kappa (1- \bar{\delta})m_n \}  } |x| d \Psi_n \le A_j + B_n +
C_n,
\]
$A_j \goto 0$ as $j \goto \infty$, $B_n \goto 0$ as $ n \goto
\infty$, and $C_n \goto 0$ as $n \goto \infty$.
\iffalse    It
follows from these properties that
\[\lim_{j \goto \infty} \limsup_{n \goto \infty}\int_{ \{ |x|>j \}  \cap \{ x >
-n^\kappa (1- \bar{\delta})m_n \}  } |x| d \Psi_n \le \lim_{j \goto
\infty}A_j + \lim_{n \goto \infty}B_n +  \lim_{n \goto \infty}C_n=0.
\]
This display yields item (3).     \fi

We now specify the quantities $A_j$, $B_n$ and $C_n$. Given positive
integers $j$ and $n$, let $R$ and $K$ be positive numbers that
satisfy $K>R$ and that will be specified below. Then
\iffalse We then partition
the set $\{|x|>j\} \cap \{x/n^\kappa > -(1-\bar{\delta})m_n\}$ into
the following three subsets:
\fi
\bea \lefteqn{ \{|x|>j\} \cap
\{x/n^\kappa > -(1-\bar{\delta})m_n\}
        } \nonumber \\
&=&     [ \{|x| > j\} \cap \{|x/n^\kappa| < R \} \cap \{x/n^\kappa > -(1-\bar{\delta})m_n\} ]  \nonumber \\
&&\cup  [ \{|x| > j\} \cap \{R \le |x/n^\kappa| < K \} \cap \{x/n^\kappa > -(1-\bar{\delta})m_n\} ]  \nonumber \\
&&\cup  [ \{|x| > j\} \cap \{|x/n^\kappa| \ge K \} \cap \{x/n^\kappa > -(1-\bar{\delta})m_n\}]  \nonumber \\
\iffalse\eea
Since for all $n$
\bea
\lefteqn{ \{|x|>j\} \cap \{x/n^\kappa > -(1-\bar{\delta})m_n\}
        } \nonumber \\
\fi
&\subset&     [ \{|x| > j\} \cap \{|x/n^\kappa| < R \} \cap \{x/n^\kappa > -(1-\bar{\delta})m_n\} ]  \nonumber \\
&&\cup  [ \{R \le |x/n^\kappa| < K \} \cap \{x/n^\kappa > -(1-\bar{\delta})m_n\} ]
\cup  \{|x/n^\kappa| \ge K \} \nonumber. \eea
Since $m_n \goto 0$, for all sufficiently large $n$
\[
\{R \le |x/n^\kappa| < K \} \cap \{x/n^\kappa > -(1-\bar{\delta})m_n\} = \{R \le x/n^\kappa < K \}.
\]
Hence for all sufficiently large $n$
\bea
\label{eqn:Aj+Bn+Cn}
\lefteqn{ \int_{ \{|x|>j\} \cap
\{x/n^\kappa > -(1-\bar{\delta})m_n\}} |x| d \Psi_n } \\
&& \le  \int_{ \{
|x|>j \} \cap \{|x/n^\kappa|<R\} \cap
\{ x > -n^\kappa (1- \bar{\delta})m_n \}  } |x| d \Psi_n \nonumber \\
&& \hspace{.15in} +
\int_{ \{ R \le x/n^\kappa<K\} } |x| d \Psi_n
+
\int_{ \{
|x/n^\kappa| \ge K\} } |x| d \Psi_n . \nonumber \eea

We next estimate each of these three integrals. By part (b) of Lemma
\ref{lem:G(xmn)miusG(mn)}, there exists $\Delta \in (0,1)$ such that
for any $\bar{\delta} \in (\Delta, 1)$ there exists $R>0$ such that
for all sufficiently large n and all $x \in \mathbb{R}$ satisfying
$|x/n^\kappa| < R$ and $x/n^\kappa > -(1-\bar{\delta})m_n$
\[ nG_{\beta_n,K_n}( x/n^\kappa + m_n) -nG_{\beta_n,K_n}(m_n) \ge
H(x)= \frac{1}{8} g^{(2)}(\bar{x})x^2.
\]
\iffalse
where $H(x) \goto \infty$ as $|x| \goto \infty$.
\fi
Since $\exp[-H(x)]$
is integrable, for all sufficiently large $n$ we estimate the first
integral on the right hand side of equation (\ref{eqn:Aj+Bn+Cn}) by
\bea
\label{eqn: Aj}
\lefteqn{  \int_{ \{ |x|>j \} \cap
\{|x/n^\kappa|<R\} \cap
\{ x > -n^\kappa (1- \bar{\delta})m_n \}  } |x| d \Psi_n  }  \\
\nonumber &&\le
A_j = \int_{ \{|x|>j\}}|x| \cdot \exp[-H(x)]dx \goto 0   \
\mbox{ as }  \
 j \goto \infty.
\eea

By part (a) of Lemma 4.4 in \cite{CosEllOtt}, there exists $K>0$ and
$D_1>0$ such that $G_{\beta_n,K_n}(x) \ge D_1 x^2$ for all $|x| >K$.
Since $m_n \goto 0$, it follows that for all sufficiently large $n$
and all $x \in \mathbb{R}$ satisfying $|x/n^\kappa| \ge K$, there
exists $D>0$ such that
\[ nG_{\beta_n,K_n}(x/n^\kappa +m_n) \ge nD_1(x/n^\kappa +m_n)^2 \ge
nD(x/n^\kappa)^2.
\]
Without loss of generality $K$ can be chosen to be larger than the
quantity $R$ specified in the preceding paragraph. By Lemma
\ref{lem:AsyOFGmn}, for all sufficiently large $n$ there exists $\ve_n \goto 0$ such that
\[ G_{\beta_n,K_n}(m_n) =
\frac{g(\bar{x})+\varepsilon_n}{n^{\alpha/\alpha_0}} \le
\frac{g(\bar{x})}{2n^{\alpha/\alpha_0}} < 0.
\]
These bounds allow us to estimate the third integral on
the right hand side of equation (\ref{eqn:Aj+Bn+Cn}) by
\bea
\label{eqn: Cn}
\lefteqn{ \int_{ \{ |x/n^\kappa| \ge K\} } |x| d \Psi_n } \\
&& \le  \int_{ \{
|x/n^\kappa| \ge K\} } |x|
\exp[-nG_{\beta_n,K_n}(x/n^\kappa+m_n)+nG_{\beta_n,K_n}(m_n)]dx
\nonumber \\
&&\le
\int_{ \{ |x/n^\kappa| \ge K\} } |x| \exp[-n D(x/n^\kappa)^2
]dx \nonumber \\
\iffalse
&&=
C_n  \nonumber \\
&&=
\frac{1}{2}\int_{ \{ |x/n^\kappa| \ge K\} } \exp[-n
D|x/n^\kappa|^2 ]d|x|^2 \nonumber \\
\nonumber
\fi
&&\le
C_n = \frac{2}{D} \cdot n^{2\kappa -1} \exp[-nDK^2 ]  \goto  0
\  \mbox{ as }  \ n \goto \infty. \nonumber \eea

With these choices of $R$ and $K$, we use Lemma
\ref{lem:intRnkexpIneq} to estimate the second integral on the right hand
side of (\ref{eqn:Aj+Bn+Cn}). There exists $c_2 > 0$ such that for all sufficiently large $n$
\bea
\label{eqn: Bn}
\lefteqn{
\int_{
\{ R \le x/n^\kappa <K\} } |x| d \Psi_n } \\
&& =
\int_{ \{ R \le
x/n^\kappa < K\} } |x|
\exp[-nG_{\beta_n,K_n}(x/n^\kappa+m_n)+nG_{\beta_n,K_n}(m_n)]dx
\nonumber \\
&&\le
Kn^\kappa \int_{\{ x/n^\kappa \ge R \}}
\exp[-nG_{\beta_n,K_n}(x/n^\kappa+m_n)+nG_{\beta_n,K_n}(m_n)]dx
\nonumber \\
&&\le
B_n = Kn^\kappa \exp[-c_2 n]  \goto 0 \
\mbox{ as }  \  n \goto \infty  \nonumber . \eea

Together equations (\ref{eqn: Aj}), (\ref{eqn: Bn}) and (\ref{eqn:
Cn}) prove (\ref{eqn:Aj+Bn+Cn}). This completes the proof of part
(a) of Lemma \ref{lem:Step1OFThm61b}.

\skp {\noindent}{\bf Proof of part (b) of Lemma
\ref{lem:Step1OFThm61b}.} Part (b) of Lemma
\ref{lem:tEgammaLowerValue} states that
\bea \lefteqn{
\tilde{E}_{n,\beta_n,K_n} \{ |S_n/n^{1-\kappa}+ W_n/n^{1/2-\kappa}
-n^{\kappa} m_n | \ \big| \ A_n(\bar{\delta})\}
        } \nonumber \\
&&=  \frac{\int_{-n^\kappa(1-\bar{\delta})m_n}^{\infty}
|x|\exp[-nG_{\beta_n,K_n}(x/n^\kappa+m_n)+ nG_{\beta_n,K_n}(m_n)]dx
          }
          {\int_{-n^\kappa(1-\bar{\delta})m_n}^{\infty}
\exp[-nG_{\beta_n,K_n}(x/n^\kappa+m_n)+ nG_{\beta_n,K_n}(m_n)]dx
          }
\nonumber. \eea
Hence by part (a) of Lemma \ref{lem:Step1OFThm61b} and
Lemma \ref{lem:Thm61aStep1} for $f(x) \equiv 1$, the last integral
has the limit
\[ \frac{ \int_{\mathbb{R}} |x|
\exp[-\frac{1}{2}g^{(2)}(\bar{x})x^2]dx
        }
        { \int_{\mathbb{R}}
\exp[-\frac{1}{2}g^{(2)}(\bar{x})x^2]dx
        } = \bar{z}.
\]
This completes the proof of part (b) of Lemma \ref{lem:Step1OFThm61b} and hence the proof of the
lemma. $\Box$

\skp Having completed Step 1 in the proof of part (b) of Theorem
\ref{thm:6.1}, we now turn to Substep 2a.

\subsection{Proof of Substep 2a in Proof of Theorem \ref{thm:6.1} (b)}
\label{subsection:step2aOFpartB}

Lemma \ref{lem:substep2aOFThm61b} proves Substep 2a of part (b)
of Theorem \ref{thm:6.1}. We recall that $W_n$ is a sequence of
normal random variables with mean 0 and variance $(2\beta_n
K_n)^{-1}$ defined on a probability space $(\Omega,\mathcal {F},Q)$.
We denote by $\tilde{E}_{n, \beta_n,K_n}$ expectation with respect
to the product measure $P_{n,\beta_n,K_n} \times Q$.

\begin{lem}
\label{lem:substep2aOFThm61b} We assume that $(\beta_n, K_n)$
satisfies the hypotheses of Theorem {\em \ref{thm:main}} for all
$0<\alpha<\alpha_0$. Denote $m_n=m(\beta_n,K_n)$. For $\delta \in
(0,1)$ define
\[ C_n= \tilde{E}_{n,\beta_n,K_n} \{
|S_n/n^{1-\kappa} -n^{\kappa} m_n | \ \big| \ S_n/n >\delta m_n\}
\]
and
\[ D_n = \tilde{E}_{n,\beta_n,K_n} \{
|S_n/n^{1-\kappa}+ W_n/n^{1/2-\kappa} -n^{\kappa} m_n | \ \big| \
S_n/n>\delta m_n \}.
\]
Then $\lim_{n \goto \infty}|C_n - D_n| =0$.
\iffalse
\bea
\label{eqn:Thm61bSub2a}
\lim_{n \goto \infty}|C_n
- D_n| =0 . \eea
\fi
\end{lem}

{\noindent}{\bf Proof.}
\iffalse
Since $(\beta_n,K_n) \goto (\beta,
K(\beta))$ with $0 < \beta \leq \beta_c$, there exists a constant $c > 0$
such that for all $n \in \mathbb{N}$
\fi
By part (b) of Lemma \ref{lem:normalrvs} there exists a constant $c_1 > 0$ such that for all $n$
\[ \tilde{E}_{n,\beta_n,K_n}  \{  |W_n/n^{1/2-\kappa}|  \}  \le \sqrt{c_1}/n^{1/2-\kappa}.
\]
By Lemma \ref{lem:Sym1/2Small}
\[ \lim_{n \goto \infty} P_{n,\beta_n,K_n}\{ S_n/n > \delta m_n
\}={1}/{2}.
\]
It follows that there exists a constant $c_2 > 0$ such that for all sufficiently large $n$
\bea
\tilde{E}_{n,\beta_n,K_n}  \{
|W_n/n^{1/2-\kappa}| \ \big| \  S_n/n >\delta m_n \}
&=&
\frac{\tilde{E}_{n,\beta_n,K_n} \{  |W_n/n^{1/2-\kappa}| \cdot
1_{\{S_n/n >\delta m_n\} }\}
     }
     {P_{n,\beta_n,K_n}\{ S_n/n > \delta m_n\}
     } \nonumber \\
&\le&
{c_2}/{n^{1/2-\kappa}} \nonumber. \eea
Since
\bea
|S_n/n^{1-\kappa}-n^\kappa m_n|-|W_n/n^{1/2-\kappa}|
&\le&
|S_n/n^{1-\kappa}+W_n/n^{1/2-\kappa}-n^\kappa m_n| \nonumber \\
&\le&
|S_n/n^{1-\kappa}-n^\kappa m_n|+|W_n/n^{1/2-\kappa}|
,\nonumber \eea
we have for all sufficiently large $n$
\bea
\lefteqn{ D_n + {c_2}/{n^{1/2-\kappa}} } \\
&=& \tilde{E}_{n,\beta_n,K_n} \{
|S_n/n^{1-\kappa}+W_n/n^{1/2-\kappa}-n^\kappa m_n| \ \big| \ S_n/n
>\delta m_n \} + c_2/n^{1/2-\kappa}
        \nonumber \\
&=&
\frac{
\int_{\Lambda^n \times \Omega}|S_n/n^{1-\kappa}+W_n/n^{1/2-\kappa}-n^\kappa m_n|
\cdot 1_{\{S_n/n >\delta m_n\} } d (P_{n,\beta_n,K_n} \times Q)
     }
     {E_{n,\beta_n,K_n}  \{  1_{ \{S_n/n >\delta m_n \} } \}
     }
+ c_2/n^{1/2-\kappa} \nonumber \\
&\ge&
\frac{1}{E_{n,\beta_n,K_n}  \{  1_{ \{S_n/n >\delta m_n\} }\}}
\cdot \Big(\int_{\Lambda^n \times \Omega}|S_n/n^{1-\kappa}-n^\kappa m_n| \cdot
1_{ \{S_n/n >\delta m_n \} } d (P_{n,\beta_n,K_n} \times Q)
\nonumber \\
&&-
\int_{\Lambda^n \times \Omega}|W_n/n^{1/2-\kappa}| \cdot 1_{ \{ S_n/n
>\delta m_n \} } d (P_{n,\beta_n,K_n} \times Q) \Big) \nonumber \\
&&+
\tilde{E}_{n,\beta_n,K_n}  \{  |W_n/n^{1/2-\kappa}|  \
\big| \  S_n/n >\delta m_n \} \nonumber \\
&=&
E_{n,\beta_n,K_n}  \{ |S_n/n^{1-\kappa}-n^\kappa m_n| \ \big| \
S_n/n >\delta m_n \} = C_n \nonumber \eea
and \iffalse
********************************************************\fi
\bea
\lefteqn{ D_n - {c_2}/{n^{1/2-\kappa}} } \\
&=& \tilde{E}_{n,\beta_n,K_n} \{
|S_n/n^{1-\kappa}+W_n/n^{1/2-\kappa}-n^\kappa m_n| \ \big| \ S_n/n
>\delta m_n \} - c_2/n^{1/2-\kappa}
       \nonumber \\
&=&
\frac{
\int_{\Lambda^n \times \Omega}|S_n/n^{1-\kappa}+W_n/n^{1/2-\kappa}-n^\kappa m_n|
\cdot 1_{ \{S_n/n >\delta m_n\} } d(P_{n,\beta_n,K_n} \times Q)
     }
     {E_{n,\beta_n,K_n}  \{  1_{ \{S_n/n >\delta m_n\} } \}
     }
- c_2/n^{1/2-\kappa} \nonumber \\
&\le&
\frac{1}{E_{n,\beta_n,K_n}  \{  1_{ \{S_n/n >\delta m_n\} }\}}
\cdot \Big( \int_{\Lambda^n \times \Omega}|S_n/n^{1-\kappa}-n^\kappa m_n| \cdot
1_{ \{S_n/n >\delta m_n\} } d (P_{n,\beta_n,K_n} \times Q)    \nonumber \\
&&+
\int_{\Lambda^n \times \Omega}|W_n/n^{1/2-\kappa}| \cdot 1_{ \{S_n/n
>\delta m_n\} } d (P_{n,\beta_n,K_n} \times Q) \Big)      \nonumber \\
&&-
\tilde{E}_{n,\beta_n,K_n}  \{  |W_n/n^{1/2-\kappa}|  \
\big| \  S_n/n >\delta m_n \} \nonumber \\
&=&
E_{n,\beta_n,K_n}  \{ |S_n/n^{1-\kappa}-n^\kappa m_n| \ \big| \
S_n/n >\delta m_n \} = C_n \nonumber.\eea
Thus we obtain for any $\delta \in (0,1)$ and all sufficiently large $n$
\[
D_n + c_2/n^{1/2 - \kappa} \geq C_n \geq D_n - c_2/n^{1/2 - \kappa}.
\]
\iffalse
\bea
\lefteqn{\tilde{E}_{n,\beta_n,K_n} \{
|S_n/n^{1-\kappa}+W_n/n^{1/2-\kappa}-n^\kappa m_n| \ \big| \
S_n/n >\delta m_n \} + c_2/n^{1/2-\kappa}            }  \nonumber  \\
&\ge&
E_{n,\beta_n,K_n}  \{ |S_n/n^{1-\kappa}-n^\kappa m_n| \
\big| \ S_n/n >\delta m_n \}   \nonumber  \\
&\ge&
\tilde{E}_{n,\beta_n,K_n}  \{
|S_n/n^{1-\kappa}+W_n/n^{1/2-\kappa}-n^\kappa m_n| \ \big| \ S_n/n
>\delta m_n \} - c_2/n^{1/2-\kappa}   \nonumber.\eea
According
to part (c) of Theorem \ref{thm:6.1},
\fi
Since $\kappa < 1/2$ [Thm.\ \ref{thm:6.1}(c)], it follows that
$\lim_{n \goto \infty}|C_n - D_n| =0$.
This completes the proof of Lemma \ref{lem:substep2aOFThm61b}.
$\Box$

\skp Having proved Substep 2a in the proof of part (b) of Theorem
\ref{thm:6.1}, we next turn to Substep 2b.

\subsection{Proof of Substep 2b in Proof of Theorem \ref{thm:6.1} (b)}
\label{subsection:step2bOFpartB}

\iffalse Lemma \ref{lem:BIneqOFtwoEvent} relates expectation
$\tilde{E}_{n,\beta_n,K_n}$ on set $\{ S_n/n > \delta m_n \}$ to
expectation $\tilde{E}_{n,\beta_n,K_n}$ on set $A_n(\delta)$. Lemma
\ref{lem:BIneqOFtwoEvent} is used in the proof of the Substep 2b of
Theorem \ref{thm:6.1}(b) (i.e. Lemma \ref{lem:Substep2bOFThm61b}).
\fi

We recall that $W_n$ is a sequence of normal random variables with
mean 0 and variance $(2\beta_n K_n)^{-1}$ defined on a probability
space $(\Omega,\mathcal {F},Q)$. We denote by $\tilde{E}_{n,
\beta_n,K_n}$ expectation with respect to the product measure
$P_{n,\beta_n,K_n} \times Q$. Substep 2b in the proof of part (b) of
Theorem \ref{thm:6.1} states the following: \bea \lefteqn{ \lim_{n
\goto \infty} \tilde{E}_{n,\beta_n,K_n} \{
|S_n/n^{1-\kappa}+W_n/n^{1/2-\kappa}- n^{\kappa}m_n| \ \big| \
S_n/n>\delta m_n \}     } \label{eqn:substep2b} \\
&&=
\lim_{n \goto \infty}  \tilde{E}_{n,\beta_n,K_n} \{
|S_n/n^{1-\kappa}+W_n/n^{1/2-\kappa}- n^\kappa m_n| \ \big| \
S_n/n+W_n/n^{1/2}>\delta m_n \} \nonumber \\
&&= \frac{1}{\int_{\mathbb{R}}\exp[-\frac{1}{2}
g^{(2)}(\bar{x})x^2]dx } \cdot \int_{\mathbb{R}} |x| \exp[\ts
-\frac{1}{2} g^{(2)}(\bar{x})x^2]dx = \bar{z}. \nonumber \eea This
will be proved in Lemma \ref{lem:Substep2bOFThm61b}.

Part (a) of Lemma \ref{lem:Substep2bOFThm61b} relates the expectation of
$|S_n/n^{1-\kappa}+W_n/n^{1/2-\kappa}- n^{\kappa}m_n|$ conditioned
on the event $\{S_n/n > \delta m_n\}$ to the expectation of the same
random variable conditioned on the event $A_n(\bar{\delta})$ for two choices of $\bar{\delta}$. Part
(b) of the next lemma proves (\ref{eqn:substep2b}).
\iffalse states Substep 2b in the proof of part (b) of Theorem \ref{thm:6.1}. \fi
The hypotheses of this lemma coincide with
the hypotheses of Lemma \ref{lem:IneqOFtwoEvent} together with the
additional condition $\zeta > \frac{1}{2}(1-\alpha/\alpha_0)$, which
is used to prove $\tilde{\Theta}_{n,1} \goto 0$,
$\tilde{\Theta}_{n,3} \goto 0$, $\tilde{\Gamma}_{n,1} \goto 0$, and
$\tilde{\Gamma}_{n,3} \goto 0$ in part (a). According to part (c) of
Theorem \ref{thm:6.1},
$\frac{1}{2}(1-\alpha/\alpha_0)+\theta\alpha=\kappa<1/2$, which
implies $\frac{1}{2}(1-\alpha/\alpha_0) < \frac{1}{2}-\theta\alpha$.
This additional condition on $\zeta$ is consistent with the
hypothesis on $\zeta$ in Lemma \ref{lem:IneqOFtwoEvent}, which is
$0<\zeta<\frac{1}{2}-\theta\alpha$.

\begin{lem}
\label{lem:Substep2bOFThm61b} We assume that $(\beta_n, K_n)$
satisfies the hypotheses of Theorem {\em \ref{thm:main}} for all
$0<\alpha<\alpha_0$. For $\bar{\delta} \in (0,1)$ define
\[ A_n(\bar{\delta})=\{ S_n/n + W_n/n^{1/2} > \bar{\delta}m_n \}
\]
where $m_n=m(\beta_n,K_n)$. Let $\Delta \in (0,1)$ be the number
determined in part {\em (b)} of Lemma {\em \ref{lem:G(xmn)miusG(mn)}}. Assume that
$0<\alpha<\alpha_0$ and choose any numbers $\delta_1, \delta,
\delta_2$ and $\zeta$ satisfying $\Delta <
\delta_1<\delta<\delta_2<1$ and $\zeta \in
(\frac{1}{2}(1-\alpha/\alpha_0), \frac{1}{2}-\theta\alpha)$.
The following conclusions hold.
\begin{itemize}
\item [\em (a)] There exists sequences $\tilde{\Theta}_{n,1} \goto 0$, $\tilde{\Theta}_{n,2} \goto 1$,
$\tilde{\Theta}_{n,3} \goto 0$, $\tilde{\Gamma}_{n,1} \goto 0$,
$\tilde{\Gamma}_{n,2} \goto 1$, and $\tilde{\Gamma}_{n,3} \goto 0$
such that for all sufficiently large $n$
\bea
\lefteqn{
\frac{\tilde{E}_{n,\beta_n,K_n} \{
|S_n/n^{1-\kappa}+W_n/n^{1/2-\kappa}-n^\kappa m_n| \ \big|  \
A_n(\delta_1) \} + \tilde{\Theta}_{n,1}}{\tilde{\Theta}_{n,2}-
\tilde{\Theta}_{n,3}} }  \nonumber \\
&&\ge
\tilde{E}_{n,\beta_n,K_n} \{
|S_n/n^{1-\kappa}+W_n/n^{1/2-\kappa}-n^\kappa m_n| \ \big|
\ S_n/n > \delta m_n \} \nonumber \\
&&\ge
\frac{\tilde{E}_{n,\beta_n,K_n} \{
|S_n/n^{1-\kappa}+W_n/n^{1/2-\kappa}-n^\kappa m_n| \ \big|  \
A_n(\delta_2) \} - \tilde{\Gamma}_{n,1}}{\tilde{\Gamma}_{n,2} +
\tilde{\Gamma}_{n,3}} \nonumber. \eea

\iffalse where the sequences $\tilde{\Theta}_{n,1}$,
$\tilde{\Theta}_{n,3}$, $\tilde{\Gamma}_{n,1}$ and
$\tilde{\Gamma}_{n,3}$ are related to $\zeta$. \fi

\item[\em (b)] We have the conditional limit
\bea
\lefteqn{ \lim_{n \goto \infty} \tilde{E}_{n,\beta_n,K_n} \{
|S_n/n^{1-\kappa}+W_n/n^{1/2-\kappa}- n^{\kappa}m_n| \ \big| \
S_n/n>\delta m_n \}     } \nonumber \\
&&=
\lim_{n \goto \infty}  \tilde{E}_{n,\beta_n,K_n} \{
|S_n/n^{1-\kappa}+W_n/n^{1/2-\kappa}- n^\kappa m_n| \ \big| \
A_n(\delta) \}  \nonumber \\
&&=
\frac{1}{\int_{\mathbb{R}}\exp[-\frac{1}{2}
g^{(2)}(\bar{x})x^2]dx } \cdot \int_{\mathbb{R}} |x| \exp[\ts
-\frac{1}{2} g^{(2)}(\bar{x})x^2]dx = \bar{z}. \nonumber \eea

\end{itemize}

\end{lem}

\skp {\noindent}{\bf Proof of part (a) of Lemma
\ref{lem:Substep2bOFThm61b}.} The hypotheses of this lemma are a
subset of the hypotheses of Lemma \ref{lem:IneqOFtwoEvent}. We start
by proving the first inequality in part (a). By the first inequality
in part (b) of Lemma \ref{lem:IneqOFtwoEvent} and the second
inequality in part (a) of Lemma \ref{lem:IneqOFtwoEvent} we have
for all sufficiently large $n$
\bea
\label{eqn:thm6.1(b)+} \lefteqn{ \tilde{E}_{n,\beta_n,K_n} \{
|S_n/n^{1-\kappa}+W_n/n^{1/2-\kappa}-n^\kappa m_n| \ \big|  \
S_n/n>\delta m_n  \}
        } \nonumber \\
&&=
\frac{\tilde{E}_{n,\beta_n,K_n} \{
|S_n/n^{1-\kappa}+W_n/n^{1/2-\kappa}-n^\kappa m_n| \cdot 1_{ \{
S_n/n>\delta m_n \} } \} }{P_{n,\beta_n,K_n} \{
S_n/n>\delta m_n \}}  \nonumber \\
&&\le
\frac{\tilde{E}_{n,\beta_n,K_n} \{
|S_n/n^{1-\kappa}+W_n/n^{1/2-\kappa}-n^\kappa m_n| \cdot
1_{A_n(\delta_1)} \} +2n^{\kappa} e^{-cn^{2\zeta}} +
c_2n^{\kappa-1/2}e^{-c n^{2\zeta}/2} }{(P_{n,\beta_n,K_n}\times Q) \{
A_n(\delta_2)  \} - e^{-cn^{2\zeta}}}
\nonumber \\
&&=
\frac{ \frac{\tilde{E}_{n,\beta_n,K_n} \{
|S_n/n^{1-\kappa}+W_n/n^{1/2-\kappa}-n^\kappa m_n| \cdot
1_{A_n(\delta_1)} \} }{(P_{n,\beta_n,K_n}\times Q) \{ A_n(\delta_1)
\}}
+
\frac{2n^{\kappa} e^{-cn^{2\zeta}}
+
c_2n^{\kappa-1/2}e^{-c n^{2\zeta}/2} }{(P_{n,\beta_n,K_n}\times Q)
\{ A_n(\delta_1)  \}}
     }
     { \frac{(P_{n,\beta_n,K_n}\times Q) \{
A_n(\delta_2)  \}}{(P_{n,\beta_n,K_n}\times Q) \{ A_n(\delta_1)  \}}
- \frac{e^{-cn^{2\zeta}}}{(P_{n,\beta_n,K_n}\times Q) \{
A_n(\delta_1) \}}
     } \nonumber \\
&&=
\frac{\tilde{E}_{n,\beta_n,K_n} \{
|S_n/n^{1-\kappa}+W_n/n^{1/2-\kappa}-n^\kappa m_n| \  \big|   \
{A_n(\delta_1)} \} + \tilde{\Theta}_{n,1} }{ \tilde{\Theta}_{n,2} \
\
-
\ \ \tilde{\Theta}_{n,3}} \nonumber, \eea
where $c > 0$ and $c_2 > 0$ are constants and
\be
\label{eqn:BarThetan1Num}
\tilde{\Theta}_{n,1}
=
\frac{2n^{\kappa}
e^{-cn^{2\zeta}}
+
c_2n^{\kappa-1/2}e^{-cn^{2\zeta}/2}
}{(P_{n,\beta_n,K_n}\times Q) \{ A_n(\delta_1)  \}}, \ee
\be
\label{eqn:BarThetan2Num}
\tilde{\Theta}_{n,2}
=
\frac{(P_{n,\beta_n,K_n}\times Q) \{ A_n(\delta_2)
\}}{(P_{n,\beta_n,K_n}\times Q) \{ A_n(\delta_1)  \}} , \ee
and
\be
\label{eqn:BarThetan3Num}
\tilde{\Theta}_{n,3}
=
\frac{e^{-cn^{2\zeta}}}{(P_{n,\beta_n,K_n}\times Q) \{ A_n(\delta_1)
\}}. \ee

We prove the first inequality in part (a) of the present lemma by
showing that, as $n \goto \infty$, $\tilde{\Theta}_{n,1} \goto 0$, $
\tilde{\Theta}_{n,2} \goto 1 $ and $ \tilde{\Theta}_{n,3} \goto 0 $.
These limits hold for any $0<\delta_1<1$ and $0<\delta_2<1$. By
(\ref{eqn:tE1An}) in part (a) of Lemma \ref{lem:tEgammaLowerValue}
with $\bar{\delta}=\delta_1$
\bea
\label{eqn:BarThetan1}
\tilde{\Theta}_{n,1}
&=&
\frac{2n^{\kappa} e^{-cn^{2\zeta}}
+
c_2n^{\kappa-1/2}e^{-c n^{2\zeta}/2} }{(P_{n,\beta_n,K_n}\times Q) \{
A_n(\delta_1) \}}
 \\
&=&
\frac{2n^{\kappa} e^{-cn^{2\zeta}}
+
c_2n^{\kappa-1/2}e^{-c n^{2\zeta}/2} }{
\int_{-n^{\kappa}(1-\delta_1)m_n}^{\infty} \exp[-n
G_{\beta_n,K_n}(x/n^{\kappa}+m_n)]dx/Z_{n, \kappa}} \nonumber \\
&=&
\frac{\exp[nG_{\beta_n, K_n}(m_n)] \cdot Z_{n, \kappa} \cdot
(2n^{\kappa} e^{-cn^{2\zeta}}
+
c_2n^{\kappa-1/2}e^{-c n^{2\zeta}/2} )
}{\int_{-n^{\kappa}(1-\delta_1)m_n}^{\infty}\exp[-n
G_{\beta_n,K_n}(x/n^{\kappa}+m_n)+nG_{\beta_n, K_n}(m_n)]dx}
 \nonumber. \eea
We now use Lemma \ref{lem:Thm61aStep1} with $f \equiv 1$ and
$\bar{\delta} =
 \delta_1$. This gives
\bea
\lefteqn{\hspace{-1in}
\lim_{n \goto
\infty}\int_{-n^{\kappa}(1-\delta_1)m_n}^{\infty}\exp[-n
G_{\beta_n,K_n}(x/n^{\kappa}+m_n)+nG_{\beta_n, K_n}(m_n)]dx } \nonumber \\
&& = \int_{\mathbb{R}} \exp[\ts -\frac{1}{2} g^{(2)}(\bar{x}) x^2]dx
\nonumber. \eea
By Lemma \ref{lem:AsyOFZnk}, for any $\ve > 0$ and all sufficiently large $n$
\bea
\lefteqn{    \exp[nG_{\beta_n, K_n}(m_n)] \cdot Z_{n,
\kappa} \cdot (2n^{\kappa} e^{-cn^{2\zeta}} +
c_2n^{\kappa-1/2}e^{-c n^{2\zeta}/2})
} \nonumber \\
&&\le 2n^{\kappa} \exp[\varepsilon n^{1-\alpha/\alpha_0} -
c n^{2 \zeta}]
+ c_2 n^{\kappa-1/2} \exp[\varepsilon
n^{1-\alpha/\alpha_0}-c n^{2 \zeta}/2]. \nonumber
\eea
Since by hypothesis $\zeta
>
 \frac{1}{2}(1-\alpha/\alpha_0)$, it follows that
\iffalse $n^{2\zeta} - n^{1-\alpha/\alpha_0} \goto \infty$. So \fi
\[ \lim_{n \goto
\infty}\exp[nG_{\beta_n, K_n}(m_n)] \cdot Z_{n, \kappa} \cdot
(2n^{\kappa} e^{-cn^{2\zeta}}+c_2 n^{\kappa-1/2}e^{-c n^{2\zeta}/2})
=0.
\]
It follows from the last line of (\ref{eqn:BarThetan1}) that
\bea
\lim_{n \goto \infty} \tilde{\Theta}_{n,1}
&=&
\lim_{n \goto \infty}
\frac{\exp[nG_{\beta_n, K_n}(m_n)] \cdot Z_{n, \kappa} \cdot
(2n^{\kappa} e^{-cn^{2\zeta}}
+
c_2n^{\kappa-1/2}e^{-c n^{2\zeta}/2} )
}{\int_{-n^{\kappa}(1-\delta_1)m_n}^{\infty}\exp[-n
G_{\beta_n,K_n}(x/n^{\kappa}+m_n)+nG_{\beta_n, K_n}(m_n)]dx}
\nonumber \\
&=&
\frac{0}{\int_{\mathbb{R}} \exp[\ts -\frac{1}{2}
g^{(2)}(\bar{x}) x^2]dx} =0
 \nonumber, \eea
as claimed.

We now prove that $\lim_{n \goto \infty} \tilde{\Theta}_{n,2} = 0$.
By (\ref{eqn:tE1An}) in part (a) of Lemma
\ref{lem:tEgammaLowerValue} with $\bar{\delta} = \delta_1$ and $\bar{\delta} = \delta_2$
\bea
\tilde{\Theta}_{n,2}
&=&
\frac{(P_{n,\beta_n,K_n}\times Q) \{ A_n(\delta_2)
\}}{(P_{n,\beta_n,K_n}\times Q) \{ A_n(\delta_1)  \}} \nonumber  \\
&=&
\frac{\int_{-n^{\kappa}(1-\delta_2)m_n}^{\infty} \exp[-n
G_{\beta_n,K_n}(x/n^{\kappa}+m_n)]dx/Z_{n,
\kappa}}{\int_{-n^{\kappa}(1-\delta_1)m_n}^{\infty} \exp[-n
G_{\beta_n,K_n}(x/n^{\kappa}+m_n)]dx/Z_{n, \kappa}} \nonumber \\
&=& \frac{\int_{-n^{\kappa}(1-\delta_2)m_n}^{\infty}
\exp[-nG_{\beta_n,K_n}(x/n^{\kappa}+m_n)+n
G_{\beta_n,K_n}(m_n)]dx}{\int_{-n^{\kappa}(1-\delta_1)m_n}^{\infty}
\exp[-n G_{\beta_n,K_n}(x/n^{\kappa}+m_n)+ n G{\beta_n,K_n}(m_n)]dx}
\nonumber.
\eea
By Lemma \ref{lem:Thm61aStep1} for $f \equiv 1$,
$\bar{\delta} = \delta_1$, and $\bar{\delta} = \delta_2$,
 both the numerator and denominator have the same limit
$\int_{\mathbb{R}} \exp [\ts -\frac{1}{2} g^{(2)}(\bar{x})
 x^2 ] dx$. It follows that $\lim_{n \goto \infty} \tilde{\Theta}_{n,2} = 1$, as claimed.

We now prove that $\lim_{n \goto \infty} \tilde{\Theta}_{n,3} = 0$. Since $\tilde{\Theta}_{n,1} \geq \tilde{\Theta}_{n,3} > 0$,
\[
\lim_{n \goto \infty} \tilde{\Theta}_{n,1} =  0 \ \mbox{ implies } \lim_{n \goto \infty} \tilde{\Theta}_{n,3} = 0,
\]
as claimed. This completes the proof of the first inequality in part (a) of the present lemma.

We now prove the second inequality in part (a) of the present
lemma. By the second inequality in part (b) of Lemma
\ref{lem:IneqOFtwoEvent} and the first inequality in part (a) of
Lemma \ref{lem:IneqOFtwoEvent} we have for all sufficiently large $n$
\bea
\label{eqn:thm6.1(b)-} \lefteqn{
\tilde{E}_{n,\beta_n,K_n} \{ |S_n/n^{1-\kappa}+W_n/n^{1/2-\kappa}-n^\kappa
m_n| \ \big|  \ S_n/n>\delta m_n \}
        } \nonumber \\
&&=
\frac{\tilde{E}_{n,\beta_n,K_n} \{
|S_n/n^{1-\kappa}+W_n/n^{1/2-\kappa}-n^\kappa m_n| \cdot 1_{ \{
S_n/n>\delta m_n \} } \} }{P_{n,\beta_n,K_n} \{
S_n/n>\delta m_n \}}  \nonumber \\
&&\ge
\frac{\tilde{E}_{n,\beta_n,K_n} \{
|S_n/n^{1-\kappa}+W_n/n^{1/2-\kappa}-n^\kappa m_n| \cdot
1_{A_n(\delta_2)} \}
-
2n^{\kappa} e^{-cn^{2\zeta}}
-
c_2n^{\kappa-1/2}e^{-c n^{2\zeta}/2} }{(P_{n,\beta_n,K_n}\times Q) \{
A_n(\delta_1)  \}
+
e^{-cn^{2\zeta}}}
\nonumber \\
&&=
\frac{ \frac{\tilde{E}_{n,\beta_n,K_n} \{
|S_n/n^{1-\kappa}+W_n/n^{1/2-\kappa}-n^\kappa m_n| \cdot
1_{A_n(\delta_2)} \} }{(P_{n,\beta_n,K_n}\times Q) \{ A_n(\delta_2)
\}}
-
\frac{2n^{\kappa} e^{-cn^{2\zeta}} +
c_2n^{\kappa-1/2}e^{-c n^{2\zeta}/2}}{(P_{n,\beta_n,K_n}\times Q) \{
A_n(\delta_2)  \}}
     }
     { \frac{(P_{n,\beta_n,K_n}\times Q) \{
A_n(\delta_1)  \}}{(P_{n,\beta_n,K_n}\times Q) \{ A_n(\delta_2)  \}}
+ \frac{e^{-cn^{2\zeta}}}{(P_{n,\beta_n,K_n}\times Q) \{
A_n(\delta_2) \}}
     } \nonumber \\
&&=
\frac{\tilde{E}_{n,\beta_n,K_n} \{
|S_n/n^{1-\kappa}+W_n/n^{1/2-\kappa}-n^\kappa m_n| \  \big|   \
{A_n(\delta_2)} \} - \tilde{\Gamma}_{n,1} }{ \tilde{\Gamma}_{n,2} \
\ + \ \ \tilde{\Gamma}_{n,3}} \nonumber, \eea
where $c > 0$ and $c_2 > 0$ are constants and
\[ \tilde{\Gamma}_{n,1} = \frac{2n^{\kappa} e^{-cn^{2\zeta}} +c_2n^{\kappa-1/2}e^{-c n^{2\zeta}/2}}{(P_{n,\beta_n,K_n}\times Q) \{
A_n(\delta_2)  \}} ,
\]
\[ \tilde{\Gamma}_{n,2} = \frac{(P_{n,\beta_n,K_n}\times Q) \{
A_n(\delta_1)  \}}{(P_{n,\beta_n,K_n}\times Q) \{ A_n(\delta_2)  \}}
,
\]
and
\[ \tilde{\Gamma}_{n,3} = \frac{e^{-cn^{2\zeta}}}{(P_{n,\beta_n,K_n}\times Q) \{ A_n(\delta_2)
\}}  .
\]
The sequences $\tilde{\Gamma}_{n,1}$, $\tilde{\Gamma}_{n,2}$, and
$\tilde{\Gamma}_{n,3}$ are obtained from $\tilde{\Theta}_{n,1}$,
$\tilde{\Theta}_{n,2}$, and $\tilde{\Theta}_{n,3}$ in
(\ref{eqn:BarThetan1Num})--(\ref{eqn:BarThetan3Num}) by
interchanging $\delta_1$ and $\delta_2$. Hence the limits
$\tilde{\Gamma}_{n,1} \goto 0$, $\tilde{\Gamma}_{n,2} \goto 1$, and
$\tilde{\Gamma}_{n,3} \goto 0$ follow from the limits
$\tilde{\Theta}_{n,1} \goto 0$, $\tilde{\Theta}_{n,2} \goto 1$, and
$\tilde{\Theta}_{n,3} \goto 0$, which hold for any $0<\delta_1<1$
and $0<\delta_2<1$.
The proof of part (a) of Lemma \ref{lem:Substep2bOFThm61b} is
complete.

\skp {\noindent}{\bf Proof of part (b) of Lemma
\ref{lem:Substep2bOFThm61b}.} We know from part (b) of Lemma
\ref{lem:Step1OFThm61b} that, as $n \goto \infty$, for any $\bar{\delta}=\delta_1$ and
$\bar{\delta}=\delta_2$, \ $\tilde{E}_{n,\beta_n,K_n} \{
|S_n/n^{1-\kappa}+W_n/n^{1/2-\kappa}-n^\kappa m_n| \ \big|  \
A_n(\bar{\delta}) \}$ has the same limit
\[ \frac{1}{\int_{\mathbb{R}}
\exp[-\frac{1}{2}g^{(2)}(\bar{x})x^2]dx} \nonumber \cdot
\int_{\mathbb{R}} |x| \exp[\ts -\frac{1}{2}g^{(2)}(\bar{x})x^2]dx =
\bar{z}.
\]
By sending $n \goto \infty$ in the inequality in part (a), we have
\bea
\bar{z}
&=&
\lim_{n \goto \infty}  \tilde{E}_{n,\beta_n,K_n}
\{ |S_n/n^{1-\kappa}+W_n/n^{1/2-\kappa}- n^\kappa m_n| \ \big| \
A_n(\delta_1)
\}   \nonumber \\
&\ge&
\limsup_{n \goto \infty} \tilde{E}_{n,\beta_n,K_n} \{
|S_n/n^{1-\kappa}+W_n/n^{1/2-\kappa}- n^{\kappa}m_n| \ \big| \
S_n/n>\delta m_n \}     \nonumber \\
&\ge&
\liminf_{n \goto \infty} \tilde{E}_{n,\beta_n,K_n} \{
|S_n/n^{1-\kappa}+W_n/n^{1/2-\kappa}- n^{\kappa}m_n| \ \big| \
S_n/n>\delta m_n \}     \nonumber \\
&=&
\lim_{n \goto \infty}  \tilde{E}_{n,\beta_n,K_n} \{
|S_n/n^{1-\kappa}+W_n/n^{1/2-\kappa}- n^\kappa m_n| \ \big| \
A_n(\delta_2)
\}  \nonumber \\
&=&
\bar{z}. \nonumber \eea
Because the first and last terms in this
display are the same, it follows that
\bea
\lefteqn{ \lim_{n \goto
\infty} \tilde{E}_{n,\beta_n,K_n} \{
|S_n/n^{1-\kappa}+W_n/n^{1/2-\kappa}- n^{\kappa}m_n| \ \big| \
S_n/n>\delta m_n \}     } \nonumber \\
&&=
\frac{1}{\int_{\mathbb{R}}
\exp[-\frac{1}{2}g^{(2)}(\bar{x})x^2]dx} \nonumber \cdot
\int_{\mathbb{R}} |x| \exp[\ts -\frac{1}{2}g^{(2)}(\bar{x})x^2]dx =
\bar{z}. \nonumber \eea
On the other hand, by part (b) of Lemma
\ref{lem:Step1OFThm61b} with $\bar{\delta}=\delta$
\bea
\lefteqn{
\lim_{n \goto \infty} \tilde{E}_{n,\beta_n,K_n} \{
|S_n/n^{1-\kappa}+W_n/n^{1/2-\kappa}- n^\kappa m_n| \ \big| \
A_n(\delta)
\}    } \nonumber \\
&&=
\frac{1}{\int_{\mathbb{R}}
\exp[-\frac{1}{2}g^{(2)}(\bar{x})x^2]dx} \nonumber \cdot
\int_{\mathbb{R}} |x| \exp[\ts -\frac{1}{2}g^{(2)}(\bar{x})x^2]dx =
\bar{z}. \nonumber \eea
The proof of part (b) of Lemma
\ref{lem:Substep2bOFThm61b} is complete.
 $\Box$

\skp We now put together the pieces to complete the proof of part (b) of Theorem \ref{thm:6.1}.
Let $\delta$ be any number satisfying $\Delta < \delta <1$,
where $\Delta \in (0,1)$ is determined in part (b) of Lemma
\ref{lem:G(xmn)miusG(mn)}. \iffalse Just before subsection
\ref{subsection:step1OFpartB} we sketch the proof of part (b) of
Theorem \ref{thm:6.1}. \fi
The proof of part (b) of Theorem \ref{thm:6.1} is divided into
Step 1, Substep 2a, and Substep 2b. Step 1 is done in part (b) of
Lemma \ref{lem:Step1OFThm61b}. There we prove that with $\bar{\delta} = \delta$
\bea
\lefteqn{ \hspace{-.25in} \lim_{n \goto \infty} \tilde{E}_{n,\beta_n,K_n} \{
|S_n/n^{1-\kappa}+W_n/n^{1/2-\kappa}- n^{\kappa}m_n| \ \big| \
S_n/n+W_n/n^{1/2}>\delta m_n
\} } \nonumber \\
&&=
\frac{1}{\int_{\mathbb{R}}
\exp[-\frac{1}{2}g^{(2)}(\bar{x})x^2]dx} \nonumber \cdot
\int_{\mathbb{R}} |x| \exp[\ts -\frac{1}{2}g^{(2)}(\bar{x})x^2]dx =
\bar{z}. \nonumber \eea %%%*********
Substep 2a is done in Lemma \ref{lem:substep2aOFThm61b}. There we
prove that $\lim_{n \goto \infty}|C_n - D_n| =0$, where
\[ C_n= E_{n,\beta_n,K_n} \{
|S_n/n^{1-\kappa} -n^{\kappa} m_n | \ \big| \ S_n/n >\delta m_n\}
\]
and
\[ D_n = \tilde{E}_{n,\beta_n,K_n} \{
|S_n/n^{1-\kappa}+ W_n/n^{1/2-\kappa} -n^{\kappa} m_n | \ \big| \
S_n/n>\delta m_n \}.
\]
Substep 2b is done in part (b) of Lemma
\ref{lem:Substep2bOFThm61b}. There we prove that
\bea \lefteqn{ \lim_{n \goto
\infty} D_n } \nonumber \\
&&= \lim_{n \goto \infty} \tilde{E}_{n,\beta_n,K_n} \{
|S_n/n^{1-\kappa}+W_n/n^{1/2-\kappa}- n^{\kappa}m_n| \ \big| \
S_n/n+W_n/n^{1/2}>\delta m_n \} \nonumber
\\
&&=
\frac{1}{\int_{\mathbb{R}}\exp[-\frac{1}{2}
g^{(2)}(\bar{x})x^2]dx } \cdot \int_{\mathbb{R}} |x| \exp[ \ts
-\frac{1}{2} g^{(2)}(\bar{x})x^2]dx = \bar{z}. \nonumber \eea
Combining these limits yields
\bea
\lefteqn{ \lim_{n \goto
\infty} E_{n,\beta_n,K_n} \{ |S_n/n^{1-\kappa}- n^{\kappa}m_n| \
\big| \ S_n/n>\delta m_n
\} } \nonumber \\
&&=
\lim_{n \goto \infty}  \tilde{E}_{n,\beta_n,K_n} \{
|S_n/n^{1-\kappa}+W_n/n^{1/2-\kappa}- n^{\kappa}m_n| \ \big| \
S_n/n+W_n/n^{1/2}>\delta m_n
\} \nonumber \\
&&=
\frac{1}{\int_{\mathbb{R}}\exp[-\frac{1}{2}
g^{(2)}(\bar{x})x^2]dx } \cdot \int_{\mathbb{R}} |x| \exp[\ts
-\frac{1}{2} g^{(2)}(\bar{x})x^2]dx = \bar{z}.\nonumber \eea
This gives the conditional limit stated in part (b) of Theorem
\ref{thm:6.1}:
\bea
\lefteqn{ \lim_{n \goto \infty} n^\kappa
E_{n,\beta_n,K_n} \{ |S_n/n-m(\beta_n,K_n)|  \  \big| \
S_n/n>\delta m(\beta_n,K_n) \} }
\nonumber \\
&&=
\lim_{n \goto \infty}  E_{n,\beta_n,K_n} \{ |S_n/n^{1-\kappa} -
n^\kappa m(\beta_n,K_n)|  \  \big| \  S_n/n>\delta m(\beta_n,K_n) \}
=\bar{z} \nonumber. \eea
The proof of part (b) of Theorem \ref{thm:6.1} is complete. $\Box$

\iffalse \skp Combined part (b) of Lemma \ref{lem:Step1OFThm61b}
which gives Step 1, Lemma \ref{lem:substep2aOFThm61b} which gives
Substep 2a and part (b) of Lemma \ref{lem:Substep2bOFThm61b} which
gives Substep 2b , we have completed the proof of part (b) of
Theorem \ref{thm:6.1}. \fi

%===================== end section =====================================

%===================== begin appendix =====================================
%\appendix

%\section{Appendix A}
%\label{appendix:A}

\vspace{.5in} \noi \LARGE {\bf Appendix} \vspace{-.2in} \normalsize
\appendix

\renewcommand{\thesection}{\Alph{section}}
\renewcommand{\theequation}
{\Alph{section}.\arabic{equation}}
\renewcommand{\thedefn}
{\Alph{section}.\arabic{defn}}
\renewcommand{\theass}
{\Alph{section}.\arabic{ass}}

\section{Proof That Sequences 1a--5a Satisfy the Limits in Hypothesis (iii$^\prime$) of Theorem \ref{thm:main}}
\label{section:appendixA}
\setcounter{equation}{0}

In this appendix we prove that sequences 1a--5a satisfy the limits
in hypothesis (iii$^\prime$) of Theorem \ref{thm:main}. These limits
take the following form.

\begin{itemize}
\item[(a)] Assume that $g$ has degree 4. For $\forall \alpha
\in (0, \alpha_0)$ and for $j=2,3,4$
\[ \lim_{n \goto \infty} n^{\alpha/\alpha_0 - j \theta\alpha}
  G_{\beta_n,K_n}^{(j)}(m(\beta_n,K_n)) = g^{(j)}(\bar{x}).
\]

\item[ (b)] Assume that $g$ has degree $6$. For $\forall \alpha \in (0,
\alpha_0)$ and for $j=2,3,4,5,6$
\[ \lim_{n \goto \infty} n^{\alpha/\alpha_0 - j \theta\alpha}
  G_{\beta_n,K_n}^{(j)}(m(\beta_n,K_n)) = g^{(j)}(\bar{x}).
\]
\end{itemize}
We do this by verifying the limits (\ref{eqn: 1assumeLem51}) and
(\ref{eqn: 2assumeLem51}) in Lemma \ref{lem:5.1}. Let
$\varepsilon_n$ denote a sequence that converges to 0 and that
represents the various error terms arising in the proof. We use the
same notation $\varepsilon_n$ to represent different error terms.
\iffalse Throughout this section, let $u=1-\alpha/\alpha_0$,
$\gamma=\theta\alpha$, and $m_n=m(\beta_n, K_n)$. \fi

\begin{lem}
\label{lem:5.1} We assume the hypotheses of Theorem {\em
\ref{thm:3.1}}. We also assume {\em (\ref{eqn: 1assumeLem51})} when
the degree of the Ginzburg-Landau polynomial $g$ is $4$ and {\em
(\ref{eqn: 2assumeLem51})} when the degree of $g$ is $6$.

\begin{itemize}
\item[{\em(a)}]Assume that $g$ has degree $4$ and that for $\alpha \in (0, \alpha_0)$
and for $j= 2,3,4$ \be \label{eqn: 1assumeLem51} \lim_{n \goto
\infty} n^{\alpha/\alpha_0 - j \theta\alpha}
  G_{\beta_n,K_n}^{(j)}(x/n^{\theta\alpha}) = g^{(j)}(x)
\ee uniformly for $x$ in compact subsets of $\mathbb{R}$. Then we
have \be \label{eqn: 1concluLem51}
 \lim_{n \goto \infty} n^{\alpha/\alpha_0 - j \theta\alpha}
  G_{\beta_n,K_n}^{(j)}(m(\beta_n,K_n)) = g^{(j)}(\bar{x}).
\ee
\item[{\em (b)}]
Assume that $g$ has degree $6$ and that for $ \alpha \in (0,
\alpha_0)$ and for $j=2,3,4,5,6$ \be \label{eqn: 2assumeLem51}
\lim_{n \goto \infty} n^{\alpha/\alpha_0 - j \theta\alpha}
  G_{\beta_n,K_n}^{(j)}(x/n^{\theta\alpha}) = g^{(j)}(x)
\ee uniformly for $x$ in compact subsets of $\mathbb{R}$. Then we
have \be \label{eqn: 2concluLem51} \lim_{n \goto \infty}
n^{\alpha/\alpha_0 - j \theta\alpha}
  G_{\beta_n,K_n}^{(j)}(m(\beta_n,K_n)) = g^{(j)}(\bar{x}).
\ee
\end{itemize}
\end{lem}

\noindent{\bf Proof.} We write $m_n=m(\beta_n, K_n)$. When $g$ has
degree $4$, we have for $j= 2,3,4$, and when $g$ has degree $6$, we
have for $j=2,3,4,5,6$ \bea \label{eqn: proofLem51}
|n^{\alpha/\alpha_0-j\theta\alpha} G_{\beta_n K_n}^{(j)}(m_n) -
g^{(j)}(\bar{x}) | &\le& |n^{\alpha/\alpha_0-j\theta\alpha}
G_{\beta_n K_n}^{(j)}(n^{\theta\alpha}m_n/n^{\theta\alpha}) -
g^{(j)}(n^{\theta\alpha}m_n)|   \nonumber
\\ &&+
|g^{(j)}(n^{\theta\alpha}m_n)-g^{(j)}(\bar{x})|. \eea Let $\Xi$ be
any compact subset of $\mathbb{R}$. By hypothesis (\ref{eqn:
1assumeLem51}) for $j= 2, 3, 4$ when $g$ has degree 4 and by
hypothesis (\ref{eqn: 2assumeLem51}) for $j= 2, 3, 4, 5, 6$ when $g$
has degree 6
\[ \lim_{n \goto
\infty} \sup_{x \in \Xi}  | n^{\alpha/\alpha_0 - j \theta\alpha}
  G_{\beta_n,K_n}^{(j)}(x/n^{\theta\alpha})    -    g^{(j)}(x)  |
  =0.
\]
According to Theorem \ref{thm:3.1}, $n^{\theta\alpha}m_n \goto
\bar{x}$, and so for any $\varepsilon>0$ the sequence
$n^{\theta\alpha}m_n$ lies in the compact set $[\bar{x}-\varepsilon,
\bar{x}+\varepsilon ]$ for all sufficiently large $n$. It follows
that the first term on the right-hand side of (\ref{eqn:
proofLem51}) converges to 0 as $n \goto \infty$. Because of the
limit $n^{\theta\alpha}m_n \goto \bar{x}$ and the continuity of
$g^{(j)}$, the second term on the right-hand side of (\ref{eqn:
proofLem51}) also converges to 0 as $n \goto \infty$. We conclude
that for $j= 2, 3, 4$ when $g$ has degree 4 and for $j= 2, 3, 4, 5,
6$ when $g$ has degree 6
\[ | n^{\alpha/\alpha_0-j\theta\alpha}
G_{\beta_n, K_n}^{(j)}(m_n) - g^{(j)}(\bar{x}) | \goto 0    \  \
\mbox{as} \ \ n \goto \infty.
\]
This completes the proof of the lemma. $\Box$

\iffalse ******************************************************* \fi
\iffalse These limits take the following two forms (\ref{eqn:
Appendix_1assumeLem51}) and (\ref{eqn: Appendix_2assumeLem51})
depending on the degree of Ginzburg-Landau polynomial $g$.

\begin{itemize}
\item[(a)]When $g$ has degree $4$, the limit (\ref{eqn: 1assumeLem51}) holds; i.e., for $j= 2,3,4$
\be
\label{eqn: Appendix_1assumeLem51}
\lim_{n \goto \infty}
n^{\alpha/\alpha_0 - j \theta\alpha}
G_{\beta_n,K_n}^{(j)}(x/n^{\theta\alpha})
=
g^{(j)}(x) \ee
uniformly
for $x$ in compact subsets of $\mathbb{R}$. By part (a) of Lemma
\ref{lem:5.1} we then have the first form of the limits in
hypothesis (iii$^\prime$) of Theorem \ref{thm:main}; i.e., for
$j=2,3,4$
\[
 \lim_{n \goto \infty} n^{\alpha/\alpha_0 - j \theta\alpha}
  G_{\beta_n,K_n}^{(j)}(m(\beta_n,K_n)) = g^{(j)}(\bar{x}).
\]
\item[ (b)]
When $g$ has degree $6$, The limit (\ref{eqn: 2assumeLem51}) holds;
i.e., for $j=2,3,4,5,6$
\be
\label{eqn: Appendix_2assumeLem51}
\lim_{n \goto \infty} n^{\alpha/\alpha_0 - j \theta\alpha}
  G_{\beta_n,K_n}^{(j)}(x/n^{\theta\alpha}) = g^{(j)}(x)
\ee
uniformly for $x$ in compact subsets of $\mathbb{R}$. By part
(b) of Lemma \ref{lem:5.1} we then have the second form of the
limits in hypothesis (iii$^\prime$) of Theorem \ref{thm:main}; i.e.,
for $j=2,3,4, 5, 6$
\[\lim_{n \goto \infty} n^{\alpha/\alpha_0 - j
\theta\alpha}
  G_{\beta_n,K_n}^{(j)}(m(\beta_n,K_n)) = g^{(j)}(\bar{x}).
\]
\end{itemize}
\fi

\iffalse ******************************************************* \fi

\skp

The main point of this section is to justify rigorously the limits
in (\ref{eqn: 1assumeLem51}) and (\ref{eqn: 2assumeLem51}) for
sequences 1a--5a. We start by doing some preparatory work involving
the Taylor expansion of $G_{\beta_n, K_n}^{(j)}(x/{n^{\gamma}})$ for
$\gamma
>0$.

\iffalse Formulas (\ref{eqn: 4j2})-- (\ref{eqn: 6j6}) will be used
in verifying hypothesis (iii$^\prime$) of Theorem \ref{thm:main}.
\fi

\skp
%%%------------ pre deg(g)=4  begin---------------%%
\noindent {\it Case 1: g has degree {\em 4}, j = {\em 2, 3, 4}}.
This case arises for sequences 1 and 2, which 
converge to a second-order point $(\beta, K(\beta))$ for $0<\beta<\beta_c$.
We consider the Taylor expansions of $G_{\beta_n,
K_n}^{(j)}(x/n^\gamma)$ to order 4 with error terms. Since
$K(\beta)= (4\beta+2)/4\beta$ is continuous and $(\beta_n, K_n)$
converges to $(\beta, K(\beta))$, we have $\beta_n K_n / K(\beta_n)
\goto \beta$. Thus the coefficients in Taylor expansion of
$G^{(j)}_{\beta_n, K_n}(x/ n^{\gamma})$ are given by
\[G_{\beta_n, K_n}^{(2)}(0) = 2 \beta_n K_n - \frac{8 \beta_n^2
K_n^2}{e^{\beta_n}+2}= \frac{2\beta_n K_n (K(\beta_n)-
K_n)}{K(\beta_n)}=2 \beta (K(\beta_n)-K_n)(1+\varepsilon_n),
\]
\[ G_{\beta_n, K_n}^{(3)}(0) =  0,
\]
\[ G_{\beta_n, K_n}^{(4)}(0) = \frac{2(2\beta_n
K_n)^4(4-e^{\beta_n})}{(e^{\beta_n}+2)^2}.
\]
Let $c_4(\beta)=(e^{\beta}+2)^2(4-e^{\beta})/(8 \cdot 4! )$. Since
$2 \beta_n K_n \goto 2 \beta K(\beta)=(e^{\beta}+2)/2$, we have
\[ G_{\beta_n, K_n}^{(4)}(0) = (e^{\beta}+2)^2(4-e^{\beta})(1+\varepsilon_n)/8 = c_4(\beta) (1+\varepsilon_n) \cdot 4!.
\]
Thus for all $n \in \mathbb{N}$, any $\gamma >0$, any $R>0$, and all
$x \in \mathbb{R}$ satisfying $|x/n^\gamma|<R$, we have the Taylor
expansion
%%%%%%%%%% j=2 %%%%%%%%%%%%%%%%
\[G_{\beta_n, K_n}^{(2)}(x/n^\gamma) = G_{\beta_n, K_n}^{(2)}(0) +
\frac{G_{\beta_n, K_n}^{(4)}(0)}{2!} \cdot \frac{x^2}{n^{2\gamma}} +
\mbox{O} \! \left( \frac{1}{n^{3\gamma}} \right)x^3.
\]
Multiplying both sides by $n^{1-u-2\gamma}$ for $u>0$ yields
\bea
\label{eqn:4j2}
n^{1-u-2\gamma}G^{(2)}_{\beta_n, K_n}(x/n^\gamma)
&=&
\frac{1}{n^{2\gamma-1+u}} G_{\beta_n, K_n}^{(2)}(0)
+
\frac{1}{n^{4\gamma-1+u}} \cdot \frac{G_{\beta_n,
K_n}^{(4)}(0)}{2!}x^2  \\
&&+
\ \mbox{O} \! \left(\frac{1}{n^{5\gamma-1+u}}\right)x^3 \nonumber \\
&=&
\frac{1}{n^{2\gamma-1+u}} \cdot 2\beta(K(\beta_n)-K_n)(1+\varepsilon_n) \nonumber \\
&&+
\frac{1}{n^{4\gamma-1+u}} \cdot
\frac{c_4(\beta)(1+\varepsilon_n)\cdot 4!}{2!}x^2
+
\ \mbox{O} \!
\left(\frac{1}{n^{5\gamma-1+u}}\right)x^3. \nonumber \eea

%%%%%%%%% j=3 %%%%%%%%%%%%%%%
For all $n \in \mathbb{N}$, any $\gamma >0$, any $R>0$, and all $x
\in \mathbb{R}$ satisfying $|x/n^\gamma|<R$, we have the Taylor
expansion
\[G_{\beta_n, K_n}^{(3)}(x/n^\gamma) =
G_{\beta_n, K_n}^{(4)}(0) \cdot \frac{x}{n^{\gamma}} + \mbox{O} \!
\left(\frac{1}{n^{2\gamma}}\right)x^2.
\]
Multiplying both sides by $n^{1-u-3\gamma}$ for $u>0$ yields
\bea
\label{eqn:4j3}
n^{1-u-3\gamma}G^{(3)}_{\beta_n,
K_n}(x/n^\gamma)
&=&
\frac{1}{n^{4\gamma-1+u}} \cdot G_{\beta_n,
K_n}^{(4)}(0) \cdot x + \ \mbox{O} \! \left(\frac{1}{n^{5\gamma-1+u}}\right)x^2 \\
&=&
\frac{1}{n^{4\gamma-1+u}} \cdot c_4(\beta)(1+\varepsilon_n)\cdot
4! \cdot x + \ \mbox{O} \! \left(\frac{1}{n^{5\gamma-1+u}}\right)x^2. \nonumber 
\eea

%%%%%%%%% j=4 %%%%%%%%%%%%%%%
For all $n \in \mathbb{N}$, any $\gamma >0$, any $R>0$, and all $x
\in \mathbb{R}$ satisfying $|x/n^\gamma|<R$, we have the Taylor
expansion
\[G_{\beta_n, K_n}^{(4)}(x/n^\gamma) =
G_{\beta_n, K_n}^{(4)}(0) + \mbox{O} \!
\left(\frac{1}{n^{\gamma}}\right)x.
\]
Multiplying both sides by $n^{1-u-4\gamma}$ for $u>0$ yields
\bea
\label{eqn:4j4}
n^{1-u-4\gamma}G^{(4)}_{\beta_n,
K_n}(x/n^\gamma)
&=&
\frac{1}{n^{4\gamma-1+u}} \cdot G_{\beta_n,
K_n}^{(4)}(0) + \ \mbox{O} \! \left(\frac{1}{n^{5\gamma-1+u}}\right)x \\
&=&
\frac{1}{n^{4\gamma-1+u}} \cdot c_4(\beta)(1+\varepsilon_n)\cdot
4!
+
\ \mbox{O} \! \left(\frac{1}{n^{5\gamma-1+u}}\right)x. \nonumber \eea

In formulas (\ref{eqn:4j2})--(\ref{eqn:4j4}) the big-oh terms are
uniform for $x \in (-Rn^\gamma, Rn^\gamma)$. We will use
(\ref{eqn:4j2})--(\ref{eqn:4j4}) to verify hypothesis (\ref{eqn:
1assumeLem51}) for sequences 1a and 2a.

%%%------------ pre deg(g)=4  done---------------%%
\skp
%%%------------ pre deg(g)=6  begin---------------%%
\noindent {\it Case 2: g has degree {\em 6}, j = {\em 2, 3, 4, 5,
6}}. This case arises for sequences 3, 4 and 5, which converge to
the tricritical point $(\beta_c, K(\beta_c))$. \iffalse
Specifically, although the Ginzburg-Landau polynomial $g$ of
sequence 6 has degree 4, the following preparatory work still can be
used to verify assumption (\ref{eqn: 1assumeLem51}) for sequence 6.
\fi  We consider the Taylor expansions of $G_{\beta_n,
K_n}^{(j)}(x/n^\gamma)$ to order 6 with error terms. Since
$K(\beta)= (4\beta+2)/4\beta$ is continuous and $(\beta_n, K_n)$
converges to $(\beta_c, K(\beta_c))$, we have $\beta_n K_n /
K(\beta_n) \goto \beta_c$. Thus the coefficient in Taylor expansion
of $G^{(j)}_{\beta_n, K_n}(x/ n^{\gamma})$ are given by
\[G_{\beta_n, K_n}^{(2)}(0) = 2 \beta_n K_n - \frac{8 \beta_n^2
K_n^2}{e^{\beta_n}+2}= \frac{2\beta_n K_n (K(\beta_n)-
K_n)}{K(\beta_n)}=2 \beta_c (K(\beta_n)-K_n)(1+\varepsilon_n),
\]
\[ G_{\beta_n, K_n}^{(3)}(0) =  0,
\]
\[ G_{\beta_n, K_n}^{(4)}(0) = \frac{2(2\beta_n
K_n)^4(4-e^{\beta_n})}{(e^{\beta_n}+2)^2}.
\]
Let $c_4 =3/16$. Since $2 \beta_n K_n \goto 2 \beta_c
K(\beta_c)=(e^{\beta_c}+2)/2=3$ and $e^{\beta_n}+2 \goto
e^{\beta_c}+2 = 6$, we have
\[ G_{\beta_n, K_n}^{(4)}(0) = 2 \cdot 3^4 (4-e^{\beta_n})(1+\varepsilon_n)/6^2 = c_4 (4-e^{\beta_n}) (1+\varepsilon_n) \cdot
4!,
\]
\[ G_{\beta_n, K_n}^{(5)}(0) =  0.
\]
Let $c_6=9/40$, since $ G_{\beta_n, K_n}^{(6)}(0) \goto G_{\beta_c,
K(\beta_c)}^{(6)}(0) = 2 \cdot 3^4$, we have
\[ G_{\beta_n, K_n}^{(6)}(0)= 2 \cdot 3^4 (1+\varepsilon_n) = c_6(1+\varepsilon_n)
\cdot 6!.
\]
Thus for all $n \in \mathbb{N}$, any $\gamma >0$, any $R>0$, and all
$x \in \mathbb{R}$ satisfying $|x/n^\gamma|<R$, we have the Taylor
expansion
%%%%%%%%%% j=2 %%%%%%%%%%%%%%%%
\[G_{\beta_n, K_n}^{(2)}(x/n^\gamma) = G_{\beta_n, K_n}^{(2)}(0) +
\frac{G_{\beta_n, K_n}^{(4)}(0)}{2!} \cdot \frac{x^2}{n^{2\gamma}} +
\frac{G_{\beta_n, K_n}^{(6)}(0)}{4!} \cdot \frac{x^4}{n^{4\gamma}} +
\mbox{O} \! \left(\frac{1}{n^{5\gamma}}\right)x^5.
\]
Multiplying both sides by $n^{1-u-2\gamma}$ for $u>0$ yields
\bea
\label{eqn:6j2}
n^{1-u-2\gamma}G^{(2)}_{\beta_n, K_n}(x/n^\gamma)
&=&
\frac{1}{n^{2\gamma-1+u}} G_{\beta_n, K_n}^{(2)}(0)
+
\frac{1}{n^{4\gamma-1+u}} \cdot \frac{G_{\beta_n,
K_n}^{(4)}(0)}{2!}x^2 \\
&&+
\frac{1}{n^{6\gamma-1+u}} \cdot \frac{G_{\beta_n,
K_n}^{(6)}(0)}{4!}x^4
+
\ \mbox{O} \! \left(\frac{1}{n^{7\gamma-1+u}}\right)x^5 \nonumber \\
&=&
\frac{1}{n^{2\gamma-1+u}} \cdot 2\beta_c(K(\beta_n)-K_n)(1+\varepsilon_n) \nonumber \\
&&+
\frac{1}{n^{4\gamma-1+u}} \cdot
\frac{c_4(4-e^{\beta_n})(1+\varepsilon_n)\cdot 4!}{2!}x^2
+
\frac{c_6(1+\varepsilon_n)\cdot 6!}{4!}x^4 \nonumber \\
&&+
\
\mbox{O} \! \left(\frac{1}{n^{7\gamma-1+u}}\right)x^5. \nonumber \eea

%%%%%%%%% j=3 %%%%%%%%%%%%%%%
For all $n \in \mathbb{N}$, any $\gamma >0$, any $R>0$, and all $x
\in \mathbb{R}$ satisfying $|x/n^\gamma|<R$, we have the Taylor
expansion
\[G_{\beta_n, K_n}^{(3)}(x/n^\gamma) =
G_{\beta_n, K_n}^{(4)}(0) \cdot \frac{x}{n^{\gamma}} +
\frac{G_{\beta_n, K_n}^{(6)}(0)}{3!}\cdot \frac{x^3}{n^{3\gamma}} +
\mbox{O} \! \left(\frac{1}{n^{4\gamma}}\right)x^4.
\]
Multiplying both sides by $n^{1-u-3\gamma}$ for $u>0$ yields
\bea
\label{eqn:6j3}
n^{1-u-3\gamma}G^{(3)}_{\beta_n, K_n}(x/n^\gamma)
&=&
\frac{1}{n^{4\gamma-1+u}} \cdot G_{\beta_n, K_n}^{(4)}(0) \cdot
x
+
\frac{1}{n^{6\gamma-1+u}} \cdot \frac{G_{\beta_n,
K_n}^{(6)}(0)}{3!}\cdot x^3 \\
&&+
\ \mbox{O} \! \left(\frac{1}{n^{7\gamma-1+u}}\right)x^4 \nonumber \\
&=&
\frac{1}{n^{4\gamma-1+u}} \cdot
c_4(4-e^{\beta_n})(1+\varepsilon_n)\cdot 4! \cdot x
+\frac{1}{n^{6\gamma-1+u}} \cdot \frac{c_6(1+\varepsilon_n) \cdot
6!}{3!}x^3 \nonumber \\
&&+
\ \mbox{O} \! \left(\frac{1}{n^{7\gamma-1+u}}\right)x^4. \nonumber \eea

%%%%%%%%% j=4 %%%%%%%%%%%%%%%
For all $n \in \mathbb{N}$, any $\gamma >0$, any $R>0$, and all $x
\in \mathbb{R}$ satisfying $|x/n^\gamma|<R$, we have the Taylor
expansion
\[G_{\beta_n, K_n}^{(4)}(x/n^\gamma) =
G_{\beta_n, K_n}^{(4)}(0) + \frac{G_{\beta_n, K_n}^{(6)}(0)}{2!}
\cdot \frac{x^2}{n^{2\gamma}} + \mbox{O} \!
\left(\frac{1}{n^{3\gamma}}\right)x^3.
\]
Multiplying both sides by $n^{1-u-4\gamma}$ for $u>0$ yields
\bea
\label{eqn:6j4}
n^{1-u-4\gamma}G^{(4)}_{\beta_n, K_n}(x/n^\gamma)
&=&
\frac{1}{n^{4\gamma-1+u}} \cdot G_{\beta_n, K_n}^{(4)}(0)
+
\frac{1}{n^{6\gamma-1+u}} \cdot \frac{G_{\beta_n,
K_n}^{(6)}(0)}{2!}\cdot x^2
\\
&&+
\ \mbox{O} \! \left(\frac{1}{n^{7\gamma-1+u}}\right)x^3 \nonumber \\
&=&
\frac{1}{n^{4\gamma-1+u}} \cdot
c_4(4-e^{\beta_n})(1+\varepsilon_n)\cdot 4!
+
\frac{1}{n^{6\gamma-1+u}} \cdot \frac{c_6(1+\varepsilon_n) \cdot
6!}{2!}\cdot x^2 \nonumber
\\
&&+
\ \mbox{O} \! \left(\frac{1}{n^{7\gamma-1+u}}\right)x^3. \nonumber \eea

%%%%%%%%% j=5 %%%%%%%%%%%%%%%
For all $n \in \mathbb{N}$, any $\gamma >0$, any $R>0$, and all $x
\in \mathbb{R}$ satisfying $|x/n^\gamma|<R$, we have the Taylor
expansion
\[G_{\beta_n, K_n}^{(5)}(x/n^\gamma) =
G_{\beta_n, K_n}^{(6)}(0) \cdot \frac{x}{n^\gamma} + \mbox{O} \!
\left(\frac{1}{n^{2\gamma}}\right)x^2.
\]
Multiplying both sides by $n^{1-u-5\gamma}$ for $u>0$ yields
\bea
\label{eqn:6j5}
n^{1-u-5\gamma}G^{(5)}_{\beta_n, K_n}(x/n^\gamma)
&=&
\frac{1}{n^{6\gamma-1+u}} \cdot G_{\beta_n,
K_n}^{(6)}(0) x
+
\ \mbox{O} \! \left(\frac{1}{n^{7\gamma-1+u}}\right)x^2  \\
&=&
\frac{1}{n^{6\gamma-1+u}} \cdot c_6(1+\varepsilon_n) \cdot
6!\cdot x + \ \mbox{O} \! \left(\frac{1}{n^{7\gamma-1+u}}\right)x^2. \nonumber
\eea

%%%%%%%%% j=6 %%%%%%%%%%%%%%%
For all $n \in \mathbb{N}$, any $\gamma >0$, any $R>0$, and all $x
\in \mathbb{R}$ satisfying $|x/n^\gamma|<R$, we have the Taylor
expansion
\[G_{\beta_n, K_n}^{(6)}(x/n^\gamma) =
G_{\beta_n, K_n}^{(6)}(0) + \mbox{O} \!
\left(\frac{1}{n^{\gamma}}\right)x.
\]
Multiplying both sides by $n^{1-u-6\gamma}$ for $u>0$ yields
\bea
\label{eqn:6j6}
n^{1-u-6\gamma}G^{(6)}_{\beta_n, K_n}(x/n^\gamma)
&=&
\frac{1}{n^{6\gamma-1+u}} \cdot G_{\beta_n,
K_n}^{(6)}(0)
+
\ \mbox{O} \! \left(\frac{1}{n^{7\gamma-1+u}}\right)x \\
&=&
\frac{1}{n^{6\gamma-1+u}} \cdot c_6(1+\varepsilon_n) \cdot
6!\cdot + \ \mbox{O} \! \left(\frac{1}{n^{7\gamma-1+u}}\right)x. \nonumber \eea

In formulas (\ref{eqn:6j2})--(\ref{eqn:6j6}) the big-oh terms are
uniform for $x \in (-Rn^\gamma, Rn^\gamma)$. We will use
(\ref{eqn:6j2})--(\ref{eqn:6j6}) to verify assumption (\ref{eqn:
2assumeLem51}) for sequences 3a, 4a, 5a. \iffalse We also will use
(\ref{eqn:6j2})--(\ref{eqn:6j4}) to verify assumption (\ref{eqn:
1assumeLem51}) for sequence 6. \fi
%%%------------ pre deg(g)=6  done---------------%%

\skp

%%------------- sequence 1 -----------------%%

\noindent {\bf Sequence 1a}

\noindent This sequence is defined in (\ref{eqn:seq1}). For sequence
1a, $g$ has degree 4. Since $K(\beta_n)-K_n =
(K'(\beta)b-k)/n^\alpha + \mbox{O}(1/n^{2\alpha})$, it follows from
(\ref{eqn:4j2}) that for all $n \in \mathbb{N}$, any
$u>0$, any $\gamma
>0$, any $R>0$, and all $x \in \mathbb{R}$ satisfying
$|x/n^\gamma|<R$, we have
\bea
\label{eqn:1-4j2star}
n^{1-u-2\gamma} G_{\beta_n,K_n}^{(2)}(x/n^\gamma)
&=&
\frac{1}{n^{2\gamma-1+u+\alpha}} 2\beta
(K'(\beta)b-k)(1+\varepsilon_n) \\
&&+
\frac{1}{n^{4\gamma-1+u}} \cdot
\frac{c_4(\beta)(1+\varepsilon_n)\cdot 4!}{2!} x^2 \nonumber \\
&&+
\ \mbox{O} \! \left(\frac{1}{n^{2\gamma-1+u+2\alpha}} \right)
+
\ \mbox{O} \! \left(\frac{1}{n^{5\gamma-1+u}}\right)x^3. \nonumber \eea
We now
define $\gamma=\theta\alpha$ and $u=1-\alpha/\alpha_0$, and we
recall that $\alpha_0=1/2$, $\theta=1/2$. With these choices of
$\gamma$ and $u$, the powers of $n$ appearing in the first two terms
in (\ref{eqn:1-4j2star}) are 0, and the powers of $n$ appearing in
the last two terms in (\ref{eqn:1-4j2star}) are positive. Letting
$n \goto \infty$ in (\ref{eqn:1-4j2star}), we have uniformly for
$x$ in compact subsets of $\mathbb{R}$
\bea
\lim_{n \goto \infty}
n^{1-u-2\gamma} G_{\beta_n,K_n}^{(2)}(x/n^\gamma)
&=&
\lim_{n \goto
\infty} n^{\alpha/\alpha_0-2\theta\alpha}
G_{\beta_n,K_n}^{(2)}(x/n^\gamma)
\nonumber \\
&=&
2\beta (K'(\beta)b-k) + \frac{c_4(\beta) \cdot 4!}{2!} x^2 =
g^{(2)}(x). \nonumber \eea       %%%%%%%%%%%%%%%%%%%%%%%%%%%%%%%%%%%%%

The same choices of $\gamma$ and $u$ ensure that the powers of $n$
appearing in the first term in (\ref{eqn:4j3}) and (\ref{eqn:4j4})
are 0, and the powers of $n$ appearing in the last term in
(\ref{eqn:4j3}) and (\ref{eqn:4j4}) are positive. Taking $n \goto
\infty$ in (\ref{eqn:4j3}) and (\ref{eqn:4j4}) gives
\bea
\lim_{n
\goto \infty} n^{1-u-3\gamma} G_{\beta_n,K_n}^{(3)}(x/n^\gamma)
&=&
\lim_{n \goto \infty} n^{\alpha/\alpha_0-3\theta\alpha}
G_{\beta_n,K_n}^{(3)}(x/n^\gamma)
\nonumber \\
&=&
c_4(\beta) \cdot 4! \cdot x = g^{(3)}(x) \nonumber .\eea   %%%%%%%%%%
and
\bea
\lim_{n \goto \infty} n^{1-u-4\gamma}
G_{\beta_n,K_n}^{(4)}(x/n^\gamma)
&=&
\lim_{n \goto \infty}
n^{\alpha/\alpha_0-4\theta\alpha} G_{\beta_n,K_n}^{(4)}(x/n^\gamma)
\nonumber \\
&=&
c_4(\beta) \cdot 4!
=
g^{(4)}(x) \nonumber. \eea %%%%%%%%%%%%%%%%%%%%%
uniformly for $x$ in compact subsets of $\mathbb{R}$. Thus sequence
1 satisfies hypothesis (\ref{eqn: 1assumeLem51}) in Lemma
\ref{lem:5.1}, and so the conclusion (\ref{eqn: 1concluLem51}) in
Lemma \ref{lem:5.1} follows for $j=2, 3, 4$. This is the convergence
in hypothesis (iii$^\prime$) of Theorem \ref{thm:main}.
%%-------------  1 end  --------------------%%

\skp

%%------------- sequence 2 -----------------%%
\noindent {\bf Sequence 2a}

\noindent This sequence is defined in (\ref{eqn:seq2}). For sequence
2a, $g$ has degree 4. Since $K(\beta_n)-K_n =
(K^{(p)}(\beta)-\ell)b^p/p!n^{p\alpha} +
\mbox{O}(1/n^{\alpha(p+1)})$, it follows from 
(\ref{eqn:4j2}) that for all $n \in \mathbb{N}$, any $u>0$, any
$\gamma >0$, any $R>0$, and all $x \in \mathbb{R}$ satisfying
$|x/n^\gamma|<R$, we have 
\bea 
n^{1-u-2\gamma} G_{\beta_n,K_n}^{(2)}(x/n^\gamma) &=&
\frac{1}{n^{2\gamma-1+u+\alpha}} \cdot \frac{1}{p!} \cdot 2\beta
(K^{(p)}(\beta)-\ell)b^p(1+\varepsilon_n) \label{eqn:2-4j2star} \\
&&+
\frac{1}{n^{4\gamma-1+u}} \cdot
\frac{c_4(\beta)(1+\varepsilon_n)\cdot 4!}{2!} x^2 \nonumber  \\
&&+
\ \mbox{O}\left(\frac{1}{n^{2\gamma-1+u+(p+1)\alpha}}\right)
+
\
\mbox{O}\left(\frac{1}{n^{5\gamma-1+u}}\right)x^3. \nonumber
\eea
We now define
$\gamma=\theta\alpha$ and $u=1-\alpha/\alpha_0$, and we recall that
$\alpha_0=1/2p$, $\theta=p/2$. With these choices of $\gamma$ and
$u$, the power of $n$ appearing in the first two terms in 
(\ref{eqn:2-4j2star}) are 0, and the power of $n$ appearing in the last two
terms in (\ref{eqn:2-4j2star}) are positive. Letting $n \goto
\infty$ in (\ref{eqn:2-4j2star}), we have uniformly for $x$ in
compact subsets of $\mathbb{R}$
\bea
\lim_{n \goto \infty}
n^{1-u-2\gamma} G_{\beta_n,K_n}^{(2)}(x/n^\gamma)
&=&
\lim_{n \goto
\infty} n^{\alpha/\alpha_0-2\theta\alpha}
G_{\beta_n,K_n}^{(2)}(x/n^\gamma) \nonumber
\\
&=&
\frac{1}{p!} \cdot 2\beta (K^{(p)}(\beta)-\ell)b^p +
\frac{c_4(\beta) \cdot 4!}{2!}
x^2 = g^{(2)}(x). \nonumber \eea                %%%%%%%%%%%%%%%%%%%%%%%%%%%%%

The same choices of $\gamma$ and $u$ ensure that the powers of $n$
appearing in the first term in (\ref{eqn:4j3}) and (\ref{eqn:4j4})
are 0, and the powers of $n$ appearing in the last term in
(\ref{eqn:4j3}) and (\ref{eqn:4j4}) are positive. Taking $n \goto
\infty$ in (\ref{eqn:4j3}) and (\ref{eqn:4j4}) gives
\bea
\lim_{n
\goto \infty} n^{1-u-3\gamma} G_{\beta_n,K_n}^{(3)}(x/n^\gamma)
&=&
\lim_{n \goto \infty} n^{\alpha/\alpha_0-3\theta\alpha}
G_{\beta_n,K_n}^{(3)}(x/n^\gamma) \nonumber
\\
&=&
c_4(\beta) \cdot 4! \cdot x = g^{(3)}(x) \nonumber \eea  %%%%%%%%%%%%%
and
\bea
\lim_{n \goto \infty} n^{1-u-4\gamma}
G_{\beta_n,K_n}^{(4)}(x/n^\gamma)
&=&
\lim_{n \goto \infty}
n^{\alpha/\alpha_0-4\theta\alpha} G_{\beta_n,K_n}^{(4)}(x/n^\gamma)
\nonumber
\\
&=&
c_4(\beta) \cdot 4!  = g^{(4)}(x). \nonumber \eea %%%%%%%%%%%%%%%%
uniformly for $x$ in compact subsets of $\mathbb{R}$. Thus sequence
2 satisfies hypothesis (\ref{eqn: 1assumeLem51}) in Lemma
\ref{lem:5.1}, and so the conclusion (\ref{eqn: 1concluLem51}) in
Lemma \ref{lem:5.1} follows for $j=2, 3, 4$. This is the convergence
in hypothesis (iii$^\prime$) of Theorem \ref{thm:main}.
%%-------------  2 end  --------------------%%

\skp

%%------------- sequence 3 -----------------%%
\noindent {\bf Sequence 3a}

\noindent This sequence is defined in (\ref{eqn:seq3}). For sequence
3a, $g$ has degree 6. Since $K(\beta_n)-K_n =
(K'(\beta_c)b-k)/n^\alpha + \mbox{O}(1/n^{2\alpha})$, it follows
(\ref{eqn:6j2}), (\ref{eqn:6j3}), and
(\ref{eqn:6j4}) that for all $n \in \mathbb{N}$, any $u>0$, any
$\gamma >0$, any $R>0$, and all $x \in \mathbb{R}$ satisfying
$|x/n^\gamma|<R$, we have the following: 
\bea \label{eqn:3-6j2star}
n^{1-u-2\gamma} G_{\beta_n,K_n}^{(2)}(x/n^\gamma) &=&
\frac{1}{n^{2\gamma-1+u+\alpha}} \cdot 2\beta_c
(K'(\beta_c)b-k)(1+\varepsilon_n)  \\
&&+
\frac{1}{n^{4\gamma-1+u+\alpha}} \cdot
\frac{c_4(-4 b)(1+\varepsilon_n)\cdot 4!}{2!} x^2 \nonumber \\
&&+
\frac{1}{n^{6\gamma-1+u}} \cdot \frac{c_6(1+\varepsilon_n) \cdot 6!}{4!} \cdot x^4
+ \ \mbox{O}\left(\frac{1}{n^{2\gamma-1+u+2\alpha}}\right) \nonumber \\
&&+
\ \mbox{O}\left(\frac{1}{n^{4\gamma-1+u+2\alpha}}\right)x^2
+ \
\mbox{O}\left(\frac{1}{n^{7\gamma-1+u}}\right)x^5, \nonumber \eea   %%%%%%%%%%%%%%%%%%%%%%%%%%%%%%%%%
\bea
\label{eqn:3-6j3star}
n^{1-u-3\gamma}
G_{\beta_n,K_n}^{(3)}(x/n^\gamma)
&=&
\frac{1}{n^{4\gamma-1+u+\alpha}} \cdot
c_4(-4 b)(1+\varepsilon_n)\cdot 4! \cdot x  \\
&&+
\frac{1}{n^{6\gamma-1+u}} \cdot \frac{c_6(1+\varepsilon_n) \cdot 6!}{3!} \cdot x^3 \nonumber \\
&&+
\ \mbox{O}\left(\frac{1}{n^{4\gamma-1+u+2\alpha}}\right)x
+
\
\mbox{O}\left(\frac{1}{n^{7\gamma-1+u}}\right)x^4 \nonumber, \eea   %%%%%%%%%%%%%%%%%%%%%%%%%%%%%%%%%%%
and
\bea
\label{eqn:3-6j4star}
n^{1-u-4\gamma}
G_{\beta_n,K_n}^{(4)}(x/n^\gamma)
&=&
\frac{1}{n^{4\gamma-1+u+\alpha}} \cdot
c_4(-4 b)(1+\varepsilon_n)\cdot 4! \\
&&+
\frac{1}{n^{6\gamma-1+u}} \cdot \frac{c_6(1+\varepsilon_n) \cdot 6!}{2!} \cdot x^2 \nonumber  \\
&&+
\ \mbox{O}\left(\frac{1}{n^{4\gamma-1+u+2\alpha}}\right)x
+
\
\mbox{O}\left(\frac{1}{n^{7\gamma-1+u}}\right)x^3.\nonumber \eea  %%%%%%%%%%%%%%%%%%%%%%%%%%%%%%%%%
We now define $\gamma=\theta\alpha$ and $u=1-\alpha/\alpha_0$, and
we recall that $\alpha_0=2/3$, $\theta=1/4$. With these choices of
$\gamma$ and $u$, the powers of $n$ appearing in the first term and
the third term in (\ref{eqn:3-6j2star}) are 0, and the powers of
$n$ appearing in the second term and the last three terms in
(\ref{eqn:3-6j2star}) are positive. Letting $n \goto \infty$ in
(\ref{eqn:3-6j2star}), we have uniformly for $x$ in compact subsets
of $\mathbb{R}$
\bea
\lim_{n \goto \infty} n^{1-u-2\gamma}
G_{\beta_n,K_n}^{(2)}(x/n^\gamma)
&=&
\lim_{n \goto \infty}
n^{\alpha/\alpha_0-2\theta\alpha} G_{\beta_n,K_n}^{(2)}(x/n^\gamma) \nonumber \\
&=&
2\beta_c (K'(\beta_c)b-k) + \frac{c_6 \cdot 6!}{4!} x^4 =
g^{(2)}(x). \nonumber \eea      %%%%%%%%%%%%%%%%%%%%%%%%%%%%%

The same choices of $\gamma$ and $u$ ensure that the powers of $n$
appearing in the second term in (\ref{eqn:3-6j3star}) and
(\ref{eqn:3-6j4star}) are 0, and the powers of $n$ appearing in the
first term and last two terms in (\ref{eqn:3-6j3star}) and
(\ref{eqn:3-6j4star}) are positive. Taking $n \goto \infty$ in
(\ref{eqn:3-6j3star}) and (\ref{eqn:3-6j4star}) gives that
\bea
\lim_{n \goto \infty} n^{1-u-3\gamma}
G_{\beta_n,K_n}^{(3)}(x/n^\gamma)
&=&
\lim_{n \goto \infty}
n^{\alpha/\alpha_0-3\theta\alpha}G_{\beta_n,K_n}^{(3)}(x/n^\gamma) \nonumber \\
&=&
\frac{c_6 \cdot 6!}{3!} x^3 = g^{(3)}(x) \nonumber \eea %%%%%%%%%%%%%%%%%%%
and
\bea
\lim_{n \goto \infty} n^{1-u-4\gamma}
G_{\beta_n,K_n}^{(4)}(x/n^\gamma)
&=&
\lim_{n \goto \infty}
n^{\alpha/\alpha_0 - 4\theta\alpha}
G_{\beta_n,K_n}^{(4)}(x/n^\gamma) \nonumber \\
&=&
\frac{c_6 \cdot 6!}{2!} x^2  = g^{(4)}(x). \nonumber \eea %%%%%%%%%%%%%%%%%%%%%%%%%
uniformly for $x$ in compact subsets of $\mathbb{R}$. Similarly,
the powers of $n$ appearing in the first term in the
expansions (\ref{eqn:6j5}) and (\ref{eqn:6j6}) are 0, and the powers
of $n$ appearing in the last term in the expansions (\ref{eqn:6j5})
and (\ref{eqn:6j6}) are positive. Letting $n \goto \infty$ in
(\ref{eqn:6j5}) and (\ref{eqn:6j6}), we have uniformly for $x$ in
compact subsets of $\mathbb{R}$
\bea
\lim_{n \goto \infty}
n^{1-u-5\gamma} G_{\beta_n,K_n}^{(5)}(x/n^\gamma)
&=&
\lim_{n \goto
\infty} n^{\alpha/\alpha_0-5\theta\alpha}
G_{\beta_n,K_n}^{(5)}(x/n^\gamma) \nonumber \\
&=&
c_6 \cdot 6!
x=g^{(5)}(x) \nonumber
\eea %%%%%%%%%%%%%%%%
and
\bea
\lim_{n \goto \infty} n^{1-u-6\gamma}
G_{\beta_n,K_n}^{(6)}(x/n^\gamma)
&=&
\lim_{n \goto \infty}
n^{\alpha/\alpha_0-6\theta\alpha}
G_{\beta_n,K_n}^{(6)}(x/n^\gamma) \nonumber \\
&=&
c_6 \cdot 6! = g^{(6)}(x).
\nonumber \eea %%%%%%%%%%%%%%%%%%%%%%
Thus sequence 3 satisfies hypothesis (\ref{eqn: 2assumeLem51}) in
Lemma \ref{lem:5.1}, and so the conclusion (\ref{eqn: 2concluLem51})
in Lemma \ref{lem:5.1} follows for $j=2, 3, 4, 5, 6$. This is the
convergence in hypothesis (iii$^\prime$) of Theorem \ref{thm:main}.
%%-------------  3 end  --------------------%%

\skp

%%------------- sequence 4 -----------------%%
\noindent {\bf Sequence 4a}

\noindent This sequence is defined in (\ref{eqn:seq4}). For sequence
4a, $g$ has degree 6. Since \bea K(\beta_n)-K_n &=&
K(\beta_c + 1/n^\alpha)-K_n \nonumber \\
&=&
K(\beta_c) + K'(\beta_c)\cdot 1/n^\alpha + K''(\beta_c) \cdot
1/2! n^{2\alpha} + K'''(\beta_c) \cdot 1/3! n^{3\alpha}
\nonumber \\
&&+
\mbox{O}
(1/n^{4\alpha}) - K_n \nonumber \\
&=&
(K''(\beta_c)-\ell)/2n^{2\alpha} + (K'''(\beta_c) -
\tilde{\ell})/6n^{3\alpha} + \mbox{O}(1/n^{4\alpha})\nonumber \eea
and
\[ 4- e^{\beta_n}= -4 /n^\alpha + \mbox{O}(1/n^{2\alpha}), \]
it follows from (\ref{eqn:6j2}), (\ref{eqn:6j3}), and
(\ref{eqn:6j4}) that for all $n \in \mathbb{N}$, any $u>0$, any
$\gamma >0$, any $R>0$, and all $x \in \mathbb{R}$ satisfying
$|x/n^\gamma|<R$, we have the following:
\bea
\label{eqn:4-6j2star}
n^{1-u-2\gamma} G_{\beta_n,K_n}^{(2)}(x/n^\gamma)
&=&
\frac{1}{n^{2\gamma-1+u+2\alpha}} \cdot 2\beta_c \cdot
(K''(\beta_c)-\ell)/2 \cdot (1+\varepsilon_n)  \\
&&+
\frac{1}{n^{4\gamma-1+u+\alpha}} \cdot
\frac{c_4(-4 )(1+\varepsilon_n)\cdot 4!}{2!} x^2  \nonumber \\
&&+
\frac{1}{n^{6\gamma-1+u}} \cdot \frac{c_6(1+\varepsilon_n) \cdot 6!}{4!} \cdot x^4  \nonumber \\
&&+
\ \mbox{O} \! \left(\frac{1}{n^{2\gamma-1+u+3\alpha}}\right)
+
\mbox{O} \! \left(\frac{1}{n^{2\gamma-1+u+4\alpha}}\right) \nonumber \\
&&+
\ \mbox{O} \! \left(\frac{1}{n^{4\gamma-1+u+2\alpha}}\right)x^2
+
\
\mbox{O} \! \left(\frac{1}{n^{7\gamma-1+u}}\right)x^5, \nonumber \eea %%%%%%%%%%%
\bea
\label{eqn:4-6j3star}
n^{1-u-3\gamma}
G_{\beta_n,K_n}^{(3)}(x/n^\gamma)
&=&
\frac{1}{n^{4\gamma-1+u+\alpha}} \cdot
c_4(-4)(1+\varepsilon_n)\cdot 4! \cdot x  \\
&&+
\frac{1}{n^{6\gamma-1+u}} \cdot \frac{c_6(1+\varepsilon_n) \cdot 6!}{3!} \cdot x^3 \nonumber \\
&&+
\mbox{O} \! \left(\frac{1}{n^{4\gamma-1+u+2\alpha}}\right)x
+
\
\mbox{O} \! \left(\frac{1}{n^{7\gamma-1+u}}\right)x^4  \nonumber, \eea %%%%%%%%%%
and
\bea
\label{eqn:4-6j4star}
n^{1-u-4\gamma}
G_{\beta_n,K_n}^{(4)}(x/n^\gamma)
&=&
\frac{1}{n^{4\gamma-1+u+\alpha}} \cdot
c_4(-4)(1+\varepsilon_n)\cdot 4! \\
&&+
\frac{1}{n^{6\gamma-1+u}} \cdot \frac{c_6(1+\varepsilon_n) \cdot 6!}{2!} \cdot x^2 \nonumber  \\
&&+
\ \mbox{O} \! \left(\frac{1}{n^{4\gamma-1+u+2\alpha}}\right)
+
\mbox{O} \! \left(\frac{1}{n^{7\gamma-1+u}}\right)x^3. \nonumber \eea  %%%%%%%%%
We now define $\gamma=\theta\alpha$ and $u=1-\alpha/\alpha_0$, and
we recall that $\alpha_0=1/3$, $\theta=1/2$. With these choices of
$\gamma$ and $u$, the powers of $n$ appearing in the first three
terms in (\ref{eqn:4-6j2star}) are 0, and the powers of $n$
appearing in the last four terms in (\ref{eqn:4-6j2star}) are
positive. Letting $n \goto \infty$ in (\ref{eqn:4-6j2star}), we
have uniformly for $x$ in compact subsets of $\mathbb{R}$
\bea
\lim_{n \goto \infty} n^{1-u-2\gamma}
G_{\beta_n,K_n}^{(2)}(x/n^\gamma)
&=&
\lim_{n \goto \infty}
n^{\alpha/\alpha_0-2\theta\alpha}
G_{\beta_n,K_n}^{(2)}(x/n^\gamma)  \nonumber \\
&=&
2\beta_c (K''(\beta_c)-\ell)/2
+
\frac{c_4 (-4) \cdot 4!}{2!}
x^2 + \frac{c_6 \cdot 6!}{4!} x^4
=
g^{(2)}(x). \nonumber \eea  %%%%%%%%%%%%%%%%%%%%%%

The same choices of $\gamma$ and $u$ ensure that the powers of $n$
appearing in the first two terms in (\ref{eqn:4-6j3star}) and
(\ref{eqn:4-6j4star}) are 0 and the powers of $n$ appearing in the
last two terms in (\ref{eqn:4-6j3star}) and (\ref{eqn:4-6j4star})
are positive. Taking $n \goto \infty$ in (\ref{eqn:4-6j3star}) and
(\ref{eqn:4-6j4star}) gives
\bea
\lim_{n \goto \infty}
n^{1-u-3\gamma} G_{\beta_n,K_n}^{(3)}(x/n^\gamma)
&=&
\lim_{n \goto
\infty}
n^{\alpha/\alpha_0-3\theta\alpha} G_{\beta_n,K_n}^{(3)}(x/n^\gamma) \nonumber \\
&=&
c_4(-4) \cdot 4! x + \frac{c_6 \cdot 6!}{3!} x^3 = g^{(3)}(x)
\nonumber \eea  %%%%%%%%%%%%%%%%%%%%%%
and
\bea
\lim_{n \goto \infty} n^{1-u-4\gamma}
G_{\beta_n,K_n}^{(4)}(x/n^\gamma)
&=&
\lim_{n \goto \infty}
n^{\alpha/\alpha_0-4\theta\alpha}
G_{\beta_n,K_n}^{(4)}(x/n^\gamma) \nonumber \\
&=&
c_4(-4) \cdot 4! + \frac{c_6 \cdot 6!}{2!} x^2
=
g^{(4)}(x). \nonumber \eea  %%%%%%%%%%%%%%%%%%%%%%
uniformly for $x$ in compact subsets of $\mathbb{R}$. Similarly,
the powers of $n$ appearing in the first term in the
expansions (\ref{eqn:6j5}) and (\ref{eqn:6j6}) are 0 and the powers
of $n$ appearing in the last term in the expansions (\ref{eqn:6j5})
and (\ref{eqn:6j6}) are positive. Letting $n \goto \infty$ in
(\ref{eqn:6j5}) and (\ref{eqn:6j6}), we have uniformly for $x$ in
compact subsets of $\mathbb{R}$
\bea
\lim_{n \goto \infty}
n^{1-u-5\gamma} G_{\beta_n,K_n}^{(5)}(x/n^\gamma)
&=&
\lim_{n \goto
\infty} n^{\alpha/\alpha_0-5\theta\alpha}
G_{\beta_n,K_n}^{(5)}(x/n^\gamma) \nonumber \\
&=&
c_6 \cdot 6! x
=
g^{(5)}(x) \nonumber \eea  %%%%%%%%%%%%%%%%%%%%%%
and
\bea
\lim_{n \goto \infty} n^{1-u-6\gamma}
G_{\beta_n,K_n}^{(6)}(x/n^\gamma)
&=&
\lim_{n \goto \infty}
n^{\alpha/\alpha_0-6\theta\alpha}
G_{\beta_n,K_n}^{(6)}(x/n^\gamma) \nonumber \\
&=&
c_6 \cdot 6!
=
g^{(6)}(x).  \nonumber \eea  %%%%%%%%%%%%%%%%%%%%%%
Thus sequence 4 satisfies hypothesis (\ref{eqn: 2assumeLem51}) in
Lemma \ref{lem:5.1}, and so the conclusion (\ref{eqn: 2concluLem51})
in Lemma \ref{lem:5.1} follows for $j=2, 3, 4, 5, 6$. This is the
convergence in hypothesis (iii$^\prime$) of Theorem \ref{thm:main}.
%%-------------  4 end  --------------------%%

\skp

%%------------- sequence 5 -----------------%%
\noindent {\bf Sequence 5a}

\noindent This sequence is defined in (\ref{eqn:seq5}). For sequence
5a, $g$ has degree 6. Since $
 K(\beta_n)-K_n = (K''(\beta_c)-\ell)/2n^{2\alpha} + \mbox{O}(1/n^{3\alpha})$
and $4- e^{\beta_n}= 4 /n^\alpha + \mbox{O}(1/n^{2\alpha})$, it
follows from (\ref{eqn:6j2}), (\ref{eqn:6j3}), and
(\ref{eqn:6j4}) that for all $n \in \mathbb{N}$, any $u>0$, any
$\gamma >0$, any $R>0$, and all $x \in \mathbb{R}$ satisfying
$|x/n^\gamma|<R$, we have the following:
\bea
\label{eqn:5-6j2star}
n^{1-u-2\gamma} G_{\beta_n,K_n}^{(2)}(x/n^\gamma)
&=&
\frac{1}{n^{2\gamma-1+u+2\alpha}} \cdot 2\beta_c \cdot
(K''(\beta_c)-\ell)/2 \cdot (1+\varepsilon_n)  \\
&&+
\frac{1}{n^{4\gamma-1+u+\alpha}} \cdot
\frac{c_4 \cdot 4 \cdot (1+\varepsilon_n)\cdot 4!}{2!} x^2 \nonumber \\
&&+
\frac{1}{n^{6\gamma-1+u}} \cdot \frac{c_6(1+\varepsilon_n) \cdot 6!}{4!} \cdot x^4  \nonumber \\
&&+
\ \mbox{O} \! \left(\frac{1}{n^{2\gamma-1+u+3\alpha}}\right)  \nonumber \\
&&+
\ \mbox{O} \! \left(\frac{1}{n^{4\gamma-1+u+2\alpha}}\right)x^2
+ \
\mbox{O} \! \left(\frac{1}{n^{7\gamma-1+u}}\right)x^5, \nonumber \eea  %%%%%%%%%%%%%%%%%%%%%
\bea
\label{eqn:5-6j3star}
n^{1-u-3\gamma}
G_{\beta_n,K_n}^{(3)}(x/n^\gamma)
&=&
\frac{1}{n^{4\gamma-1+u+\alpha}} \cdot
c_4 \cdot 4 \cdot (1+\varepsilon_n)\cdot 4! \cdot x  \\
&&+
\frac{1}{n^{6\gamma-1+u}} \cdot \frac{c_6(1+\varepsilon_n) \cdot 6!}{3!} \cdot x^3 \nonumber \\
&&+
\ \mbox{O} \! \left(\frac{1}{n^{4\gamma-1+u+2\alpha}}\right)x
+
\
\mbox{O} \! \left(\frac{1}{n^{7\gamma-1+u}}\right)x^4  \nonumber, \eea  %%%%%%%%%%%%%%%%%%%%%%%
and
\bea
\label{eqn:5-6j4star}
n^{1-u-4\gamma}
G_{\beta_n,K_n}^{(4)}(x/n^\gamma)
&=&
\frac{1}{n^{4\gamma-1+u+\alpha}} \cdot
c_4 \cdot 4 \cdot (1+\varepsilon_n)\cdot 4!  \\
&&+
\frac{1}{n^{6\gamma-1+u}} \cdot \frac{c_6(1+\varepsilon_n) \cdot 6!}{2!} \cdot x^2 \nonumber \\
&&+
\ \mbox{O} \! \left(\frac{1}{n^{4\gamma-1+u+2\alpha}}\right)
+ \
\mbox{O} \! \left(\frac{1}{n^{7\gamma-1+u}}\right)x^3.\nonumber \eea %%%%%%%%%%%%%%%%%%%%
We now define $\gamma=\theta\alpha$ and $u=1-\alpha/\alpha_0$, and
we recall that $\alpha_0=1/3$, $\theta=1/2$. With these choices of
$\gamma$ and $u$, the powers of $n$ appearing in the first three
terms in (\ref{eqn:5-6j2star}) are 0, and the powers of $n$
appearing in the last three terms in (\ref{eqn:5-6j2star}) are
positive. Letting $n \goto \infty$ in (\ref{eqn:5-6j2star}), we
have uniformly for $x$ in compact subsets of $\mathbb{R}$
\bea
\lim_{n \goto \infty}
n^{1-u-2\gamma}
G_{\beta_n,K_n}^{(2)}(x/n^\gamma)
&=&
\lim_{n \goto \infty}
n^{\alpha/\alpha_0-2\theta\alpha}
G_{\beta_n,K_n}^{(2)}(x/n^\gamma) \nonumber \\
&=&
\beta_c (K''(\beta_c)-\ell)/2 + \frac{c_4 \cdot 4 \cdot 4!}{2!}
x^2
+
\frac{c_6 \cdot 6!}{4!} x^4 = g^{(2)}(x). \nonumber \eea %%%%%%%%%%%%%%%%%%%%%%

The same choices of $\gamma$ and $u$ ensure that the powers of $n$
appearing in the first two terms in (\ref{eqn:5-6j3star}) and
(\ref{eqn:5-6j4star}) are 0 and the powers of $n$ appearing in the
last two terms in (\ref{eqn:5-6j3star}) and (\ref{eqn:5-6j4star})
are positive. Taking $n \goto \infty$ in (\ref{eqn:5-6j3star}) and
(\ref{eqn:5-6j4star}) gives
\bea
\lim_{n \goto \infty}
n^{1-u-3\gamma} G_{\beta_n,K_n}^{(3)}(x/n^\gamma)
&=&
\lim_{n \goto
\infty}
n^{\alpha/\alpha_0-3\theta\alpha} G_{\beta_n,K_n}^{(3)}(x/n^\gamma) \nonumber \\
&=&
c_4 \cdot 4 \cdot 4! x + \frac{c_6 \cdot 6!}{3!} x^3
=
g^{(3)}(x)
\nonumber \eea %%%%%%%%%%%%%%%%%%%%%%
and
\bea
\lim_{n \goto \infty} n^{1-u-4\gamma}
G_{\beta_n,K_n}^{(4)}(x/n^\gamma)
&=&
\lim_{n \goto \infty} n^{\alpha/\alpha_0-4\theta\alpha} G_{\beta_n,K_n}^{(4)}(x/n^\gamma) \nonumber \\
&=&
c_4 \cdot 4 \cdot 4! + \frac{c_6 \cdot 6!}{2!} x^2
=
g^{(4)}(x). \nonumber \eea %%%%%%%%%%%%%%%%%%%%%%
uniformly for $x$ in compact subsets of $\mathbb{R}$. Similarly,
the powers of $n$ appearing in the first term in the
expansions (\ref{eqn:6j5}) and (\ref{eqn:6j6}) are 0, and the powers
of $n$ appearing in the last term in the expansions (\ref{eqn:6j5})
and (\ref{eqn:6j6}) are positive. Letting $n \goto \infty$ in
(\ref{eqn:6j5}) and (\ref{eqn:6j6}), we have uniformly for $x$ in
compact subsets of $\mathbb{R}$
\bea
\lim_{n \goto \infty}
n^{1-u-5\gamma} G_{\beta_n,K_n}^{(5)}(x/n^\gamma)
&=&
\lim_{n \goto \infty} n^{\alpha/\alpha_0-5\theta\alpha} G_{\beta_n,K_n}^{(5)}(x/n^\gamma)  \nonumber \\
&=&
c_6 \cdot 6! x = g^{(5)}(x)  \nonumber \eea %%%%%%%%%%%%%%%%%%%%%%
and
\bea
\lim_{n \goto \infty}
n^{1-u-6\gamma}G_{\beta_n,K_n}^{(6)}(x/n^\gamma)
&=&
\lim_{n \goto \infty} n^{\alpha/\alpha_0-6\theta\alpha}G_{\beta_n,K_n}^{(6)}(x/n^\gamma) \nonumber \\
&=&
c_6 \cdot 6! = g^{(6)}(x). \nonumber \eea %%%%%%%%%%%%%%%%%%%%%%
Thus sequence 5 satisfies hypothesis (\ref{eqn: 2assumeLem51}) in
Lemma \ref{lem:5.1}, and so the conclusion (\ref{eqn: 2concluLem51})
in Lemma \ref{lem:5.1} follows for $j=2, 3, 4, 5, 6$. This is the
convergence in hypothesis (iii$^\prime$) of Theorem \ref{thm:main}.
%%-------------  5 end  --------------------%%

\skp \skp

\section{Proof of the MDP in Part (a) of Theorem \ref{thm:MDP}}
\label{section:appendixB}
\setcounter{equation}{0}

In this appendix we give the proof of part (a) of the MDP stated in
Theorem \ref{thm:MDP}. We restate the theorem here for easy
reference. Concerning the proof of part (b) of Theorem
\ref{thm:MDP}, see the comment before Lemma \ref{lem:Sym1/2Small}.

\begin{thmN}
Let $(\beta_n,K_n)$ be a positive sequence 
converging either to a second-order point $(\beta,
K(\beta))$, $0<\beta<\beta_c$, or to the tricritical point $(\beta,
K(\beta))=(\beta_c, K(\beta_c))$. We assume that $(\beta_n, K_n)$
satisfies the hypotheses of Theorem {\em\ref{thm:main}} for all
$0<\alpha<\alpha_0$. The following conclusions hold.

\begin{itemize}
 % \vspace{-.1in}
  \item[\em (a)] For all $0<\alpha<\alpha_0$,
  $S_n/n^{1-\theta\alpha}$ satisfies the MDP with respect to $P_{n,
  \beta_n,K_n}$ with exponential speed $n^{1-\alpha/\alpha_0}$ and
  rate function $\Gamma(x)=g(x)-\inf_{y \in \mathbb{R}}g(y)$; in symbols
  \[ P_{n, \beta_n, K_n} \{S_n/n^{1-\theta\alpha}\in dx\} \asymp
  \exp[-n^{1-\alpha/\alpha_0} \Gamma(x)]dx.
  \]

  \item[\em (b)] The hypotheses of this theorem are satisfied by
  sequence {\em{1a--5a}}
  defined in Table {\em 5.1} .
\end{itemize}
\end{thmN}

We work with an arbitrary $\alpha$ satisfying $0<\alpha<\alpha_0$.
To ease the notation we write $\gamma=\theta\alpha$ and
$u=1-\alpha/\alpha_0$. The hypotheses of Theorem \ref{thm:MDP}
coincide with the hypotheses of Theorem \ref{thm:main}, which in
turn consist of hypothesis (iii$^\prime$) and the hypotheses of
Theorem \ref{thm:3.1} for all $0<\alpha<\alpha_0$. Clearly we have
$0<u=1-\alpha/\alpha_0<1$ and by hypothesis (iii$^\prime$)
$0<\gamma=\theta\alpha<\theta\alpha_0 <1/2$. In addition,
$1-2\gamma-u=(1-2\theta\alpha_0)\alpha/\alpha_0>0$, which implies
$1-2\gamma>u$.

\skp

The proof of Theorem \ref{thm:MDP} is analogous to the proof of
Theorem 8.1 in \cite{CosEllOtt}. Let $W_n$ be a sequence of normal
random variables with mean 0 and variance $\sigma_n^2=(2\beta_n
K_n)^{-1}$ defined on a probability space $(\Omega,\mathcal {F},Q)$.
Theorem \ref{thm:MDP} is proved in two steps.

\skp \noi {\bf Step 1.} $W_n/n^{1/2-\gamma}$ is superexponentially
small relative to $\exp (n^{-v})$; i.e., for any $\delta>0$ \be
\label{eqn:AppendixBsmallWn} \limsup_{n \goto \infty}
\frac{1}{n^{-v}} \log Q \{ |W_n/n^{1/2-\gamma}|>\delta \} = - \infty
\ee

\skp \noi {\bf Step 2.} With respect to $P_{n,\beta_n, K_n} \times
Q$, $S_n/n^{1-\gamma}+W_n/n^{1/2-\gamma}$ satisfies the Laplace
principle with exponential speed $n^{-v}$ and rate function
$\Gamma$.

\skp According to Theorem 1.3.3 in \cite{DupEll}, if we prove Step 1
and Step 2, then with respect to $P_{n,\beta_n,K_n}$,
$S_n/n^{1-\gamma}$ satisfies the Laplace principle with speed $n^u$
and rate function $\Gamma$; i.e., for any bounded, continuous
function $\psi$
\[ \lim_{n \goto \infty} \frac{1}{n^u} \log \int_{\Lambda^n} \exp
[n^u \psi(S_n/n^{1-\gamma})] d P_{n,\beta_n,K_n} = \sup_{x \in
\mathbb{R}} \{ \psi(x)- \Gamma(x) \}.
\]
Since the Laplace principle implies the MDP (Thm 1.2.3 in
\cite{DupEll}), Theorem \ref{thm:MDP} follows.

\skp Next, we prove Step 1 and Step 2.

\skp {\noindent}{\bf Proof of Step 1.} Since $\beta_n$ and $K_n$ are
bounded and uniformly positive over $n$, the sequence $\sigma_n^2$
is bounded and uniformly positive over $n$. We have \bea
Q\{|W_n/n^{1/2-\gamma}|>\delta\} &=& Q \{ | N(0, \sigma_n^2) | >
n^{1/2-\gamma} \delta \} \nonumber \\
&\le&
\frac{\sqrt{2} \sigma_n}{\sqrt{\pi} n^{1/2-\gamma}\delta}
\cdot \exp(-n^{1-2\gamma} \delta^2 / (2\sigma_n^2)) \nonumber. \eea
\[ \frac{1}{n^u} \log Q \{ |W_n/n^{1/2-\gamma}|>\delta \}
\le \frac{1}{n^u} \left[ \log \frac{\sqrt{2} \sigma_n}{\sqrt{\pi}
\delta} + \log (n^{\gamma-1/2})-
\frac{n^{1-2\gamma}\delta^2}{2\sigma_n^2} \right].
\]
The limit of the right hand side of the last inequality is $-\infty$
since $u>0$ and $1-2\gamma
> u$. Thus (\ref{eqn:AppendixBsmallWn}) follows. The proof of Step
1 is done.

\skp {\noindent}{\bf Proof of Step 2.}  Let $\psi$ be an arbitrary
bounded, continuous function. Choosing $\varphi = \exp[n^u \psi]$
and $\bar{\gamma} = \gamma$ in Lemma \ref{lem:RepresentF} yields
\bea \label{eqn:RepresentApply}
\lefteqn{ \int_{\Lambda^n \times \Omega} \exp \left[n^{u}\psi \left(\frac{S_n}{n^{1-\gamma}}+\frac{W_n}{n^{1/2-\gamma}}\right)\right] d(P_{n,\beta_n,K_n}\times Q) } \\
&&= \frac{1}{\int_{\mathbb{R}} \exp[-n
G_{\beta_n,K_n}(x/n^{\gamma})]dx} \cdot \int_{\mathbb{R}}
\exp[n^{u}\psi(x)-n G_{\beta_n,K_n}(x/n^{\gamma})]dx \nonumber. \eea
The proof of Step 2 rests on the following three properties of $n
G_{\beta_n,K_n}(x/n^{\gamma})$.
\begin{itemize}
\item [1.] By hypothesis (iv) of Theorem \ref{thm:3.1} for $0<\alpha<\alpha_0$, there exists a polynomial $H$ satisfying $H(x) \goto
\infty$ as $|x| \goto \infty$ together with the following property:
$\exists R>0$ such that for $\forall n \in \N$ sufficiently large
and $\forall x \in \R$ satisfying $|x/n^{\gamma}|<R$
\[ n G_{\beta_n,K_n}(x/n^{\gamma}) \ge n^{u} H(x).
\]
\item [ 2.] Let $\Delta = \sup_{x \in \mathbb{R}} \{ \psi(x) - g(x)
\}$. Since $H(x) \goto \infty$ and $g(x) \goto \infty$ as $|x| \goto
\infty$, there exists $M >0$ with the following three properties:
\[ \sup_{|x| > M} \{ \psi(x) - H(x) \} \le -|\Delta|-1,
\]
the supremum of $\psi-g$ on $\mathbb{R}$ is attained on the interval
$[-M,M]$, and the supremum of $-g$ on $\mathbb{R}$ is attained on
the interval $[-M,M]$. In combination with item 1, we have that for
all $n \in \mathbb{N}$ satisfying $R n^{\gamma} > M$ \bea
\label{eqn:appBitem2} \lefteqn{ \sup_{M<|x|<Rn^{\gamma}} \{
n^{u}\psi(x) - n
G_{\beta_n,K_n}(x/n^{\gamma}) \} } \\
&& \le \sup_{M<|x|<Rn^{\gamma}} \{ n^{u}\psi(x) - n^{u}H(x) \}
\nonumber \\
&& \le -n^{u} (|\Delta|+1). \nonumber \eea
\item [3.] Let $M$ be the number selected in item 2. By hypothesis (iii)(a) of Theorem \ref{thm:3.1} for $0<\alpha<\alpha_0$, for all
$x \in \mathbb{R}$ satisfying $|x| \le M$, $n^{1-u}
G_{\beta_n,K_n}(x/n^{\gamma})$ converges uniformly to $g(x)$ as $n
\goto \infty$.
\end{itemize}

Item 3 implies that for any $\delta >0$ and all sufficiently large
$n$
\bea \lefteqn{\exp(-n^{u} \delta) \int_{ \{ |x|\le M \} } \exp
[n^{u} ( \psi(x) - g(x) )]dx } \nonumber
\\
&& \le \int_{ \{ |x|\le M \} } \exp
[n^{u}  \psi(x) - n G_{\beta_n,K_n}(x/n^{\gamma}) ]dx \nonumber  \\
&& \le \exp(n^{u} \delta) \int_{ \{ |x|\le M \} } \exp [n^{u} (
\psi(x) - g(x) )]dx \nonumber. \eea
In addition, item 2 implies that
\[ \int_{ \{ M<|x|<Rn^{\gamma} \} } \exp
[n^{u}  \psi(x) - n G_{\beta_n,K_n}(x/n^{\gamma}) ]dx \le
2Rn^{\gamma} \exp[-n^{u}(|\Delta|+1)].
\]
Since $\psi$ is bounded, the last two displays show that there
exists $a_1 > 0$ and $a_2 \in \mathbb{R}$ such that for all
sufficiently large $n$
\[ \int_{ \{ |x|<Rn^{\gamma} \} } \exp
[- n G_{\beta_n,K_n}(x/n^{\gamma}) ]dx
\le
a_1 \exp(n^{u}a_2).
\]
Since $u \in (0,1)$, by part (d) of Lemma 4.4 in \cite{CosEllOtt}
there exists $a_3 >0$ such that for all sufficiently large $n$
\[  \int_{ \{ |x| \ge Rn^{\gamma} \} } \exp
[- n G_{\beta_n,K_n}(x/n^{\gamma}) ]dx \le 2a_1 \exp(-n a_3).
\]

Together these three estimates show that for all sufficiently large
$n$
\bea \lefteqn{ \int_{ \mathbb{R} } \exp
[n^{u}  \psi(x) - n G_{\beta_n,K_n}(x/n^{\gamma}) ]dx } \nonumber \\
&=& \int_{ \{ |x|\le M \} } \exp
[n^{u}  \psi(x) - n G_{\beta_n,K_n}(x/n^{\gamma}) ]dx \nonumber \\
&&+ \int_{ \{ M<|x|<Rn^{\gamma} \} } \exp
[n^{u}  \psi(x) - n G_{\beta_n,K_n}(x/n^{\gamma}) ]dx \nonumber \\
&&+ \int_{ \{ |x|\ge Rn^{\gamma} \} } \exp
[n^{u}  \psi(x) - n G_{\beta_n,K_n}(x/n^{\gamma}) ]dx \nonumber \\
&\le& \exp(n^{u} \delta) \int_{ \{ |x|\le M \} } \exp [n^{u} (
\psi(x) - g(x) )]dx \nonumber \\
&&+  2Rn^{\gamma} \exp[-n^{u}(|\Delta|+1)] + \,  2a_1 \exp(-n
a_3+n^{u} ||\psi||_{\infty}) \nonumber.
\eea
Hence for all
sufficiently large $n$ we have \bea \lefteqn{\exp(-n^{u} \delta)
\int_{ \{ |x|\le M \} } \exp [n^{u} ( \psi(x) - g(x) )]dx }
\nonumber
\\
&& \le \int_{ \mathbb{R}} \exp
[n^{u}  \psi(x) - n G_{\beta_n,K_n}(x/n^{\gamma}) ]dx \nonumber  \\
&& \le \exp(n^{u} \delta) \int_{ \{ |x|\le M \} } \exp [n^{u} (
\psi(x) - g(x) )]dx + \delta_n \nonumber,
\eea
where
\bea
\label{eqn:deltaN} \delta_n &\le& 2Rn^{\gamma}
\exp[-n^{u}(|\Delta|+1)] + 2a_1 \exp(-n a_3+n^{u}
||\psi||_{\infty})  \\
&\le& 4Rn^{\gamma} \exp[-n^{u}(|\Delta|+1)]. \nonumber
\eea
It
follows that \bea \label{eqn:appBliminflimsup} \lefteqn{ \liminf_{n
\goto \infty} \frac{1}{n^{u}} \log \left[ \exp(-n^{u} \delta) \int_{
\{ |x|\le M \} } \exp [n^{u} ( \psi(x) - g(x) )]dx \right] }
\\
&& \le \liminf_{n \goto \infty} \frac{1}{n^{u}} \log \int_{
\mathbb{R}} \exp
[n^{u}  \psi(x) - n G_{\beta_n,K_n}(x/n^{\gamma}) ]dx \nonumber  \\
&& \le \limsup_{n \goto \infty} \frac{1}{n^{u}} \log \int_{
\mathbb{R}} \exp
[n^{u}  \psi(x) - n G_{\beta_n,K_n}(x/n^{\gamma}) ]dx \nonumber  \\
&& \le \limsup_{n \goto \infty} \frac{1}{n^{u}} \log
\left[\exp(n^{u} \delta) \int_{ \{ |x|\le M \} } \exp [n^{u} (
\psi(x) - g(x) )]dx + \delta_n \right]. \nonumber
\eea

By Laplace's method applied to the continuous function $\psi - g$ on
$|x| \le M$ and the fact that the supremum of $\psi - g$ on
$\mathbb{R}$ is attained on the interval $[-M , M]$, we have
\bea
\label{eqn:appBLap} \lefteqn{ \lim_{n \goto \infty} \frac{1}{n^{u}}
\log \int_{ \{ |x|\le M \} } \exp [n^{u} ( \psi(x) - g(x) )]dx } \\
&&= \sup_{|x| \le M} \{\psi(x)-g(x)\} = \sup_{x \in \mathbb{R}} \{
\psi(x) -g(x) \} \nonumber.
\eea
Hence the first line of
(\ref{eqn:appBliminflimsup}) equals
\bea \label{eqn:appBlower}
\lefteqn{ \liminf_{n \goto \infty} \frac{1}{n^{u}} \left[
-n^{u}\delta + \log \int_{ \{ |x|\le M \} }
\exp [n^{u} ( \psi(x) - g(x) )]dx \right] } \\
&&= -\delta + \liminf_{n \goto \infty} \frac{1}{n^{u}} \log \int_{
\{ |x|\le M \} }
\exp [n^{u} ( \psi(x) - g(x) )]dx  \nonumber \\
&&= -\delta + \sup_{x \in \mathbb{R}} \{ \psi(x) -g(x) \} \nonumber.
\eea

We have to work harder to evaluate the last line of
(\ref{eqn:appBliminflimsup}). At the end of the proof we will show
that the term $\delta_n$ can be neglected in evaluating the last
line of (\ref{eqn:appBliminflimsup}); i.e.,
\bea
\label{eqn:appBwithoutDeltaN} \lefteqn{\limsup_{n \goto \infty}
\frac{1}{n^{u}} \log \left[\exp(n^{u} \delta) \int_{ \{ |x|\le M \}
} \exp [n^{u} (\psi(x) - g(x) )]dx + \delta_n \right]}  \\
&&= \limsup_{n \goto \infty} \frac{1}{n^{u}} \log \left[\exp(n^{u}
\delta) \int_{ \{ |x|\le M \} } \exp [n^{u} ( \psi(x) - g(x) )]dx
\right]. \nonumber
\eea
Under the assumption that this is true, by
(\ref{eqn:appBLap}) the last line of (\ref{eqn:appBliminflimsup})
equals
\bea
\label{eqn:appBsupper} \lefteqn{\limsup_{n \goto \infty}
\frac{1}{n^{u}} \log \left[\exp(n^{u} \delta) \int_{ \{ |x|\le M \}
} \exp [n^{u} ( \psi(x) - g(x) )]dx
\right]} \\
&&= \limsup_{n \goto \infty} \frac{1}{n^{u}} \left[ n^{u}\delta +
\log \int_{ \{ |x|\le M \} }
\exp [n^{u} ( \psi(x) - g(x) )]dx \right]  \nonumber \\
&&= \delta + \lim_{n \goto \infty} \frac{1}{n^{u}} \log \int_{ \{
|x|\le M \} }
\exp [n^{u} ( \psi(x) - g(x) )]dx  \nonumber \\
&&= \delta + \sup_{x \in \mathbb{R}} \{ \psi(x) -g(x) \} \nonumber.
\eea
Since $\delta>0$ is arbitrary, combining
(\ref{eqn:appBliminflimsup}), (\ref{eqn:appBlower}), and
(\ref{eqn:appBsupper}) yields
\[ \lim_{n \goto \infty} \frac{1}{n^{-v}} \log \int_{
\mathbb{R}} \exp [n^{-v}  \psi(x) - n G_{\beta_n,K_n}(x/n^{\gamma})
]dx = \sup_{x \in \mathbb{R}} \{ \psi(x) -g(x) \}.
\]

Using the fact that the supremum of $g$ is attained on the interval
$[-M,M]$ (see item 2 in the proof of Step 2), we apply the limit in
the last display to $\psi=0$. We conclude from
(\ref{eqn:RepresentApply}) that \bea
\lefteqn{ \lim_{n \goto \infty} \frac{1}{n^{u}} \log \int_{\Lambda^n \times \Omega} \exp \left[n^{u}\psi \left(\frac{S_n}{n^{1-\gamma}}+\frac{W_n}{n^{1/2-\gamma}}\right)\right] d(P_{n,\beta_n,K_n}\times Q) } \nonumber \\
&=& \lim_{n \goto \infty}  \frac{1}{n^{u}} \log \int_{\mathbb{R}}
\exp[n^{u}\psi(x)-n G_{\beta_n,K_n}(x/n^{\gamma})]dx \nonumber
\\
&&- \lim_{n \goto \infty}  \frac{1}{n^{u}} \log \int_{\mathbb{R}}
\exp[-n G_{\beta_n,K_n}(x/n^{\gamma})]dx \nonumber \\
&=& \sup_{x \in \mathbb{R}} \{ \psi(x) -g(x) \} + \sup_{x \in
\mathbb{R}} \{ -g(x) \} \nonumber \\
&=&  \sup_{x \in \mathbb{R}} \{ \psi(x) -g(x) \} + \inf_{y \in
\mathbb{R}} g(y) \nonumber \\
&=&  \sup_{x \in \mathbb{R}} \{ \psi(x) -\Gamma(x) \} \nonumber.
\eea
Except for the proof of (\ref{eqn:appBwithoutDeltaN}) we have
completed the proof of Step 2, which show that with respect to
$P_{n,\beta_n, K_n} \times Q$, $S_n/n^{1-\gamma}+W_n/n^{1/2-\gamma}$
satisfies the Laplace principle with exponential speed $n^{-v}$ and
rate function $\Gamma$.

To prove (\ref{eqn:appBwithoutDeltaN}) we define
\[ A_n =\exp(n^{u} \delta) \int_{ \{ |x|\le M \} } \exp [n^{u} (
\psi(x) - g(x) )]dx.
\]
It suffices to show that $\delta_n / A_n \goto 0$. To see this, we
rewrite (\ref{eqn:appBwithoutDeltaN}) as follows:
\bea
\lefteqn{\limsup \frac{1}{n^u} \log (A_n + \delta_n) }
\nonumber \\
&&= \limsup \frac{1}{n^u} \log \left[ A_n \left(1+ \frac{\delta_n}{A_n} \right) \right] \nonumber \\
&&= \limsup \left[ \frac{1}{n^u} \log A_n + \frac{1}{n^u} \log \left(1+\frac{\delta_n}{A_n} \right) \right] \nonumber \\
&&= \limsup \frac{1}{n^u} \log A_n \nonumber. \eea
Now we prove that
$\delta_n / A_n  \goto 0$. By (\ref{eqn:appBLap}) we have
\bea \lim
\frac{1}{n^u} \log A_n &=& \delta  + \lim_{n \goto \infty}
\frac{1}{n^{u}} \log \int_{\{|x|\le M \} } \exp [n^{u} ( \psi(x) - g(x) )]dx \nonumber \\
&=& \delta + \sup_{x \in \mathbb{R}} \{ \psi(x) -g(x) \}  \nonumber \\
&=& \delta + \Delta  \nonumber ,
\eea
which implies that for all
sufficiently large $n$
\[ A_n \ge \exp \left[n^u \left( \frac{\delta}{2} +
\Delta \right) \right].
\]
Since by (\ref{eqn:deltaN}) we have for all sufficiently large $n$
\[ \delta_n \le 4Rn^{\gamma} \exp[-n^{u}(|\Delta|+1)],
\]
it follows that for any $0< \varepsilon <1$ and all sufficiently
large $n$
\[ \delta_n \le \exp[n^{u}(-|\Delta|-1 + \varepsilon)]
\]
and thus
\be \label{eqn:appBdeltaNoverAn} 0 \le \frac{\delta_n}{A_n}
\le \exp \left[ n^{u} \left( -|\Delta|-1 + \varepsilon-
\frac{\delta}{2} - \Delta \right) \right].
\ee
If $\Delta \ge 0$,
then
\[ -|\Delta|-1 + \varepsilon- \frac{\delta}{2} - \Delta = -1 + \varepsilon- \frac{\delta}{2} -
2\Delta <0.
\]
If $\Delta <0$,
then
\[ -|\Delta|-1 + \varepsilon- \frac{\delta}{2} - \Delta = -1 + \varepsilon- \frac{\delta}{2}  <0.
\]
Thus in all cases the limit of the right hand side of
(\ref{eqn:appBdeltaNoverAn}) is 0. This completes the proof of
(\ref{eqn:appBwithoutDeltaN}).

\skp Together Step 1 and Step 2 prove that with respect to
$P_{n,\beta_n,K_n}$, $S_n/n^{1-\gamma}$ satisfies the Laplace
principle with speed $n^{-v}$ and rate function $\Gamma$. The proof
of Theorem \ref{thm:MDP} is complete.

%===================== end appendix =====================================

\end{document}